\documentclass[a4paper,twoside]{article}

\newcommand{\heuteIst}{January 8, 2007 }

\usepackage{amsmath,amsthm,amsfonts,latexsym,amscd,amssymb,enumerate}
\usepackage[ backref,hypertex]{hyperref}

\swapnumbers
\theoremstyle{plain}
\newtheorem{theorem}{Theorem}[section]
\newtheorem{lemma}[theorem]{Lemma}
\newtheorem{corollary}[theorem]{Corollary}
\newtheorem{proposition}[theorem]{Proposition}

\theoremstyle{definition}
\newtheorem{definition}[theorem]{Definition}
\newtheorem{example}[theorem]{Example}
\newtheorem{notation}[theorem]{Notation}

\theoremstyle{remark}
\newtheorem{remark}[theorem]{Remark}

\newcommand{\profinQ}{\mathfrak{Q}}

\newcommand{\complexs}{\mathbb{C}}
\newcommand{\naturals}{\mathbb{N}}
\newcommand{\integers}{\mathbb{Z}}
\newcommand{\rationals}{\mathbb{Q}}

\DeclareMathOperator{\id}{id}

\newcommand{\abs}[1]{\left\lvert#1\right\rvert} 

\newcommand{\generate}[1]{\langle#1\rangle}
\newcommand{\tensor}{\otimes}
\newcommand{\into}{\hookrightarrow}
\newcommand{\onto}{\twoheadrightarrow}
\newcommand{\iso}{\cong}

\newcommand{\subgroup}{\leq}
\newcommand{\supergroup}{\geq}
\newcommand{\superset}{\supset}
\newcommand{\normalsubgroup}{\lhd}
\newcommand{\semiProd}{\rtimes}
\newcommand{\semiprod}{\semiProd}
\newcommand{\dividesnot}{\nmid}
\newcommand{\divides}{\mid}
\DeclareMathOperator{\Aut}{Aut}
\DeclareMathOperator{\im}{im}      

\DeclareMathOperator{\tr}{tr}
\DeclareMathOperator{\pr}{pr}

\DeclareMathOperator*{\dirlim}{dirlim}
\DeclareMathOperator*{\invlim}{invlim}
\DeclareMathOperator{\rank}{rank}
\DeclareMathOperator{\coind}{coind}
\DeclareMathOperator{\res}{res}

\newcommand{\forget}[1]{}

\newcommand{\innerprod}[1]{\langle #1 \rangle}

{\catcode`@=11\global\let\c@equation=\c@theorem}

\allowdisplaybreaks[2]

\newcommand{\extendedC}{\mathcal{D}}
\newcommand{\extendableGroups}{\mathcal{F}}
\newcommand{\tfsubgroups}{{\mathcal{A}}}

\DeclareMathOperator{\lcm}{lcm}

\newcommand{\NeumannN}{\mathcal{N}}
\newcommand{\universalU}{\mathcal{U}}
\newcommand{\amenableGroups}{{\mathcal{Y}}}

\begin{document}
\date{Last compiled \today; last edited  \heuteIst\ or later}

\title{Finite group extensions and the Atiyah conjecture}
\author{Peter Linnell
\thanks{
\protect\href{mailto:linnell@math.vt.edu}{email: linnell@math.vt.edu}\protect\\
\protect\href{http://www.math.vt.edu/people/linnell/}{www: http://www.math.vt.edu/people/linnell/}\protect\\
Fax: ++1 540 231 5960\protect\\
Partially supported by SFB 478, M{\"u}nster}
\and
Thomas Schick\thanks{
\protect\href{mailto:schick@uni-math.gwdg.de}{e-mail: schick@uni-math.gwdg.de}
\protect\\
\protect\href{http://www.uni-math.gwdg.de/schick/}{www http://www.uni-math.gwdg.de/schick/}
\protect\\
Fax: ++49 -551/39 7766, ++49 -551/39 7766\protect\\
Research funded by DAAD (German Academic Exchange Agency)
}\\
}
        
\maketitle

\begin{abstract}
  The Atiyah conjecture for a discrete group $G$ states that the
 $L^2$-Betti numbers of a finite CW-complex with fundamental group $G$
 are integers if $G$ is torsion-free, and in general that they are
 rational numbers with
 denominators determined by the finite subgroups of $G$.

  Here we establish conditions under which the Atiyah conjecture for a
 torsion-free group $G$ implies the Atiyah conjecture for every finite
 extension of $G$.
 The most important requirement is that $H^*(G,\integers/p)$ is
 isomorphic to the 
 cohomology of the $p$-adic completion of $G$ for every prime number $p$. An 
 additional assumption is necessary, e.g.~that the
 quotients of the lower central series or of the derived series are
 torsion-free.

  We prove that these conditions are fulfilled for a certain class of
  groups, which contains in particular Artin's pure braid
  groups (and more generally fundamental groups of fiber-type
  arrangements), free groups,
  fundamental groups of orientable
  compact surfaces, certain knot and link groups, certain
  primitive one-relator groups, and
  products of these. Therefore every finite, in
  fact every elementary amenable extension of these groups satisfies
  the Atiyah conjecture, provided the group does.

  As a consequence, if such an extension $H$ is torsion-free then the
  group ring $\complexs H$ contains no non-trivial zero divisors,
  i.e.\ $H$ fulfills the zero-divisor conjecture.

  In the course of the proof we prove that if these extensions are
  torsion-free, then they have
  plenty of non-trivial torsion-free quotients which are virtually
  nilpotent. All of this applies in particular to Artin's full braid
  group, therefore
  answering question B6 on \href{http://www.grouptheory.info}{www.grouptheory.info}.

  Our methods also apply to the Baum-Connes conjecture.  This is discussed in 
\cite{Schick(2000a)}, where for example the Baum-Connes conjecture is
proved for the full braid groups. 

\end{abstract}

\tableofcontents

\section{Introduction}

In 1976, Atiyah \cite{Atiyah(1976)} constructed the $L^2$-Betti numbers of
a compact
Riemannian manifold. They are defined in terms of the spectrum of the
Laplace operator on the universal covering of $M$. By Atiyah's
$L^2$-index theorem \cite{Atiyah(1976)}, they can be used e.g.~to
compute the Euler
characteristic of $M$.

Dodziuk \cite{Dodziuk(1977)} proved an $L^2$-Hodge de Rham theorem
which gives a
combinatorial interpretation of the $L^2$-Betti numbers, now in terms
of the spectrum of the combinatorial Laplacian. This also generalized
the definition of $L^2$-Betti numbers to arbitrary finite
CW-complexes. The combinatorial Laplacian can be
considered to be a matrix over the integral group ring of the
fundamental group $G$ of $M$. These $L^2$-invariants have been
successfully used to give new results in differential geometry,
topology, and algebra. For example, one can prove certain cases of the 
Hopf conjecture about the sign of the Euler characteristic of
negatively curved manifolds \cite{Gromov(1991)}, or, in a completely
different direction, certain cases of the zero divisor conjecture,
which asserts that $\rationals[G]$ has no non-trivial zero divisors if 
$G$ is torsion-free, following from the Atiyah conjecture about
$L^2$-Betti numbers. This paper is concerned with this algebraic
aspect of the theory of $L^2$-Betti numbers. The methods 
employed will also be algebraic.

The Atiyah conjecture is the central theme of this paper. It is
a prediction of the possible values of the $L^2$-Betti numbers
mentioned above. We start with an algebraic formulation of its statement.

\begin{definition}\label{def:Atiyah}
  For a discrete group $G$, set 
  \begin{equation*}
\lcm(G):=\text{least common multiple of }\{\abs{F} \mid
  F\subgroup G\text{ and }\abs{F}<\infty\},
\end{equation*}
where we adopt the convention that the least common multiple of an
unbounded set is $\infty$.

  If $\lcm(G)<\infty$, the \emph{strong Atiyah conjecture over $KG$} says
  \begin{equation*}
   \lcm(G)\cdot  \dim_G (\ker A)\in \integers \qquad\forall A\in M_n(K G);
\end{equation*}
  where $K\subset \complexs$ will in this paper always be a
 subring of the complex numbers which
  is closed under complex conjugation.
We have used the adjective ``strong" so as to distinguish it from the
conjecture that Atiyah made in \cite{Atiyah(1976)}, which in the case
$K = \mathbb {Q}$ is equivalent to  $\dim_G (\ker A)\in \mathbb {Q}$
for all $A \in M_n(\mathbb {Q}G)$.

  Here $A$ is considered as a bounded operator $A\colon l^2(G)^n\to
  l^2(G)^n$ (acting by the convolution product).  Let $\pr_{\ker A}$ be the
  projection onto $\ker A$. Then
\begin{equation*}
\dim_G(\ker A):=\tr_G(\pr_{\ker A}):=\sum_{i=1}^n \innerprod{\pr_{\ker A}
  e_i,e_i}_{l^2(G)^n},
\end{equation*}
where $e_i\in l^2(G)^n$ is the vector with the identity element of
$G\subset l^2(G)$
at the $i^{\text{th}}$-position and zeros elsewhere.

  The ring $KG$ embeds into the group von Neumann algebra
  $\NeumannN G$, and this into von Neumann's ring of affiliated operators
  $\universalU G$. We
  denote by $DG$ the division closure of $KG$ in $\universalU G$, and
  we let $\Sigma(G)$ be the set of all matrices over $KG$ which become
  invertible in $DG$. For more details about all these objects compare
  e.g.~\cite{Reich(1999)} or \cite[Section 8]{Linnell(1998)}.
\end{definition}

\begin{remark}
  By \cite[Lemma 2.2]{Lueck(1997)} the strong Atiyah
  conjecture as stated in Definition \ref{def:Atiyah} for $K=\rationals$ is
  equivalent to the statement that for every finite CW-complex $X$ with
  fundamental group $G$, the $L^2$-Betti numbers fulfill
  $\lcm(G)b^p_{(2)}(X)\in\integers$. 
\end{remark}

\begin{remark}
  There were more general versions of the Atiyah conjecture for groups
  with $\lcm(G)=\infty$. The geometric version then says that at least
  all relevant $L^2$-Betti numbers should be contained in
  the additive subgroup of $\rationals$ generated by the numbers
  $\abs{F}^{-1}$, $F\subgroup G$ with $\abs{F}<\infty$.

  However, this conjecture is false, as was shown by counterexamples in
  \cite{Grigorchuk-Linnell-Schick-Zuk(2000)} and \cite{Dicks-Schick(2002)}.
\end{remark}

\begin{remark}
  If $G$ is torsion-free, the strong Atiyah conjecture for $KG$ is
  particularly interesting: it implies that $KG$
is a domain, i.e.\ has no non-trivial
  zero-divisors \cite[Lemma 2.4]{Lueck(1997)}. More precisely, the ring
  $DG$ which contains $KG$ is a skew field (compare e.g.~\cite[Lemma
  3]{Schick(1999)}).

  There are of course also purely algebraic methods available to
  attack the above-mentioned zero-divisor problem (and generalizations
  thereof). We will
  give partial results also in such a more general setting. More precisely:
\end{remark}

\begin{definition}\label{def:crossed_product}
  Let $R$ be a division ring and $G$ a group. A \emph{crossed product} $R*G$ is
an associative ring which has underlying set
  the free $R$-module $\bigoplus_{g\in G} R\cdot \overline g$ with basis
  $\{\overline g\mid g\in
  G\}$. 
The multiplication is determined by set maps $\phi \colon G\times G
\to R \setminus \{0\}$, $\sigma \colon G \to \Aut(R)$ satisfying
certain conditions, and is given by $\overline{g_1}\cdot \overline{g_2} =
  \phi(g_1,g_2)\overline{g_1g_2}$, $\overline {g} \cdot r = (\sigma
(g) r) \cdot \overline {g}$; see \cite[p.~2]{Passman(1989)} for further
information.
\end{definition}

We will be interested in the question of when such crossed products are
domains, and in particular when they admit
embeddings into skew fields.

We will throughout the paper only deal with groups $G$ with
$\lcm(G)<\infty$. Note that $G$ is torsion-free if and only if $\lcm(G)=1$.

The strong Atiyah conjecture is true for elementary amenable groups, and
more generally for free groups and their elementary
amenable extensions \cite{Linnell(1993)}. Note that this class of
groups is closed under finite extensions.

The strong Atiyah conjecture is also
known for $\rationals G$ if $G$ is a
residually torsion-free elementary
amenable group (i.e.\ the case $K=\mathbb {Q}$)
\cite[Corollary 4]{Schick(1999)}. A generalization to
$K=\overline\rationals$, the field of algebraic numbers over
$\rationals$ in $\complexs$, is given in
\cite{Dodziuk-Linnell-Mathai-Schick-Yates(2001)}. This applies
e.g.~for pure 
braid groups and
certain positive one-relator groups, but
in these cases with known methods no general
statement about finite extensions can be made.

Observe that the last statement applies in particular to Artin's pure
braid groups
$P_n$,
which are residually torsion-free nilpotent
\cite[Theorem 2.6]{Falk-Randell(1988)}.
Our work was motivated by the question of whether the full braid groups
$B_n$, which are finite torsion-free extensions of $P_n$, also
satisfy the strong Atiyah conjecture.  It was also motivated by the result
of
Rolfsen and Zhu \cite[Theorem 1]{Rolfsen-Zhu(1998)} where
using the property that $B_n$ is left orderable
\cite{Fenn-Greene(2000)}, they proved
that $\mathbb {C}B_n$ has no non-trivial zero divisors.

In this paper, we establish this and much more general results.

The basic idea is: it is hard to deal with extensions with finite
quotient $F$, but easy to handle torsion-free elementary amenable
quotients. We
recall in Section \ref{sec:amenableext} (see Corollary
\ref{corol:atiyah_amext} for a more general statement):
\begin{theorem}\label{theo:tfext}
Let $H$ be a torsion-free group which fulfills the strong
Atiyah conjecture over $KH$.  Suppose we have an
  exact sequence $1\to H\to G\to A\to 1$, where $A$ is torsion-free
  and elementary amenable. Then the strong
Atiyah conjecture for $KG$ is true.
\end{theorem}

For a finite extension $1\to H\to G\to F\to 1$, our goal therefore is to find a
factorization $G\overset{\pi}{\onto}
A\onto F$, where $A$ is elementary amenable and
torsion-free, and
apply Theorem \ref{theo:tfext} to the extension $1\to
\ker \pi \to G\to A\to 1$.

One can further reduce to the case where $F$ is a finite $p$-group for
a prime number $p$. So we have to answer the question: given a torsion-free
group $H$ which fulfills the strong
Atiyah conjecture, under which conditions
on $H$ does every extension $1\to H\to G\to F\to 1$, where $F$ is a
finite $p$-group and $G$ is torsion-free, admit a factorization $G\onto
A\onto F$ with $A$ torsion-free and elementary amenable.

We give such a criterion in Section \ref{sec:braid}. In particular we
prove
\begin{theorem}
  Assume that $H$ has a finite classifying space and 
  that for each prime number $p$, the canonical map
  \begin{equation*}
    H^*({\hat H}^p,\integers/p)\to H^*(H,\integers/p)
  \end{equation*}
from the continuous cohomology of the pro-$p$
  completion of $H$ to the cohomology of $H$ is an isomorphism.
  Suppose infinitely many
  quotients of the derived
  series of $H$ or of the lower central series of $H$ are torsion-free.
Then the above factorization property is fulfilled.
\end{theorem}

The proof uses the pro-$p$ completions as an intermediate
step to show that if no such torsion-free quotient exists, then the
projection $G\to F$ induces a split injective map in
cohomology (we might have to replace $F$ with a non-trivial
subgroup first). Using an Atiyah-Hirzebruch spectral sequence
argument, we
show that the same is true for stable cohomotopy. But then a
fixed-point theorem of Jackowski \cite{Jackowski(1988)} implies that
the map $G\to F$ itself splits, which contradicts the assumption that
$G$ is torsion-free.

All of this can be generalized to extensions $G$ which contain
torsion. Then the requirement is to find a factorization through an
elementary amenable group $A$ which has as little torsion as
possible. The theorem extends to this case, and moreover, the
assumption on the torsion-free quotients can be somewhat relaxed.

The cohomological condition is rather strong. However, we show in
Section \ref{sec:examples} that all conditions are satisfied by
Artin's pure braid group. More generally, consider the following
classes of groups: finitely generated
free groups, surface groups, primitive one-relator groups, knot
groups and primitive link groups.
Using these building blocks, construct $H$ as
an iterated semi-direct product, such that always the quotient acts trivially
on the abelianization of the kernel. Then $H$ fulfills the
conditions of
Theorem \ref{theo:tfext}. By \cite{Falk-Randell(1988)} the
pure braid groups can be constructed in this way, and therefore the
theorem can be applied to them.

We remark that, in the course of the proof, we construct torsion-free
non-trivial elementary amenable quotients of the torsion-free
extensions of our class of groups. This applies in particular to the
full braid groups $B_n$, thus answering question B6 of
\href{http://grouptheory.info/}{http://grouptheory.info/} for all $n\ge 1$.

Explicitly (compare Corollary \ref{corol:Atiyah_for_braid}) we have
\begin{theorem}
   Let $B_n$ be the $n$th classical full Artin braid group. Then the
  strong Atiyah conjecture for $\overline{\rationals}B_n$ is
  true, where $\overline{\rationals}$ is the field of algebraic numbers
    in $\complexs$. Moreover, $B_n$ is residually torsion-free nilpotent-by-finite.
\end{theorem}

The algebraic methods of this paper can also be applied to the
Baum-Connes conjecture, to establish this conjecture in more new
cases. This is discussed in 
\cite{Schick(2000a)}, where for example the Baum-Connes conjecture is
proved for the full braid groups. 

\section{Atiyah conjecture and extension with elementary amenable
   quotient}
\label{sec:amenableext}

In this section, we recall some general results about the Atiyah
conjecture, proved essentially in \cite{Linnell(1998),Linnell(1993)}.
We will frequently make use of results in
\cite{Linnell(1998),Linnell(1993)}. There, these results are stated
only for $K=\complexs$. However, everything generalizes immediately to
the more general rings $K$ we are considering.

\begin{notation}
We let $\mathcal{A}$ denote the class of finitely generated
abelian-by-finite groups.  If $\mathcal {X}$ and $\mathcal {Y}$ are
classes of groups, then $G \in L\mathcal {X}$ means that every
finitely generated subgroup of $G$ is contained in an $\mathcal
{X}$-group,  and $G \in \mathcal {X} \mathcal{Y}$ means that $G$ has
a normal $\mathcal{X}$-subgroup $H$ such that $G/H$ is a $\mathcal
{Y}$-group.

  A group is called \emph{elementary amenable}, if it belongs to
  the smallest class of groups which contains all solvable groups and
  all finite groups, and
  which is closed under extensions and under directed unions.

Also $\mathbb {N}$ will denote the positive
integers (so 0 is excluded), and $\naturals_0:=\naturals\cup\{0\}$.
\end{notation}

\begin{proposition}\label{prop:artinian_equivalence}
  Assume that $G$ is a discrete group with $\lcm(G)<\infty$ and  $DG$, the
  division closure of $KG$ in $\universalU G$, is Artinian.

  Then $DG$ is semi-simple Artinian. Moreover, if $L$ is a positive integer
  then 
  \begin{equation*}
L\cdot\dim_G(\ker
  A)\in\integers \quad\text{for each } A\in M_n(KG)
\end{equation*}
is
  equivalent to 
  \begin{equation*}
L\cdot\tr_G(e)\in\integers\quad\text{for each projection } e\in DG.
\end{equation*}

  In particular, $G$
  fulfills the strong Atiyah conjecture if and only if for each
  projection $e\in DG$ we have $\lcm(G)\cdot
  \tr_G(e)\in\integers$.
\end{proposition}
\begin{proof}
  Semi-simplicity follows from \cite[Lemma 9.4]{Linnell(1998)}.

  The remaining statements  are also standard. The strong
Atiyah conjecture
  implies the statement about the trace of the projections in $DG$
  using Cramer's rule \cite[Proposition 7.1.3]{Cohn(1985)} as
  in the proof of \cite[Lemma 3.7]{Linnell(1993)}. 

  The converse follows by noting that the semi-simple Artinian subring $DG$ of
  $\universalU G$ is $*$-regular, and hence for $A\in M_n(KG)\subset
  M_n(DG)$ we have $A M_n(DG)=eM_n(DG)$ for a unique projection $e\in
  M_n(DG)$. Thus the $G$-dimension of the range of $A$
is equal to the $G$-trace of the projection $e$, the $G$-dimension of
the kernel of $A$ is equal to the $G$-trace of the projection $1-e$,
and we can apply
  \cite[Corollary 13.11]{Linnell(1998)} to see that $L
  \tr_G(e)\in\integers$.
\end{proof}

\begin{lemma}\label{lem:trivial_finite_ext_bound}
  Assume that $H$ is a discrete group, $DH$ is Artinian and there is a
  positive integer
  $L$ such that 
  \begin{equation*}
    L \cdot \tr_H(e)\in\integers\qquad\text{if $e\in DH$ is a projection.}
  \end{equation*}
  Suppose $H\subgroup G$ and the index $[G:H]<\infty$. Then $DG$ is
  semi-simple Artinian and
  \begin{equation*}
    [G:H]\cdot L\cdot \tr_G(e)\in\integers\qquad\text{if $e\in DG$ is a
    projection.}
  \end{equation*}
\end{lemma}
\begin{proof}
  This follows immediately from \cite[Lemma 9.4]{Linnell(1998)}, Proposition
  \ref{prop:artinian_equivalence} and \cite[Proposition
  4]{Schick(1999)}. Alternatively, a direct argument is given in
  \cite[Proof of Lemma 13.12]{Linnell(1998)}.
\end{proof}

\begin{lemma}\label{lem:psubgroup_check}
  Assume that $H$ is a group with $\lcm(H)<\infty$.
Let $G$ be a finite extension of $H$ with quotient group $G/H$ and
projection $\pi \colon G\to G/H$.

  If there is
  $L\in\naturals$ such that for the inverse images $H_p\subgroup G$ of
  all Sylow
  subgroups of $G/H$
  \begin{equation*}
    L\cdot \lcm(H_p) \cdot\dim_{H_p}(\ker A)\in\integers\qquad\forall A\in
    M_n(KH_p),
  \end{equation*}
  then
  \begin{equation*}
    L\cdot \lcm(G)\cdot \dim_G(\ker A)\in\integers\qquad\forall A\in M_n(KG).
  \end{equation*}
  In particular, the strong
Atiyah conjecture for all the $H_p$ implies the
  strong Atiyah conjecture for $G$.
\end{lemma}
\begin{proof}
  Write $\abs{G/H}=p_1^{n_1}\dots p_k^{n_k}$ where the $p_j$ are
  different primes. Let $H_{p_j}$ be the inverse image under $\pi$ of
  a $p_j$-Sylow 
  subgroup of $G/H$. This implies that $p_j$ does not divide
  $[G:H_{p_j}]=\abs{G/H}/p_j^{n_j}$, i.e.~the greatest common divisor of
  all $[G:H_{p_j}]$ is one.

For trivial reasons we have the following divisibility
  relation:
  \begin{equation*}
    \lcm(H_{p_j}) \divides \lcm(G) \qquad\forall j.
  \end{equation*}
  If $A\in M_n(KG)$ then, by assumption and by \cite[Proposition
  4]{Schick(1999)},
  \begin{equation*}
    \begin{split}
      & L \cdot \lcm(H_{p_j})\cdot [G:H_{p_j}] \cdot \dim_G(\ker(A))
      \in \integers,\\
\text{consequently}\quad & L \cdot [G:H_{p_j}] \cdot \lcm(G)\cdot \dim_G(\ker(A))
      \in \integers.
    \end{split}
  \end{equation*}
  Passing to  the greatest common divisor we get
  \begin{equation*}
    L \cdot \lcm(G)\cdot \dim_G(\ker A) =
    g.c.d.\{[G:H_{p_j}]\}\cdot L \cdot\lcm(G)\cdot \dim_G(\ker A)\in
    \integers. \hfill\qed
  \end{equation*}
\renewcommand{\qed}{}
\end{proof}

The following lemma is needed to arrive at the purely algebraic
results ---beyond the Atiyah conjecture--- we want to obtain.
Information on Ore domains can be found in \cite[\S
4]{Linnell(1998)}.

\begin{lemma}\label{Tvirabelian}
   Let $G$ be an elementary amenable group, let $k$
    be an Ore domain, and let $k*G$ be a crossed product. If $k*F$ is a
    domain for all finite subgroups $F$ of $G$, then $k*G$ is an Ore
    domain, consequently it has a right quotient ring which is a skew field.
\end{lemma}
\begin{proof}
First suppose $G$ is finitely generated abelian-by-finite.
If $k$ is a skew field, this is the statement
  of \cite[Corollary 4.5]{Linnell(1998)}; it also
 follows from \cite[Lemma 4.1]{Kropholler-Linnell-Moody(1988)}.
For an arbitrary Ore domain $k$, let $S = k \setminus \{0\}$.  Then
$kS^{-1}$ is a skew field and
$kS^{-1} * H \cong (k*H)S^{-1}$ for all subgroups $H$ of $G$.  From
the first sentence, $kS^{-1}*G$ is an Ore domain, and we deduce that
$k*G$ is an Ore domain.

  The general case is proved by transfinite induction. Let
  $\mathcal {A}$ be the
  class of all finitely generated abelian-by-finite groups. Let $\amenableGroups_0$ denote the
  class of finite groups, and define inductively
$\amenableGroups_{\alpha} := (L
  \amenableGroups_{\alpha-1})\mathcal {A}$ if $\alpha$ is a
successor ordinal, and $\amenableGroups_{\alpha} := \bigcup_{\beta <
\alpha} \amenableGroups_{\beta}$ if $\alpha$
is a limit ordinal (cf.\ \cite[Lemma 4.9]{Linnell(1993)}).
Then, by the proof of \cite[Lemma 4.9(ii)]{Linnell(1993)}, the class of
elementary amenable groups is $\bigcup_{\alpha \ge 0}
\amenableGroups_{\alpha}$. 

 Let $\alpha$ be the least ordinal
such that $G \in \amenableGroups_{\alpha}$.  Then $\alpha$ cannot be a
limit ordinal, and the result is certainly true if $\alpha = 0$.
Therefore, we may write $\alpha = \beta + 1$, and assume that the
result is true for all groups in $\amenableGroups_{\beta}$.  If $U \in
L\amenableGroups_{\beta}$, then since $k*E$ is an Ore domain for every
finitely generated subgroup $E$ of $U$, it is clear that
$k*U$ is an Ore domain.
Hence, we may assume that the result is true for groups in
$L\amenableGroups_{\beta}$.  Now $G$ has a normal subgroup $D$ such that $D
\in L\amenableGroups_{\beta}$ and $G/D \in \mathcal {A}$.  Furthermore, if
$E/D$ is a finite subgroup of $G/D$, then $E \in
L\amenableGroups_{\beta}$ by the proof of \cite[Lemma
4.9(iv)]{Linnell(1993)}. The result follows from the case first dealt
with, where we replace $G$ with $G/D$ and $k$ with $k*D$.
\end{proof}

\begin{proposition}\label{prop:amext}
  Let $1\to H\to G\to A\to 1$ be an exact sequence of groups. 
Fix $K=\bar K\subset\complexs$ and $L\in\naturals$.
Suppose $A$ is elementary amenable and $DH$ is Artinian.
For a subgroup $E\subgroup A$, let $H_E$ denote the inverse
  image of $E$ in $G$. 
\begin{enumerate} [\normalfont (i)]
\item
Assume that we can find $L\in\naturals$ such that for all finite subgroups
$E\subgroup A$,
  \begin{equation*}
    L\cdot \tr_{H_E}(e)\in\integers\qquad\text{for each projection }e\in DH_E.
  \end{equation*}
  Then $DG$ is semi-simple Artinian and
  \begin{equation*}
    L\cdot\tr_{G}(e)\in\integers\qquad\text{for each projection }e\in DG. 
  \end{equation*}
\item
With the setup in (i), assume further that the identity map on $KH$
  induces an isomorphism 
  \begin{equation*}
KH_{\Sigma(H)}\to DH,\quad\text{where $\Sigma(H)$ is defined in Definition
  \ref{def:Atiyah}}.
\end{equation*}
  Then the identity map on $KG$ induces an isomorphism
  $KG_{\Sigma(G)}\to DG$.
\item
  Let $k$ be a domain and $k*G$ be a crossed product. Assume that 
  $k*H$ embeds into a skew field $D_H$ and the twisted action of $A$
  on $k*H$ extends to a twisted action on $D_H$ (in
  particular, we can form $D_H*A$). Assume that $D_H*E$ is a domain for all
  finite subgroups $E$ of $A$. 
Then $k*G$ can be embedded in a skew field.
\end{enumerate}
\end{proposition}
\begin{proof}
We shall prove parts (i) and (ii) by the standard transfinite
induction argument as used in Lemma \ref{Tvirabelian}, and we shall
adopt the same notation as defined in the second paragraph there.
The proof of (i) is similar to 
\cite[Proposition 8]{Schick(1999)}.
Thus we have two cases to consider:
\begin{enumerate}
\item $A$ has a normal subgroup $B$ such that $A/B \in \mathcal {A}$
and the result is known with $B$ in place
of $A$.  Then in (i) we use \cite[Lemma 13.10, parts (i) and
(iii)]{Linnell(1998)}, while in (ii) we use
\cite[Lemma 13.10(ii)]{Linnell(1998)}, noting
that $\mathbb {C}$ can be replaced by any subfield of $\mathbb {C}$
which is closed under complex conjugation.

\item $A$ is the directed union of groups $A_i$, and the result is
known with $A_i$ in place of $A$ for all $i$.  In (i) we apply
\cite[Lemma 13.5, parts (i) and (ii)]{Linnell(1998)}, while in (ii)
we apply
\cite[Lemma 13.5(iii)]{Linnell(1998)}.  Again we note that $\mathbb
{C}$ can be replaced by any subfield of $\mathbb {C}$ which is closed
under complex conjugation.
\end{enumerate}

  For the proof of part (iii) of the proposition, observe that
  $(k*H)*G/H=k*G$. By Lemma \ref{Tvirabelian} the
  ring $D_H *A$, and hence also $(k*H)*A=k*G$, embeds into a skew field.
\end{proof}

\begin{corollary}\label{corol:atiyah_amext}
  Under the assumptions of Proposition \ref{prop:amext}, if
  $L=\lcm(G)$ then the strong Atiyah conjecture is true for $KG$. 

  If, for every finite subgroup $E$ of $A$, $KH_E$ fulfills the strong Atiyah
  conjecture, then
we may take
  $L=\lcm(G)$. Hence under this
  condition the strong Atiyah conjecture for $KG$ holds.
\end{corollary}
\begin{proof}
  This follows immediately from Proposition \ref{prop:amext} and
  Proposition \ref{prop:artinian_equivalence}.
\end{proof}

\begin{corollary}
  If $G$ is elementary amenable and $\lcm(G)<\infty$, then the strong
  Atiyah conjecture holds for $\complexs G$.
\end{corollary}
\begin{proof}
  This follows since the strong Atiyah conjecture is true for
  $\complexs G$ if $G$ is a trivial group.
\end{proof}

\section{More positive results about the Atiyah conjecture}
\label{sec:positive_results}

In this section, we recall a few more known results about the
Atiyah conjecture, proved in particular in \cite{Schick(1999)} and \cite{Dodziuk-Linnell-Mathai-Schick-Yates(2001)}.

\begin{definition}
  Let $\extendedC$ be the smallest non-empty class of groups such that:
  \begin{enumerate}
  \item \label{aex} If $G$ is torsion-free and $A$ is elementary
    amenable, and we have a projection
    $p\colon G\to A$ such that $p^{-1}(E)\in\extendedC$ for every
    finite subgroup $E$ of $A$, then $G\in\extendedC$.
   \item \label{sgr} $\extendedC$ is subgroup closed.
   \item\label{lim} Let $G_i\in \extendedC$ be a
  directed system of groups with homomorphisms $\phi_{ij}\colon G_i\to
  G_j$ or $\phi_{ij}\colon G_j\to G_i$, respectively, for $i<j$,
    and assume that $G$ is its direct limit or inverse limit,
    respectively. Then
    $G\in\extendedC$.
\end{enumerate}
\end{definition}

\begin{theorem}\cite[Theorem 1]{MR1894160},\cite[Theorem 1.4]{Dodziuk-Linnell-Mathai-Schick-Yates(2001)}\label{theo:Atiyah_for_D}\\
  Assume that $\overline \rationals$ is the algebraic closure of
  $\rationals$ in $\complexs$ and $G\in\extendedC$. Then the strong
  Atiyah conjecture is true for $\overline\rationals G$.
\end{theorem}

\begin{corollary}\label{corol:Atiyah_for_restfnilpotent}
  The strong Atiyah conjecture is true for $\overline\rationals G$, if
  $G$ is residually torsion-free nilpotent or if $G$ is residually
  torsion-free solvable.
\end{corollary}

\section{Finite extensions and the Atiyah conjecture}
\label{sec:braid}

In this section, we give an abstract criterion when the strong Atiyah
conjecture for a torsion-free group $G$ implies the strong Atiyah conjecture
for all of its extensions with finite quotient, and more generally for all
extensions with elementary amenable quotient.

The strategy is the following:
making use of Lemma \ref{lem:psubgroup_check}, instead of
investigating all of the finite extensions at once, we study
them one prime at a time (concentrating on those extensions where the
quotient is a finite $p$-group). 

The crucial condition is cohomological completeness, defined in
Definition~\ref{def:cohom_complete}. This is expressed in terms of
pro-$p$ completions and their continuous cohomology, following
the definitions and
notation of \cite{Serre(1997)}.

Our condition is used to show that each finite extension $G$
of the given group $H$ is an extension of a suitable subgroup of $H$
with elementary amenable quotient $A$,
such that $A$ contains as little torsion as theoretically
possible. This allows us
to apply standard results to conclude the Atiyah conjecture for $G$.

Before we obtain these results, we start with a comparably easy other
condition, which is also of cohomological type.

\subsection{Groups with Euler characteristic equal to one}

\begin{theorem}\label{theo:eulerchar}
Assume that $H$ is a group with finite classifying space $BH$, and with
Euler characteristic
\begin{equation*}
  \chi(H)=\chi(BH)\in \{1,-1\}.
\end{equation*}
Suppose the strong Atiyah conjecture is true for $KH$ and we have an exact
sequence
\begin{equation*}
  1\to H\to G\to Q\to 1,
\end{equation*}
such that $Q$ is elementary amenable and $\lcm(G)<\infty$. 

Then the strong Atiyah conjecture is true for $KG$.
\end{theorem}
\begin{proof}
  Since $H$ has a finite classifying
  space, it is torsion-free. Because of Proposition \ref{prop:amext}, Corollary
  \ref{corol:atiyah_amext} and Lemma
  \ref{lem:psubgroup_check}, it is sufficient to assume that $Q$ is a
  finite $p$-group.

  In this situation, we will prove that $\abs{Q}$ divides $\lcm(G)$. It
  then follows from Lemma \ref{lem:trivial_finite_ext_bound} that $KG$
  fulfills the strong Atiyah conjecture.

  Since the Euler characteristic is multiplicative \cite[Proposition
  7.3]{Brown(1982)}, $\chi(G)=\pm \frac{1}{\abs{Q}}$. By \cite[Theorem
  9.3]{Brown(1982)}, $\lcm(G)\cdot\chi(G)\in\integers$, therefore
  $\abs{Q}$ divides $\lcm(G)$. This concludes the proof.
\end{proof}

This result is less interesting than it might seem to be because a
non-trivial 
finite extension of a group $H$ with $\chi(H)=\pm 1$ can not be
torsion-free (which is well known and follows from the argument in the
proof of Theorem \ref{theo:eulerchar}).

\subsection{Cohomology of groups and their completions}

\begin{definition}
  For a discrete group $G$ and a prime number $p$, let ${\hat G}^p$ denote the
  \emph{pro-$p$ completion} of $G$. This is the inverse limit of the
  system of finite $p$-group quotients of $G$.  Then ${\hat G}^p$ is a
  compact totally disconnected group, and we have a (not necessarily
  injective) canonical homomorphism $G\to{\hat G}^p$. The cohomology of
  a profinite group like ${\hat G}^p$ in this paper will always be
  its \emph{Galois cohomology}  \cite{Serre(1997)} (sometimes also called
  \emph{continuous cohomology}).
\end{definition}

\begin{definition}\label{def:cohom_complete}
  A discrete group $G$ is called \emph{cohomologically complete} if for each
  prime number $p$ the homomorphism 
  \begin{equation}\label{eq:compare}
H^*({\hat G}^p,\integers/p)\to
  H^*(G,\integers/p)
\end{equation}
is an isomorphism. Here the action on $\integers/p$ is assumed
  to be trivial. Define \emph{cohomological $p$-completeness} by requiring
  the same only for $p$.
\end{definition}

\begin{lemma}\label{lem:low_completeness}
  Let $G$ be any group. Then in low degrees we get isomorphisms
  \begin{equation*}
    \begin{split}
      & H^0(\hat{G}^p,\integers/p)\xrightarrow{\iso}
      H^0(G,\integers/p)\\
      & H^1(\hat{G}^p,\integers/p)\xrightarrow{\iso}
      H^1(G,\integers/p),
  \end{split}
\end{equation*}
and an injection
\begin{equation*}
  H^2(\hat{G}^p,\integers/p)\into H^2(G,\integers/p).
\end{equation*}
\end{lemma}
\begin{proof}
  This is well known; for the convenience of the reader we give an
  argument. The statement for $H^0$ is obvious from the
  definition. Next, assume $V$ is a normal subgroup of $G$ such that
  $G/V$ is a finite $p$-group. Then we get the following exact
  sequence in low cohomology degrees
  \begin{multline*}
    0\to H^1(G/V,\integers/p)\to H^1(G,\integers/p) \\ \to
    H^1(V,\integers/p)^G \to H^2(G/V,\integers/p)\to H^2(G,\integers/p).
  \end{multline*}
  Passing to the direct limit, as $V$ varies over all normal subgroups 
  of $p$-power index, we obtain the following exact sequence (since
  taking direct limits is exact)
  \begin{multline*}
    0\to H^1(\hat{G}^p,\integers/p)\to H^1(G,\integers/p)\\
    \to \lim_{V} 
    H^1(V,\integers/p)^G \to H^2(\hat{G}^p,\integers/p)\to H^2(G,\integers/p).
  \end{multline*}
  To prove the remaining statements, we only have to prove that the
  limit in the middle is zero. Now, given any such $V$, an element
  $x\in H^1(V,\integers/p)$ is a homomorphism $x\colon
  V\to\integers/p$. By Lemma \ref{lem:subgroupintersect}, the kernel of 
  $x$ (which has index $1$ or $p$ in $V$) contains a subgroup $W$ which 
  is normal in $G$, and such that $V/W$ is a finite $p$-group. Then
  also $G/W$ is a finite $p$-group, and by construction of $W$ the
  image of $x$ (which is just the restriction of the homomorphism) in
  $H^1(W,\integers/p)$, and therefore also in the limit, is
  zero. Since $V$ and  $x$ were arbitrary, $\lim
  H^1(V,\integers/p)=0$, as desired.
\end{proof}

\begin{example}\label{ex:free_completeness}
  Let $F$ be an arbitrary free group (not necessarily finitely
  generated). Then $F$ is cohomologically complete.
\end{example}
\begin{proof}
  The case where $F$ is finitely generated is
  contained in Example \ref{ex:primitive_one_relator} and Proposition
  \ref{prop:primitive_one_relator}. The general case is proved in
  \cite{Linnell-Schick(2000c)}, just observing that $F$ and $\hat F^p$ have
  cohomological dimension $1$ and that \eqref{eq:compare} always is an
  isomorphism in degrees $0$ and $1$.
\end{proof}

\begin{lemma}\label{lemma:coeffiso}
  Assume that $G$ is cohomologically $p$-complete and $A$ is a finite discrete
  $p$-primary ${\hat G}^p$-module. Then the homomorphism $H^*({\hat
    G}^p,A)\to H^*(G,A)$ is an isomorphism.

  Assume that $\profinQ$ is a profinite group such that
  $H^k(\profinQ,\integers/p)$ is
  finite for every $k$. If $A$ is a discrete finite $p$-primary
  $\profinQ$-module, then $H^k(\profinQ,A)$ is finite for every
  $k\in\naturals_0$.
\end{lemma}
\begin{proof}
  Because of the corollary to \cite[Proposition 20]{Serre(1997)} $A$
  has a composition series whose successive quotients are isomorphic
  to $\integers/p$. We can now use induction on the length of this
  composition series, and the five lemma applied to the long exact
  cohomology sequences 
  associated to short exact sequences of coefficients (on the discrete
  as well as on the pro-$p$ side) to conclude the proof. Existence of
  these long exact sequences 
  is proved in \cite[Theorem 9.3.3]{Wilson(1998)} and
  \cite[Proposition 6.1]{Brown(1982)}. Compatibility of the boundary
  maps with the comparison maps follows from the
  constructions.
\end{proof}

\begin{definition}
  A set of subgroups $\mathcal{U}$ of a group $G$ is cofinal among a
  set of subgroups $\mathcal{V}$ of $G$, if for each $V\in\mathcal{V}$
  there is $U\in\mathcal{U}$ with $U\subset V$.
\end{definition}

\begin{definition}
  The subgroups in the \emph{lower central series} of a group $G$ are
  called $\gamma_n(G)$. They are inductively defined by
  $\gamma_1(G):=G$ and $\gamma_{n+1}(G):=[\gamma_n(G),G]$. For two
  subgroups $U,V$ of $G$, $[U,V]$ denotes the subgroup generated by
  commutators $[u,v]:=u^{-1}v^{-1}uv$ with $u\in U$ and $v\in V$.

  The group $G$ is called \emph{nilpotent} (\emph{of class $\le n-1$}),
  if $\gamma_n(G)=\{1\}$ for
  some $n\in\naturals$.

  We define the \emph{derived series} $G^{(n)}$ inductively by
  $G^{(1)}:=G$ and $G^{(n+1)}:=[G^{(n)},G^{(n)}]$. A group is called
  \emph{solvable} (\emph{of class $\le n-1$}) if $G^{(n)}=\{1\}$.
\end{definition}

The following lemma will be used in the construction of cofinal sets of
subgroups at various places.

\begin{lemma}\label{lem:product_of_finite_groups}
  Let $F$ be a finite group, $I$ any index set. Then the product
  $G:=\prod_{i\in I} F$ is locally finite and in particular, elementary amenable.
\end{lemma}
\begin{proof}
  We show that every finitely generated subgroup of $G$ is 
  finite. Let $g_1,\dots,g_n$ with $g_k=(f^k_i)_{i\in I}$, $f^k_i\in
  F$, be a finite collection of elements of $G$, and let $H$ be the
  subset of $G$ generated by $g_1,\dots,g_n$. Define the map
  $\mu\colon I\to F^n$ by $\mu(i):= (f_i^1,\dots, f_i^n)$. The target
  of $\mu$ is a finite set. Therefore, there is a finite subset
  $J\subset I$ such that for each $i\in I$ there is $j\in J$ with
  $\mu(i)=\mu(j)$.

  Let $p\colon G=\prod_{i\in I} F\to \prod_{j\in J} F=:E$ be the
  obvious projection. Observe that $E$ is finite since $J$ is
  finite. We prove that the restriction of $p$ to $H$ is injective,
  which implies that the finitely generated subgroup $H$ of $G$ is
  finite.

  To do this, let $x=g_{r_1}^{n_1}\cdots g_{r_N}^{n_N}\in H$ be an
  element of the kernel of $p$. That means that
  $(f^{r_1}_j)^{n_1}\cdots (f^{r_N}_j)^{n_N} =1\in F$ for each $j\in
  J$. But for an arbitrary $i\in I$ we find $j\in J$ with
  $f^{r_k}_i=f^{r_k}_j$ for each $k=1,\dots,N$ by our choice of
  $J$. In other words, $1=g_{r_1}^{n_1}\cdots g_{r_N}^{n_N} =x \in G$, 
  i.e.~the kernel of $p|_H$ indeed is trivial.
\end{proof}

\begin{lemma}\label{lem:subgroupintersect}
  Let $\{U_i\}_{i\in I}$ be a set of normal subgroups of a group
  $G$. Define $V:=\bigcap_{i\in I} U_i$.
  \begin{enumerate}
  \item If $G/U_i$ is torsion-free, or nilpotent of class
    $\le n$, or solvable of class $\le n$, respectively, for every
    $i\in I$, then $G/V$ has the same property.
  \item If $I$ is finite and $G/U_i$ is finite, or is a finite
    $p$-group, or is elementary amenable, respectively, for every
    $i\in I$, then $G/V$ has the same property.
  \end{enumerate}
  Fix a normal subgroup $U$ of $G$ and a subset $A$ of the
  set of all automorphisms of $G$. Set $V:=\bigcap_{\alpha\in
  A}\alpha(U)$. 
\begin{enumerate}
\item If $G/U$ is torsion-free, or nilpotent, or solvable,
  respectively, then $G/V$ has the same property.
\item If $A$ is finite and $G/U$ is finite, or is a finite $p$-group,
  or is elementary amenable,
  respectively, then $G/V$ has the same property.
\item Assume that $G$ is finitely generated and $G/U$ is finite or is
  a finite $p$-group, respectively. Then $G/V$ is also
 finite or a finite $p$-group, respectively. In particular, $U$
 contains a characteristic subgroup $W$ such that $G/W$ is finite, or
 a finite $p$-group, respectively.
\item If $G$ is finitely generated and $G/U$ is solvable-by-finite
(i.e.\ has a solvable subgroup of finite index),
then $G/V$ is also solvable-by-finite.
\end{enumerate}
\end{lemma}
\begin{proof}
  Observe that $V$ is the kernel of the projection $G\onto \prod_{i\in
  I} G/U_i=:P$. Consequently, $G/V$ is a subgroup of the infinite
  product $P$. Now $P$ is torsion-free, or nilpotent of class $\le n$,
  or solvable of class $\le n$, respectively, if and only if the same
  is true for every factor $G/U_i$. If $I$ (and therefore the product)
  is finite, then we  conclude also that $P$ is finite, or a finite
  $p$-group, or elementary amenable, if every $G/U_i$ has the same property.
  Since all these properties are inherited
  by subgroups, the first part of the lemma follows.

  For the second part set $I:=A$ and $U_\alpha:=\alpha(U)$ for $\alpha\in
  I$. The 
first two statements
  now follow from the first part by observing that $G/\alpha(U)$ is
  isomorphic to $G/U$ for every $\alpha\in A$.  For the third
  statement observe that, since $G$ is finitely generated, for each
  fixed finite 
  group $E$ there are only finitely many different homomorphisms from
  $G$ to $E$. Consequently, since all the groups $G/\alpha(U)$ are
  isomorphic finite groups, there is only a finite number of different 
  subgroups $\alpha(U)$ (being the kernel of finitely many
  homomorphisms), and the statement follows from the second statement.
 Putting this together with the first statement yields the fourth
statement.
\end{proof}

\begin{lemma}\label{lemma:shortexactcompletion}
  Let $1\to H \to G\xrightarrow{\pi} Q\to 1$ be an exact
  sequence of discrete groups where $H = \ker \pi$,
  and let $p$ be a prime number. Then the following sequence of pro-$p$ completions
  is exact:
  \begin{equation*}
    {\hat H}^p\to {\hat G}^p\xrightarrow{\hat \pi} {\hat Q}^p\to 1.
  \end{equation*}
  The kernel of $\hat\pi$ is exactly the closure of $H$ in ${\hat G}^p$.
  If the subgroups $U\cap H$, where
  $U$ runs through all normal subgroups of $G$ with $p$-power index,
  are cofinal among all normal subgroups of $H$ with $p$-power index,
  then we get the exact sequence
  \begin{equation*}
    1\to     {\hat H}^p\to {\hat G}^p\xrightarrow{\hat\pi}{\hat Q}^p\to 1 .
  \end{equation*}
  This is the case in particular if $Q$ is a finite $p$-group (which
  implies $Q={\hat Q}^p$).
\end{lemma}
\begin{proof}
  The first statement is a standard fact; compare \cite[Exercise
  1.6.(10)]{Wilson(1998)} and use the property that the image of ${\hat
H}^p$ in ${\hat G}^p$ is closed, because ${\hat H}^p$ is compact and
${\hat G}^p$ is Hausdorff.  The kernel of $\hat \pi$
  is the closure in ${\hat G}^p$ of $H$, in other words is the
  completion of $H$ with respect to the system of normal subgroups
  $H\cap U$, where $U\normalsubgroup G$ and $G/U$ is a finite
  $p$-group. For the second statement, we have therefore to show that
  this completion coincides with ${\hat H}^p$. By \cite[Exercise
  1.6(4)]{Wilson(1998)} this is the case under the assumption that the
  system $H\cap U$ is cofinal among all normal subgroups of $H$ of
  $p$-power index.

  Assume now that $Q$ is a finite $p$-group.

  It remains to check that whenever $V\normalsubgroup H$ such that
  $H/V$ is a finite
  $p$-group, we can find $U\normalsubgroup G$ with  $G/U$ another
  finite $p$-group and $U\cap H\subset V$. The intersection $W$ of all
  (finitely many)
  conjugates $V_i$ of $V$ in $G$ is contained in $V$ and is normal
  in $G$. Its index
  in $G$ is a power of $p$ if and only if the same is true for its
  index in $H$ (since $G/H$ is a finite $p$-group). The latter
  follows from Lemma \ref{lem:subgroupintersect}. Therefore $W$ is the
  normal subgroup we
  had to find, and the last assertion of the lemma follows.
\end{proof}

The following results establish cohomological $p$-completeness.

\begin{proposition}\label{tfcompletion}
  Assume that $H$ is a discrete group which
  is cohomologically $p$-complete. Let 
  \begin{equation*}
    1\to H \to G \to Q \to 1
  \end{equation*}
  be an exact sequence of groups, where $Q$ is a finite
  $p$-group. Then $G$ is also cohomologically $p$-complete.
\end{proposition}
\begin{proof}
  The exact sequence $1\to H\to G\to Q\to 1$ induces by Lemma
  \ref{lemma:shortexactcompletion} an
  exact sequence
  \begin{equation*}
    1 \to {\hat H}^p\to {\hat G}^p\to Q \to 1,
  \end{equation*}
  
  Now we have two spectral sequences with $E_2^{l,q}$-terms
  \begin{equation*}
 H^l(Q,H^q({\hat H}^p,\integers/p)) 
\qquad\text{or}\qquad
    H^l(Q,H^q(H,\integers/p))
  \end{equation*}
  converging to
 $H^{l+q}({\hat G}^p,\integers/p)$
 or
 $ H^{l+q}(G,\integers/p)$,
 respectively. By
  assumption on $H$ the natural map between Galois and ordinary
  cohomology induces an isomorphism between the $E_2$-terms, and
  consequently the isomorphism between $H^*({\hat G}^p,\integers/p)$
  and $H^*(G,\integers/p)$ we had to establish. 
\end{proof}

\begin{definition}\label{def:unipotent_action}
  Let $V$ be a vector space (over any field), or let $V$ be an abelian
  group, and $A\colon V\to V$ an
  automorphism of $V$. We say $A$ acts \emph{unipotently} on $V$, if
  $(A-\id_V)^n=0$ for some $n\in\naturals$.

  If a group $G$ acts on $V$ via $\alpha\colon G\to \Aut(V)$, we call the action of $G$
  \emph{unipotent} if there is $m\in\naturals$ such that
  $(\alpha(g_1)-\id_V)\cdots(\alpha(g_m)-\id_V)=0$  for all
  $g_1,\dots,g_m\in G$.
\end{definition}

\begin{theorem}\label{theo:extensions_and_cohomological_completeness}
  Let $1\to H\to G \to Q\to 1$ be an exact sequence of groups. Assume
  that $H$ is finitely generated, and that $H$ as well as $Q$ are
  cohomologically $p$-complete for some prime number $p$. Assume that
  $H^i(H;\integers/p)$ is finite for all $i\in\naturals$ (this is true 
  e.g.~if $H$ has a finite classifying space), and that
  each element of $Q$ 
  acts unipotently on $H_1(H;\integers/p)$. Then $G$ is
  cohomologically $p$-complete.

  If, moreover, $H^k(Q,\integers/p)$ is finite for each
  $k\in\naturals$, then $H^k(G,\integers/p)$ is finite for each
  $k\in\naturals$.
\end{theorem}

This theorem with proof is modeled after \cite[Exercise 1) and 2) on
p.~13]{Serre(1994)}. In its present formulation, it was suggested to
us by Guido Mislin. To prove it, we first need a sequence of auxiliary
results.

\begin{proposition}\label{prop:completeness_and_subgroups}
  Assume that $G$ is cohomologically $p$-complete, and $U$ is a normal
  subgroup of $G$ such that $G/U$ is a finite $p$-group. Then $U$ is
  cohomologically $p$-complete, as well.
\end{proposition}
\begin{proof}
  We have the exact sequence $1\to U\to G\to G/U\to 1$ and, by Lemma
  \ref{lemma:shortexactcompletion} an exact sequence $1\to
  \hat{U}^p\to \hat{G}^p\to G/U\to 1$. Since $U$ has finite index in
  $G$, $\coind_U^G\integers/p$ is a finite $G$-module and, therefore,
  it is the restriction of $\coind_{\hat{U}^p}^{\hat{G}^p}\integers/p$ 
  to $G$, which is itself a finite discrete $\hat{G}^p$-module. The
  discrete Shapiro isomorphism \cite[Proposition 6.2]{Brown(1982)}
  \begin{equation*}
    \phi\colon H^k(G,\coind_U^G\integers/p)\to H^k(U,\integers/p)
  \end{equation*}
  is given as composition
  \begin{equation*}
  H^k(G,\coind_U^G\integers/p)\xrightarrow{i^*}H^k(U,\res_U\coind_U^G
  \integers/p)\to H^k(U,\integers/p),
\end{equation*}
the first map being induced from the inclusion $i\colon U\into G$,
and the
second induced from the coefficient homomorphism.  The profinite
Shapiro isomorphism \cite[Proposition 10]{Serre(1994)}
\begin{equation*}
  \hat{\phi} \colon
  H^k(\hat{G}^p,\coind_{\hat{U}^p}^{\hat{G}^p}\integers/p)\to H^k(\hat{U}^p,\integers/p)
\end{equation*}
is given in exactly the parallel way. Naturality implies that we get a 
commutative diagram
\begin{equation*}
  \begin{CD}
    H^k(\hat{G}^p,\coind_{\hat{U}^p}^{\hat{G}^p}\integers/p)
    @>{\hat{\phi}}>{\iso}> H^k(\hat{U}^p,\integers/p)\\
    @VV{\alpha_G}V @VV{\alpha_U}V\\
    H^k(G,\coind_{U}^G\integers/p) @>{\phi}>{\iso}> H^k(U,\integers/p).
  \end{CD}
\end{equation*}
Because $G$ is cohomologically $p$-complete and
$\coind_{\hat{U}^p}^{\hat{G}^p}\integers/p$ is a finite discrete
$\hat{G}^p$-module, by Lemma \ref{lemma:coeffiso}, $\alpha_G$ is an
isomorphism. It follows that
\begin{equation*}
  \alpha_U\colon H^k(\hat{U}^p,\integers/p)\to 
H^k(U,\integers/p)
\end{equation*}
is an isomorphism, as we had to prove.
\end{proof}

\begin{lemma}\label{lem:trivial_p_power1}
  Assume that $V$ is a $\integers/p$-vector space and $A\colon V\to V$ is a 
 unipotent automorphism of $V$, i.e.~$(A-\id_V)^n=0$ for some $n\ge
  0$. Then, $A^{p^k}=\id_V$ if $p^k\ge n$.
\end{lemma}
\begin{proof}
  If $p^k\ge n$, then 
  \begin{equation*}
    0= (A-\id_V)^{p^k} = \sum_{j=0}^{p^k} \binom{p^k}{j} (-1)^j
    A^j=\id_V-A^{p^k},
  \end{equation*}
  since $\binom{p^k}{j}$ is divisible by $p$ for $0<j<p^k$, and $V$ is 
  a $\integers/p$-vector space.
\end{proof}

\begin{definition}\label{def:p_length}
  For a discrete group $H$, we define the \emph{$p$-lower central series}
  inductively by $\gamma_1^p(H):=H$, and for $k>1$
  $\gamma_k^p(H):=\generate{[\gamma_{k-1}^p(H),H],\gamma_{k-1}^p(H)^p}$, the
  subgroup of $H$
  generated by the commutators $[\gamma_{k-1}^p(H),H]$ and all $p$-th powers of
elements in $\gamma_{k-1}^p(H)$. Observe that these are characteristic 
subgroups. Then
$\gamma_1^p(H)/\gamma_2^p(H)=H_1(H,\integers/p)$ and, more
generally, there is a projection $H_1(\gamma_{k-1}(H),\integers/p)\to
\gamma_{k-1}^p(H)/\gamma_k^p(H)$, in particular,
$\gamma_{k-1}^p(H)/\gamma_k^p(H)$ is a $\integers/p$-vector
space for each $k\ge 2$.

If $H$ is a finite group of order $p^k$, then
$\gamma_{k+1}^p(H)=\{1\}$. The smallest $l\in\naturals$ such that
$\gamma_{l+1}^p(H)=\{1\}$ is called the \emph{nilpotent $p$-length} of 
$H$. If $l$ is the nilpotent $p$-length, then $\gamma_l^p(H)$ is a
central subgroup of $H$ and each element has order $p$ (or $1$).
\end{definition}

\begin{lemma}\label{lem:trivial_p_power2}
  Assume that $H$ is a finite $p$-group of nilpotent $p$-length $k$ and
  $\alpha\colon H\to
  H$ is an automorphism such that the induced action on 
 $H_1(H,\integers/p)$ is trivial. Then $\alpha^{p^{k-1}}=\id_H$.
\end{lemma}
\begin{proof}
  This is done by induction on the nilpotent $p$-length $k$. If $k=1$
  then $H\iso H_1(H,\integers/p)$, and the assertion follows by
  assumption.

  If $k>1$, consider the group $H':= H/\gamma_k^p(H)$. This is a
  finite $p$-group of nilpotent $p$-length $k-1$. Moreover, the
  projection $H\to H'$ induces an isomorphism $H_1(H,\integers/p)\to
  H_1(H',\integers/p)$. Let $\overline\alpha\colon H'\to H'$ be the
  automorphism induced by $\alpha$. By induction,
  $(\overline\alpha)^{p^{k-2}}=\id_{H'}$. 

 Set $\beta:=\alpha^{p^{k-2}}$. We
  show $\beta^p=\id_H$, which immediately implies the
  result.

  Since $\beta$ acts
  trivially on $H/\gamma^p_{k}(H)$, for each $x\in H$ we have
  $\beta(x)=x\phi(x)$ with $\phi(x)\in \gamma_{k}^p(H)$. Now 
  \begin{equation*}
    xy\phi(xy)=\beta(xy)=\beta(x)\beta(y)=x\phi(x)y\phi(y)=xy\phi(x)\phi(y),
  \end{equation*}
  using the fact that $\gamma^p_{k}(H)$ is central in $H$. Therefore
  $\phi\colon H\to \gamma_{n-1}(H)$ is a group homomorphism. Since
  $\gamma_{k}^p(H)$ is an
  abelian $\integers/p$-vector space and lies in the kernel of the
  projection $H\to H_1(H,\integers/p)$ (here we use $k>1$), the
  restriction of $\phi$ to $\gamma_{k}^p(H)$ is trivial, and therefore
  the restriction
  of $\beta$ to $\gamma_{k}^p(H)$ is the identity. It follows that
  $\beta^r(x)=x\cdot (\phi(x))^r$. In particular
  $\beta^p(x)=x$, since $\phi(x)$ like each element of $\gamma_k^p(H)$
  has order $p$.
\end{proof}

\begin{lemma}\label{lem:baby_extension}
Let $H$ be a central subgroup of the group $G$ such that $H$ is a
finite $p$-group and $G/H$ is cohomologically $p$-complete.
Then $G$ contains a normal subgroup
  $N$ such that $G/N$ is a finite $p$-group, and $H\cap N=\{1\}$.
\end{lemma}
\begin{proof}
  We prove this by induction on the nilpotent $p$-length $k$ of
  $H$. Set $Q = G/H$ and assume first that $k=1$. Then $H=H_1(H,\integers/p)\iso
  H^1(H,\integers/p)$. By Lemma \ref{lemma:shortexactcompletion}, the
  exact sequence $1\to H\to G\to Q\to 1$ maps to the exact sequence of 
  pro-$p$-completion $1\to\overline H\to\hat{G}^p\to\hat{Q}^p\to 1$,
  where, since $H$ is a finite $p$-group, $\overline H$ is a quotient
  of $H$. The spectral sequences for the cohomology of such
  extensions give rise to exact sequences in low degrees, and because 
  of naturality we get the following diagram of exact sequences
{\small  \begin{equation*}
    \begin{CD}
       H^1(\hat{Q}^p,\integers/p) @>>> H^1(\hat{G}^p,\integers/p)
      @>>> H^1(\overline{H},\integers/p) @>>>
      H^2(\hat{Q}^p,\integers/p)@>>> H^2(\hat{G}^p,\integers/p)\\
       @VV{\iso}V @VV{\iso}V @VV{\beta}V @VV{\iso}V
      @VV{\alpha}V\\
       H^1({Q},\integers/p) @>>> H^1({G},\integers/p)
      @>>> H^1({H},\integers/p) @>>>
      H^2(Q,\integers/p)@>>> H^2(G,\integers/p)        .
    \end{CD}
  \end{equation*}
}
  Note that in general, the middle term is $H^1(H,\integers/p)^G$, the
  $G$-fixed set of $H^1(H,\integers/p)$ (and similarly $H^1(\overline
  H,\integers/p)^{\hat{G}^p}$). In our case, by assumption, the action
  of $G$ on $H$ is trivial, and this implies that $G$ acts trivially
  on $H^1(H,\integers/p)$ and
  $\hat{G}^p$ acts trivially on $H^1(\overline H,\integers/p)$, so
  that we get exactly the above diagram.
  
  By Lemma \ref{lem:low_completeness}, the first two homomorphism are
  isomorphisms, and $\beta$ is injective. $Q$ being cohomologically
  $p$-complete explains why the fourth homomorphism is an
  isomorphism. The $5$-lemma now implies that $\beta$ is an
  isomorphism.

  Since, by assumption, $H\iso H^1(H,\integers/p)$ and $\overline H$
  is a quotient of $H$, we also have $\overline H\iso H^1(\overline
  H,\integers/p) \underset{\beta}{\xrightarrow{\iso}} H^1(H,\integers/p)\iso
  H$. By Lemma \ref{lemma:shortexactcompletion} this means that for
  each normal subgroup $U$ of $H$ of $p$-power index, in particular for
  $\{1\}$, we find a normal subgroup $V$ of $G$ such that $G/V$ is a
  finite $p$-group and such that $V\cap H\subset U$. This proves the
  statement in the case $k=1$.

  If the nilpotent $p$-length $k$ is bigger than $1$, consider
  $L:=\gamma_2^p(H)= \ker(H\to H_1(H,\integers/p))$. Since this is a
  characteristic subgroup of $H$, we get an exact sequence
  \begin{equation*}
    1\to H/L\to G/L\to Q\to 1,
  \end{equation*}
  with $H_1(H/L,\integers/p)\iso H/L$, and still $G/L$ acting
  trivially on $H_1(H/L,\integers/p)\iso H_1(H,\integers/p)$. As we
  have just proved, we then find a normal subgroup $U$ of $G/L$ which
  has finite $p$-power index in $G/L$ and such that $U\cap
  H/L=\{1\}$. Let $\pi\colon G\to G/L$ be the canonical projection,
  and $V:=\pi^{-1}(U)$ the inverse image of $U$ in $G$. Then $G/V$ is
  a finite $p$-group, and $V\cap H\subset L=\gamma_2^p(H)$. Set
  $Q':=VH/H$. This is a normal subgroup of $Q$, and $Q/Q'$ is a finite 
  $p$-group. By Proposition \ref{prop:completeness_and_subgroups},
  $Q'$ is cohomologically $p$-complete. We get an exact sequence 
  \begin{equation*}
    1\to L\to V\to Q'\to 1. 
  \end{equation*}
  Of course, still $V$ acts trivially on $L$. Moreover, the nilpotent
  $p$-length of $L=\gamma_2^p(H)$ is strictly smaller than $k$, the
  nilpotent $p$-length of $H$. By induction, we find therefore a
  normal subgroup $W$ of $V$ such that $V/W$ is a finite $p$-group and 
  such that $L\cap W=\{1\}$. By Lemma \ref{lem:subgroupintersect}, we
  can replace $W$ by a smaller normal subgroup $W'$ such that still
  $V/W'$ is a finite $p$-group, such that of course also still
  $L\cap W'=\{1\}$, and such that $W'$ is normal in $G$. It follows
  that $G/W'$ is a finite $p$-group. Since $W'\subgroup V$ and $V\cap
  H\subgroup L$, we finally get from $W'\cap L=\{1\}$ that $W'\cap
  H=\{1\}$, as desired.
\end{proof}
\begin{remark}
Though the condition that $H$ is central in $G$ is rather special,
Lemma \ref{lem:baby_extension} is still an important first step
  in our proof of Theorem
  \ref{theo:extensions_and_cohomological_completeness}.
\end{remark}

\begin{lemma}\label{lem:toddler_extension}
    Assume that $1\to H\to G\to Q\to 1$ is an exact sequence of groups, $H$
  is a finite $p$-group, each element of $G$ acts nilpotently on
  $H_1(H,\integers/p)$, and $Q$ is
  cohomologically $p$-complete. Then $G$ contains a normal subgroup
  $V$ such that $G/V$ is a finite $p$-group, and $H\cap V=\{1\}$.
\end{lemma}
\begin{proof}
  The action of $G$ on $H$ is a homomorphism from $G$ to the finite
  group $\Aut(H)$. Let $N$ be the kernel of this homomorphism. We claim 
  that $G/N$ is a finite $p$-group. Since $G/N$ is finite, it suffices 
  to show that the order of each element of $G/N$ is a power of
  $p$. Fix $g\in G$. Lemma \ref{lem:trivial_p_power1} and Lemma
  \ref{lem:trivial_p_power2} together imply that $g^{p^k}$ acts
  trivially on $H$ if $k$ is sufficiently big, i.e.~$g^{p^k}\in
  N$. Since $g$ was arbitrary, the order of each element of $G/N$ is a 
  power of $p$.

  Set $Q':= NH/H$. This is a normal subgroup of $Q$, and $Q/Q'$ is a
  finite $p$-group. By Proposition
  \ref{prop:completeness_and_subgroups}, $Q'$ is cohomologically
  $p$-complete. Moreover, in the exact sequence
$    1\to H\cap N\to N \to Q' \to 1$, by construction $N$ acts
trivially on $H\cap N$. By Lemma \ref{lem:baby_extension}, we can find a
normal subgroup $U$ of $N$ such that $N/U$ is a finite $p$-group and
such that $H\cap U= (H\cap N)\cap U =\{1\}$. Moreover by Lemma
\ref{lem:subgroupintersect}, we can replace $U$ by a smaller subgroup
$V$ such that $V \lhd G$ and $N/V$ is a finite $p$-group.
Then $G/V$ is also a finite $p$-group, which finishes the proof
because obviously $H \cap V = \{1\}$.
\end{proof}

\begin{proposition}\label{prop:big_extension}
  Assume that $1\to H\to G\to Q\to 1$ is an exact sequence of groups, $H$
  is finitely generated, each element of $G$ acts unipotently on
  $H_1(H,\integers/p)$, and $Q$ is
  cohomologically $p$-complete. Then we get an exact sequence of
  pro-$p$ completions
  \begin{equation*}
    1\to \hat{H}^p\to \hat{G}^p\to \hat{Q}^p\to 1.
  \end{equation*}
\end{proposition}
\begin{proof}
  By Lemma \ref{lemma:shortexactcompletion}, we only have to prove
  that for each normal subgroup $U$ of $H$ such that $H/U$ is a finite 
  $p$-group we can find a normal subgroup $V$ of $G$ such that $G/V$ is a
  finite $p$-group and such that $V\cap H\subgroup U$. Fix such a
  $U$. Since $H$ is finitely generated, by Lemma
  \ref{lem:subgroupintersect}, we can replace $U$ by a smaller characteristic
  subgroup $V$ such that $V\subgroup U$ and $H/V$ is a finite
  $p$-group. We obtain the exact sequence $1\to H/V\to G/V \to Q\to
  1$. Since $H_1(H/V,\integers/p)$ is a quotient of
  $H_1(H,\integers/p)$, each element of $G/V$ acts unipotently on
  $H_1(H,\integers/p)$ (without loss of generality we can even assume
  that $V\subgroup \gamma_2^p(H)$, and then $H_1(H,\integers/p)\iso
  H_1(H/V,\integers/p)$). Moreover, $H/V$ is a finite $p$-group. By
  Lemma \ref{lem:toddler_extension}, we find a normal subgroup $W$ of
  $G/V$ of finite $p$-power index such that $W\cap H/V=\{1\}$. Let
  $\tilde W$ be the inverse image of $W$ in $G$. Then $\tilde W\cap
  H\subgroup V\subgroup U$, and $G/\tilde W$ is a finite $p$-group.
\end{proof}

\begin{proof}[Proof of Theorem
  \ref{theo:extensions_and_cohomological_completeness}]
  Because of Proposition 
  \ref{prop:big_extension}, the assumptions imply  that we have the
  exact sequence
  \begin{equation*}
    1\to \hat{H}^p \to \hat{G}^p\to \hat{Q}^p\to 1.
  \end{equation*}
  As in the proof of Proposition \ref{tfcompletion}, we have two
  compatible spectral sequences with 
 $E_2^{l,q}$-terms
  \begin{equation*}
 H^l(\hat{Q}^p,H^q({\hat H}^p,\integers/p)) 
\qquad\text{or}\qquad
    H^l(Q,H^q(H,\integers/p))
  \end{equation*}
  converging to
 $H^{l+q}({\hat G}^p,\integers/p)$
 or
 $ H^{l+q}(G,\integers/p)$,
 respectively. The
  assumption on $H$ means that the natural map between Galois and ordinary
  cohomology induces an isomorphism between the coefficients in the
  $E_2$-terms, which, again by assumption, are finite
  $\hat{Q}^p$-modules. Because of the assumption on $H$ and Lemma
  \ref{lemma:coeffiso}, the natural map therefore induces an
  isomorphism on the $E_2$-term, and
  consequently the isomorphism between $H^*({\hat G}^p,\integers/p)$
  and $H^*(G,\integers/p)$ we had to establish.

  If $H^l(\hat{Q}^p,\integers/p)$ is finite for each $l$, the same is true for 
  each $H^l(\hat{Q}^p,A)$, for an arbitrary finite discrete
  $\hat{Q}^p$-module $A$. This means that every $E_2^{k,l}$ is
  finite. Because the spectral sequence is a first quadrant spectral
  sequence, the same is then true for the groups
  $H^{k+l}(\hat{Q}^p,\integers/p)$ the spectral sequence is converging 
  to.
\end{proof}

\begin{lemma}\label{lem:free_prod_pgroup_complete}
  Let $P$ and $Q$ be two finite $p$-groups. Then the free product $P*Q$ is
  cohomologically $p$-complete. In particular,
  \begin{equation*}
H^k(\widehat{P*Q}^p,\integers/p)=H^k({\hat P}^p,\integers/p)\oplus
  H^k({\hat Q}^p,\integers/p)\quad\text{for } k\ge 1.
\end{equation*}
\end{lemma}
\begin{proof}
  By \cite[Lemma 5.5]{Gruenberg(1957)} we have an exact sequence
  \begin{equation*}
    1\to F\to P*Q\to P\times Q\to 1
  \end{equation*}
  where $F$ is a (finitely generated, since $P$ and $Q$ are finite)
  free group. 
  Since $P\times Q$ is a finite $p$-group and $F$ is cohomologically
  $p$-complete by Example \ref{ex:free_completeness}, Proposition
  \ref{tfcompletion} implies that $P*Q$ is cohomologically
  $p$-complete, as well. 

  Since $H^k(P*Q,\integers/p)\iso H^k(P,\integers/p)\oplus
  H^k(Q,\integers/p)$ for $k\ge 1$, the corresponding statement for
  the pro-$p$ completions follows from cohomological $p$-completeness.
\end{proof}

We note the following proposition for later reference in Section
\ref{sec:examples}. 
\begin{proposition}\label{prop:free_product_complete}
  Let $G_1$ and $G_2$ be two discrete groups. Set $G:=G_1*G_2$. The
  inclusions $G_k\into G$ induce an isomorphism
  \begin{equation}\label{eq:1}
    H^k({\hat G}^p,\integers/p) \xrightarrow{\iso} H^k({\hat
    G_1}^p,\integers/p)\oplus H^k({\hat
    G_2}^p,\integers/p)\qquad\text{for }k\ge 1.
\end{equation}
  In particular, $G_1$ and $G_2$ both are cohomologically
  $p$-complete if and only if the same is true for $G$.
\end{proposition}
\begin{proof}
  The completion ${\hat G}^p$ is the inverse limit of the inverse
  system of all quotients $G/U$ where $U$ is a normal subgroup of $G$
  of finite $p$-power index. Set $U_1:=G_1\cap U$ and $U_2:=G_2\cap
  U$. We get a projection $(G_1/U_1)* (G_2/U_2)\onto G/U$. On the other
  hand, $(G_1/U_1)*(G_2/U_2)$ is a quotient of $G$, so every finite
  $p$-group quotient of $(G_1/U_1)*(G_2/U_2)$ is a finite $p$-group
  quotient of $G$. It follows that we can view the inverse system
  which defines ${\hat G}^p$ as the inverse system (over $U$) of the inverse
  systems which define the pro-$p$ completion $X_U$ of
  $(G_1/U_1)*(G_2/U_2)$. The Galois
  cohomology of ${\hat G}^p$ by definition is the direct limit of the
  corresponding direct system of cohomology groups. With our
  reinterpretation, for $k\ge 1$,
  \begin{equation*}
    \begin{split}
      H^k({\hat G}^p,\integers/p) & = \lim_{G/U}
      H^k(X_U,\integers/p)\\
      & = \lim_{G/U} H^k(G_1/U_1,\integers/p)\oplus
      H^k(G_2/U_2,\integers/p).
  \end{split}
\end{equation*}
  For the last equality, we use Lemma \ref{lem:free_prod_pgroup_complete}.

  We take the limit along the system of finite $p$-group quotients of
  $G$. It remains to check for each pair of normal subgroups
  $U_1\subgroup G_1$ and $U_2\subgroup G_2$ of $p$-power index that
  there is a corresponding normal subgroup $U$ of $G$ with $G_1\cap
  U=U_1$ and $G_2\cap U=U_2$. But for such $U_1$ and $U_2$, the kernel
  $U$ of the projection $G\to
  (G_1/U_1)\times (G_2/U_2)$ fulfills $G_1\cap U=U_1$ and $G_2\cap
  U=U_2$. 
Therefore the limit along the system of finite $p$-group quotients of
  $G$ is the same as the limit
  along the system induced from the finite $p$-group quotients of
  $G_1$ and of $G_2$, and consequently for $k\ge 1$
  \begin{equation*}
    \begin{split}
      H^k({\hat G}^p,\integers/p) & = \lim_{G/U}
      H^k(G_1/U_1,\integers/p)\oplus
      H^k(G_2/U_2,\integers/p)\\
     & = H^k({\hat G_1}^p,\integers/p)\oplus H^k({\hat
        G_2}^p,\integers/p).
  \end{split}
\end{equation*}
  
  We get a commutative diagram
  \begin{equation*}
    \begin{CD}
      H^k({\hat G}^p,\integers/p) @>>> H^k({\hat
        G_1}^p,\integers/p)\oplus H^k({\hat G_2}^p,\integers/p)\\
      @V{\alpha_k}VV @V{\beta_k}VV\\
      H^k(G,\integers/p) @>>> H^k(G_1,\integers/p)\oplus H^k(G_2,\integers/p)
    \end{CD}
  \end{equation*}
  where the horizontal arrows are induced from the inclusions of the
  factors $G_1$ and $G_2$ in $G$, and, for $k\ge 1$, both are
  isomorphisms (this is a classical fact in the discrete case, and we
  just checked it for the pro-$p$ completions). Therefore, $\alpha_k$
  is an isomorphism if and only if
  $\beta_k$ is an isomorphism. Since, by definition, $\alpha_0$ and
  $\beta_0$ always are isomorphisms, $G$ is cohomologically
  $p$-complete if and only if $G_1$ and $G_2$ both are cohomologically
  $p$-complete. 
\end{proof}

The following theorem is the technical heart of our program to find
torsion-free quotients of our torsion-free finite extension, and for
more general
finite extensions to find quotients with as little 
torsion as possible.

\begin{theorem}\label{theorem:splitback}
  Assume that we have an exact sequence
  \begin{equation*}
    1\to {H} \to G \to Q \to 1,
  \end{equation*}
  where $Q$ is a finite $p$-group and  ${H}$ is cohomologically $p$-complete
  and has a finite classifying space $B{H}$. Assume that the induced
  exact sequence
  \begin{equation*}
    1\to {\hat {H}}^p\to {\hat G}^p \to Q \to 1
  \end{equation*}
  splits. Then the original sequence also splits, in particular we
  have an embedding $Q\into G$.
\end{theorem}

The idea of the proof is that the splitting of the pro-$p$ completions
induces a splitting of cohomology. But the cohomology for the
completions and the original groups coincide (since we assume that ${H}$
is cohomologically $p$-complete), so that the original map splits after
taking cohomology. With a spectral sequence argument, the same follows
on the level of stable cohomotopy, and this induces the splitting on
the group level.

To make this precise, some preparation is necessary.

The next theorem with proof was communicated to us by Alejandro Adem. We thank
him for his very prompt and kind answer to our questions.

\begin{theorem}\textbf{\textup{(Adem)}}\label{Adem}
Let ${H}$ denote a discrete group of finite cohomological dimension.
Suppose that we are given an extension 
\[
1\to {H}\to G\xrightarrow{f} Q\to 1
\]
where $Q$ is a finite $p$-group. For a space $Y$, let $\tilde\pi^0_S(Y)$
denote its degree zero reduced stable cohomotopy group.

The extension above splits if and only if the epimorphism
$f\colon G \to Q$ induces an injection $\tilde
f^*\colon \tilde\pi_S^0(BQ)\to\tilde\pi_S^0(BG)$.
\end{theorem}
\begin{proof}
Clearly, if the extension splits, then the map induced on reduced stable 
cohomotopy is injective. Now let us assume that this map is
injective. Then the map $f^*$ induced in unreduced cohomotopy is also
injective, since we only add the summand $\id\colon \pi_S^0(pt)\to
\pi_S^0(pt)$.

By a construction due to Serre \cite[Theorem 3.1 and following
exercise]{Brown(1982)} we can find a finite dimensional
CW-complex $X$ with an action of $G$ and such that 
$BG\simeq X/{H}\times_QEQ$. In particular, we can identify
the map $f^*$ with the map induced in stable cohomotopy by the fibration
$X/{H}\times_QEQ\to BQ$. 

The main result in \cite{Jackowski(1988)} implies that, if $f^*$ is
injective, then
the action of $Q$ on $X/{H}$ has a fixed point. In particular, the
fibration above has a \emph{section} and hence the group extension
must split.
\end{proof}

We want to apply this result to the extension given in Theorem
\ref{theorem:splitback}. As
an intermediate step we use stable cohomotopy of the pro-$p$
completions and localization of these cohomotopy groups at the prime
$p$.

\subsection{Finiteness conditions for classifying spaces}
At several points, we will have to assume finiteness conditions on the
classifying spaces. We collect a few results about this here. Since we
are only interested in (co)homology, we can do all of
our constructions up to homotopy.

\begin{definition}
  A CW-complex is called finite, if it consists of finitely many
  cells. It is called of \emph{finite type}, if each finite dimensional
  subcomplex is finite.
\end{definition}

\begin{lemma}\label{lemma:fibrationtype}
  If $F\to E\to B$ is a fibration such that $F$ and $B$ are homotopy
  equivalent to CW-complexes of finite type then the same is
  true for $E$.

    Let $1\to H\to G\to P\to 1$ be an exact sequence of groups,  where
  $BH$ and $BP$ are
  CW-complexes of finite type. Then the same is true for $BG$. 
\end{lemma}
\begin{proof}
  The first statement follows from \cite[1.3 and
  1.8]{Lueck-Schick-Thielmann(1998)}, and also
  from \cite[Chapter 14]{Lueck(1989)}. An explicit account can be
  found in \cite[2.2.1]{Bratzler(1997)}.

Then apply this to the fibration $BH\to BG\to BP$.
\end{proof}

\subsection{Generalized Galois cohomology of profinite groups}

\begin{definition}\label{def:genhompgroup}
  Assume that $h^*$ is a generalized cohomology theory and $\profinQ$ is a
  profinite group. We set
  \begin{equation*}
    h^*(\profinQ):= \dirlim_{\abs{\profinQ/N}<\infty} h^*(\profinQ/N)
  \end{equation*}
  where $h^*(\profinQ/N):=h^*(B[\profinQ/N]) $ is the cohomology of the classifying
  space of the finite
  quotient  $\profinQ/N$.
\end{definition}

\subsubsection{Filtration by skeleta}

Let $X$ be a CW-complex and $h^*$ an arbitrary
generalized 
cohomology theory. Let $X^{0}\subset X^{1}\subset\dots\subset X$
be the filtration of $X$ by its
skeleta. Let $F^*_k(X)$ for $k\in\naturals$ be the kernel of the induced map
$h^*(X)\to
h^*(X^{(k)})$. Set $F^k_0(X):= h^k(X)$. This way, we get a filtration
\begin{equation*}
  h^*(X)= F^*_0(X)\supset
  F^*_1(X)\supset\dots\supset \{0\}.
\end{equation*}
If $X$ is finite dimensional then this stabilizes after finitely many
steps. More precisely, $F^k_{\dim X}(X)=\{0\}$, since $X^{\dim X}=X$.
If $X$ is not finite dimensional the sequence does not
necessarily stabilize. It
is not even clear whether $\bigcap_{k\in\naturals} F^*_k(X)=\{0\}$. We
will later study
conditions under which the last assertion is true.

Note, however, that every continuous map (being homotopic to a cellular
map) induces maps which preserve these filtrations.
In particular, by passing to direct limits, we get, for a profinite
group $\profinQ$, a filtration
\begin{equation*}
  h^*(\profinQ) =  F_0^*(\profinQ)\supset F_1^*(\profinQ)\supset\cdots .
\end{equation*}

\begin{lemma}\label{lem:generalsplitting}
  Definition \ref{def:genhompgroup} defines a functor from the
  category of profinite groups to filtered graded abelian groups. In
  particular, a split projection $G\to\profinQ$ induces a split injection
  $h^*(\profinQ)\to h^*(G)$. It also induces split injections of the
  filtration spaces $F^*_k(\profinQ)\to F^*_k(G)$ and of the quotients
  $F^*_l(\profinQ)/F^*_k(\profinQ)\to F^*_l(G)/F^*_k(G)$ if $l\le k$.

  If for a prime number $p$ the group $\profinQ$ is the
  pro-$p$ completion of a discrete group $G$,
  then the inclusion induces a map $h^*(\profinQ)\to h^*(G)$.
\end{lemma}
\begin{proof}
  A map between profinite groups is the same as a map between the
  inverse system of finite quotients, giving maps between the direct
  systems of cohomology groups. This is functorial. Because of the
  splitting, we can also pass to the quotients $F_l/F_k$.

  The last statement, about the pro-$p$ completion, follows because we
  get a consistent system of maps to all the finite quotients of $\profinQ$,
  which induce a corresponding map from the direct system of
  cohomology groups to $h^*(G)$.
\end{proof}

\subsection{Localization}

We collect a few well known results about localization
in the following proposition (compare e.g.~\cite[Chapter
15.4]{Dummit-Foote(1999)} and \cite[p.~67 in II.2.4]{Bourbaki(1989)}), which we are going to use later.

\begin{proposition}\label{prop:localizationproperties}
  Let $p$ be a prime number and
  $\integers_{(p)}:=\{a/b\in\rationals \mid p \dividesnot b\}$ the
  localization of $\integers$ at $p$. For an arbitrary abelian group
  $M$, $M_{(p)}:=M\tensor_{\integers}\integers_{(p)}$ is the
  localization of $M$ at $p$. Then
  \begin{itemize}
  \item Localization at $p$ is an exact functor.
  \item The kernel of the map $M\to M_{(p)}\colon m\mapsto m\tensor 1$
    consists of those elements in $M$ which are of finite order not
    divisible by $p$. In particular, if $M$ is a torsion module, this
    map is the projection onto the $p$-primary component of $M$. 
  \item Localization commutes with tensor products.
  \item Localization commutes with direct limits.
  \end{itemize}
\end{proposition}

\subsection{Spectral sequences}

A (cohomological) spectral sequence $(E^{*,*}_*)$ consists of abelian
groups 
\begin{equation*}
E^{s,t}_r \quad\text{for } r\ge 2,\; s,t\in\integers
\end{equation*} 
with homomorphisms, the so called
differentials,
\begin{equation*}
d^{s,t}_r\colon E^{s,t}_r\to E^{s+r,t-r+1}_r,\qquad\text{which satisfy }
d_r^{s,t}\circ 
d_r^{s-r,t+r-1}=0,
\end{equation*}
and with isomorphisms
\begin{equation*}
\Phi^{s,t}_r\colon \ker(d_r^{s,t})/\im(d_r^{s-r,t+r-1})
\xrightarrow{\iso}  E^{s,t}_{r+1}\qquad\text{for } r\ge
  2,\; s,t\in\integers.
\end{equation*}
In other words, each page $E_r^{*,*}$ ($r\ge
  2$) consists of chain complexes, and $E^{*,*}_{r+1}$ is obtained as
  the cohomology of $(E_r^{*,*},d^{*,*}_r)$.

  A map of spectral sequences $f_*^{*,*}\colon E^{*,*}_*\to F^{*,*}_*$
  consists of homomorphisms $f_r^{s,t}\colon E^{s,t}_r\to F^{s,t}_r$ which
  are compatible with the differentials $d^{*,*}_*$ and the
  isomorphisms $\Phi^{*,*}_*$.

  If $E^{s,t}_2=0$ whenever $s<0$, then for all $r\ge 2$ we have
  $E^{s,t}_r=0$ if $s<0$. Consequently, for $r>s$ the domain
  of $d_r^{s-r,t+r-1}$ is zero and hence $\im(d^{s-r,t+r-1}_r) =0$ and
  therefore
  \begin{equation*}
E^{s,t}_{r+1}=\ker(d_r^{s,t})\subset E^{s,t}_r\qquad\text{for }r>s .
\end{equation*}
In
  this case, we define $E_\infty^{s,t}:= \bigcap_{r>s}E_r^{s,t}$.

  We denote a spectral sequence $(E^{*,*}_*)
$ a \emph{fourth quadrant}
  spectral sequence, if $E^{s,t}_2=0$ for $s<0$ or $t>0$. Note that
  for a fourth quadrant spectral sequence, $E_\infty^{s,t}$ is defined
  for all $s,t\in\integers$.

  Let $h^p = F^p_0\supset F^p_1\superset F^p_2\superset\cdots$ be a
  collection of graded groups, $p\in\integers$. We say a spectral
  sequence $(E^{*,*}_*)$ for which $E_\infty^{s,t}$ is defined (e.g.~a
  fourth quadrant spectral sequence) \emph{converges} to $h^*$, if
  \begin{equation*}
E^{s,t}_\infty\iso F^{s+t}_s/F^{s+t}_{s+1}\qquad\text{for all }
  s,t\in\integers.
\end{equation*}

  If $F^p_k=0$ for $k$ sufficiently large, this determines (up to
  extension problems) $h^p$. This is not necessarily the case if only
  $\bigcap_{k\in\naturals} F^p_k =0$.

\begin{lemma}\label{lem:finspecseq}
  Assume that $E^{s,t}_r$ is a fourth quadrant spectral sequence converging
  to the graded groups $h^{s+t}=H^{s+t}_0\superset
  H^{s+t}_1\superset\cdots$.
  Let $\mathcal{P}$ be a property of abelian groups which is preserved
  under passage to subgroups, quotient groups, and under extensions
  (i.e.~if in an exact sequence $1\to U\to G\to Q\to 1$ the groups $U$
  and $Q$ have property $\mathcal{P}$, then also $G$ has property
  $\mathcal{P}$). Examples of such properties are the property of being
  a finite group and the property of being a $p$-group ($p$ a prime number).
  \begin{enumerate}
  \item If, for some $s,t\in\integers$, $E^{s,t}_2$ has property
    $\mathcal{P}$, then $E^{s,t}_\infty$ also has property
    $\mathcal{P}$. 
\item If $s,t\in\integers$ with $s+t=n$ implies that
    $E^{s,t}_2$ has property $\mathcal{P}$, then $h^{n}/H^{n}_k$ has property $\mathcal{P}$ for
    every $k\in\naturals$.
  \item If, in addition, for the given $n\in\integers$ there is some
    $k\in\naturals$ such that $H^{n}_k=0$, then $h^n$ has property
    $\mathcal{P}$.
  \item\label{item:injectiv} Assume that $F^{s,t}_r$ is another fourth quadrant spectral
    sequence converging to the graded group
    $g^{s+t}=G^{s+t}_0\superset G^{s+t}_1\superset \cdots$. Assume that we
    have a map of spectral sequences $f^{s,t}_r\colon E^{s,t}_r\to
    F^{s,t}_r$ such that $f^{s,t}_\infty:E^{s,t}_\infty\to
    F^{s,t}_\infty$ is injective for every $s,t\in\integers$. Suppose,
    moreover, that this map is compatible with a corresponding map of
    graded groups $h^*\to g^*$. Then this induces injective maps
    \begin{equation*}
      h^n/H^n_k \into g^n/G^n_k\qquad\text{for all } k\in\naturals,\;
      n\in\integers.
    \end{equation*}
  \end{enumerate}
\end{lemma}
\begin{proof}
  First observe that by induction on $r$ every $E^{s,t}_r$ ---as a
  sub-quotient of $E^{s,t}_{r-1}$--- has property $\mathcal{P}$, and the
  same then is also true for $E^{s,t}_\infty\subset E^{s,t}_s$. Using
  the extensions
  \begin{equation}\label{eq:itextens}
    1\to    H^{s+t}_{s}/H^{s+t}_{s+1} \to h^{s+t}/H^{s+t}_{s+1}\to
    h^{s+t}/H^{s+t}_{s}  
 \to 1,
  \end{equation}
  and the isomorphism $H^{s+t}_s/H^{s+t}_{s+1}\iso E^{s,t}_\infty$,
  it follows by induction that $h^{s+t}/H^{s+t}_k$ has property $\mathcal{P}$
  for every $k\in\naturals$. If $H^{s+t}_k=0$, then
  $h^{s+t}=h^{s+t}/H^{s+t}_k$, so that under this assumption also $h^{s+t}$ has
  property $\mathcal{P}$.

  Given a map between short exact sequences of abelian groups
  \begin{equation*}
    \begin{CD}
      0 @>>> A_1@>>> A_2@>>> A_3 @>>> 0\\
      && @VV{f_1}V  @VV{f_2}V  @VV{f_3}V\\
      0 @>>> B_1 @>>> B_2 @>>> B_3 @>>> 0,
    \end{CD}
  \end{equation*}
  an easy diagram chase shows that $f_2$ is injective if 
  both $f_1$ and $f_3$ are injective. The injectivity in
  \ref{item:injectiv} of
  $h^n/H^n_k\into g^n/G^n_k$ follows inductively from this principle
  applied to the exact sequences \eqref{eq:itextens} and their
  counterparts for $g^*$.
\end{proof}

\begin{lemma}\label{lem:specseq}
  Let $E_r^{s,t}$ be a cohomological spectral sequence. 

  \begin{enumerate}
  \item \label{item:finiteprop} Given $r,s,t$, there is a finite set $V(r,s,t)\subset
    \integers\times \integers$ such that $E_r^{s,t}$ is determined by
    $E_2^{x,y}$ for $(x,y)\in V(r,s,t)$. More precisely, if
    $f_r^{s,t}\colon E_r^{s,t}\to F_r^{s,t}$ is a map of spectral
    sequences such that $f_2^{x,y}\colon E_2^{x,y}\to F_2^{x,y}$ is an
    isomorphism for all $(x,y)\in V(s,t,r)$, then $f_r^{s,t}\colon
    E_r^{s,t}\to F_r^{s,t}$ is also an isomorphism.
  \item \label{item:4quadfiniteprop}
    If we restrict to the category of fourth quadrant spectral
    sequences and if $s+t\le 0$, then we can
    choose $V(r,s,t)$ as in \ref{item:finiteprop} such that
    $(x,0)\notin V(r,s,t)$ if
    $x\ge 1$.
\end{enumerate}
\end{lemma}
\begin{proof}
  We prove \ref{item:finiteprop} by induction. The statement is trivial for
  $r=2$. Given $r>2$,
  $E^{s,t}_r$ depends only on the three groups
  $E^{s-r+1,t+r-2}_{r-1}$, $E^{s,t}_{r-1}$ and
  $E^{s+r-1,t+r-2}_{r-1}$. By induction, each of these depends only on
  finite set $V_1, V_2, V_3\subset \integers^2$, and we can choose $V$
  to be the union of these.

To prove \ref{item:4quadfiniteprop},
  assume from now on that all spectral sequences are fourth quadrant
  spectral sequences. Define $Z_2:=\{(x,0)\mid x\ge
  1\}\subset \integers\times\integers$ and inductively 
  \begin{equation*}
Z_{r+1}:= \bigl(Z_{r}\cup d_r(Z_r)\cup
  d_r^{-1}(Z_r)\bigr)\cap \bigl(\integers_{\ge 0}\times \integers_{\le 0}\bigr) \qquad
  \text{for } r\ge 2,
\end{equation*}
  with $d_r(x,y):=(x+r,y-r+1)$.
This is
  the ``zone of influence'' of
  $Z_2$: if $(s,t)\notin Z_r$ then $E^{s,t}_r$ does not change if only
  $E_2^{x,y}$ with $(x,y)\in Z_2$ is changed. 

To prove \ref{item:4quadfiniteprop}, it suffices to show that
$(s,t)\notin Z_r$ if $s+t\le
  0$. Assume that this last assertion is wrong and let $r$ be minimal
  such that there is $(s,t)\in Z_r$
  with $s+t\le 0$. It follows from the definition of $Z_r$  that there
   are sequences
   \begin{equation*}
2\le r_1<r_2<\dots<r_n\qquad\text{and}\qquad
\epsilon_1,\dots,\epsilon_n\in\{1,-1\},
\end{equation*}
 and $x\ge 1$ such that for every  $k=1,\dots,n$
 \begin{equation*}
   (x+\epsilon_1r_1+\dots +\epsilon_k r_k,
  (1-r_1)\epsilon_1+ (1-r_2)\epsilon_2+\dots+(1-r_k)\epsilon_k)\in
  Z_{r_k+1},
\end{equation*}
and
  \begin{equation*}
(s,t)=(x+\epsilon_1r_1+\dots +\epsilon_n r_n,
  (1-r_1)\epsilon_1+ (1-r_2)\epsilon_2+\dots+(1-r_n)\epsilon_n).
\end{equation*}
Let $i_1<i_2<\cdots$ be the ordered sequence of indices with
$\epsilon_{i_k}=1$, $j_1<j_2<\cdots$ similarly with $\epsilon_{j_k}=-1$.
The inequality $s+t\le 0$ translates to
\begin{equation}\label{ineq:s+t}
\epsilon_1+\dots+\epsilon_n\le -x <0.
\end{equation}
On the other hand, since we consider spectral sequences which are zero in
the upper half-plane and define $Z_r$ accordingly,
\begin{equation}\label{ineq:4quad}
(1-r_1)\epsilon_1+\dots+(1-r_k)\epsilon_k
\le 0 \qquad\text{for every }k=1,\dots,n.
\end{equation}
This implies that $\epsilon_1=1$, i.e.~$i_1=1$ and, using monotonicity
of the $r_i$, that $\epsilon_2=1$, i.e.~$i_2=1$ and, inductively, that
$i_{k+1}<j_k$ for $k=1,2,\dots$. Consequently, less than half of the
$\epsilon_1,\epsilon_2,\dots,\epsilon_n$ are negative. But then
\eqref{ineq:s+t} is violated. This gives the desired
contradiction and proves \ref{item:4quadfiniteprop}.
\end{proof}

\subsection*{Proof of Theorem \protect\ref{theorem:splitback}}

From now on, we study the generalized homology theory ``(reduced) stable
cohomotopy'' $\pi_S^*$. Consequently, if $X$ is a CW-complex or a
discrete group or a
profinite group, $F^n_r(X)$ always will denote the corresponding
filtration group  $F^n_r(X)\subset \tilde\pi_S^n(X)$.
 
In the situation of Theorem \ref{theorem:splitback}, we will consider
the following
commutative diagram:
\begin{equation}
  \label{eq:bigdiagram}
  \begin{CD}
    \tilde\pi_S^0(Q) @>>>  
    \tilde\pi^0_S(Q)/F^0_k(Q) @>{A}>{\iso}> \left
    ( \tilde\pi_S^0(Q)/F^0_k(Q)\right)_{(p)}     \\
    @VV{\id}V &&  @V{A'}V{\iso}V\\
    \tilde\pi_S^0(Q) @>>> \tilde\pi_S^0(Q)_{(p)} @>>>
    \tilde\pi_S^0(Q)_{(p)}/F_k^0(Q)_{(p)}\\
    @VV{i}V    @VV{i}V    @VV{i}V\\
    \tilde\pi_S^0({\hat G}^p)@>>> \tilde\pi_S^0({\hat G}^p)_{(p)} @>>>
    \tilde\pi_S^0({\hat G}^p)_{(p)}/ F^0_k({\hat G}^p)_{(p)}\\
    @VVV    @VVV  @V{B}VV\\
    \tilde\pi^0_S(G) @>>> \tilde\pi^0_S(G)_{(p)} @>>>
    \tilde\pi^0_S(G)_{(p)}/F^0_k(G)_{(p)}.
  \end{CD}
\end{equation}

Because of Lemma \ref{lem:generalsplitting} (which immediately
generalized to the pro-$p$ completions) all the maps denoted with $i$ are
split injective. $A'$ is an isomorphism since localization is
an exact functor.

 The map $A$ is an isomorphism because of the first
part of the following proposition:
\begin{proposition}\label{prop:AHforfinitegroup}
  Suppose $Q$ is a finite $p$-group. Then the quotients
  \begin{equation*}
    \tilde\pi_S^n(Q)/F_k^n(Q)
  \end{equation*}
  are finite $p$-groups for each $n,k\in\integers$. Moreover, 
  \begin{equation*}
\tilde\pi^n_S(Q)=\invlim_{k\to\infty} \tilde\pi^n_S(Q)/F_k^n(Q),
\end{equation*}
in particular the intersection of all the $F_k^n(Q)$ is zero for each
$n\in\integers$.
\end{proposition}

Next, we claim that the map $B$ is injective. This follows from the following Theorem
\ref{theorem:cohomotopy_completion_iso} because of  Proposition
\ref{tfcompletion}. Theorem
\ref{theorem:cohomotopy_completion_iso} applies since by assumption
$B{H}$ is of finite type, so that the same is true for $BG$, using Lemma
\ref{lemma:fibrationtype} (together with \cite[Example
  3.4]{Adem-Milgram(1994)} for the finite quotient $Q$).

\begin{theorem}\label{theorem:cohomotopy_completion_iso}
  Suppose $G$ is cohomologically $p$-complete and has a classifying
  space of finite type. If $n\le 0$, then the induced map
  \begin{equation*}
   B^n\colon \tilde\pi^n_S({\hat G}^p)_{(p)}/F^n_k({\hat G}^p)_{(p)} \to
    \tilde\pi^n_S(G)_{(p)}/F^n_k(G)_{(p)}
  \end{equation*}
  is injective for each $k\ge 0$.
\end{theorem}

Given these properties of the diagram \eqref{eq:bigdiagram}, it is
easy to conclude that the induced map $\tilde
\pi^0_S(Q)\into\tilde\pi^0_S({G})$ is injective:

If
$x\in\tilde\pi^0_S(Q)$ is mapped to zero in $\tilde\pi_S^0({G})$, it
is also mapped to zero in each of the quotients
$\tilde\pi^0_S({G})_{(p)}/F^0_k({G})_{(p)}$, therefore (using injectivity
of the rightmost column of the diagram) in each of the quotients
$\tilde\pi^0_S(Q)/F^0_k(Q)$. Hence, $x$ is contained in
$\bigcap_{k}F^0_k(Q)=\{0\}$ (where we use the second part of
Proposition \ref{prop:AHforfinitegroup}), i.e.~$x=0$, as we had to
prove. Theorem \ref{theorem:splitback} follows now by applying Theorem
\ref{Adem}.

It remains to prove Proposition \ref{prop:AHforfinitegroup} and
Theorem \ref{theorem:cohomotopy_completion_iso}.

First remember: if $X$ is a finite dimensional CW-complex, we have the
Atiyah-Hirzebruch spectral sequence \cite[Theorem
4.2.7]{Kochman(1996)} converging to
$\tilde\pi^{s+t}_S(X)$, with $E_2$-term $E_2^{s,t}(X) = \tilde
H^s(X,\pi_{-t}^S)$, where $\pi_*^S:=\pi_*^S(pt)$ denotes the
coefficients of stable homotopy. 

Recall the following \cite[Proposition
  5.1.1]{Kochman(1996)}:
\begin{theorem}\label{stable_homotopy} The coefficients of stable
  homotopy satisfy:
  $\pi_t^S$ is finite for $t>0$, $\pi_0^S\iso\integers$, and
  $\pi_t^S=0$ for $t<0$.
\end{theorem}

This implies in particular that the Atiyah-Hirzebruch spectral
sequence for stable cohomotopy is a fourth quadrant spectral sequence.

We are going to apply this to the skeleta of our (infinite
dimensional) classifying spaces $BG$ and $BQ$. First we study
arbitrary finite or profinite groups.

\begin{proposition}\label{prop:AHfinite2}
  Let $Q$ be a finite group. Then we have a fourth quadrant
  Atiyah-Hirzebruch
  spectral sequence $E^{s,t}_r(Q)$, converging to
  $\tilde\pi^{s+t}_S(Q)$, with $E_2$-term
  \begin{equation*}
E_2^{s,t}(Q)=\tilde H^s(Q,\pi_{-t}^S).
\end{equation*}
Moreover, for
  fixed $s,t\in\integers$ and then $k$ sufficiently large $E_k^{s,t}$
  does not depend on $k$
  (this defines $E_\infty^{s,t}$),
  \begin{equation*}
E_\infty^{s,t}=
  F^{s+t}_{s}(Q)/F^{s+t}_{s+1}(Q)\quad\text{and also}\quad
  \tilde\pi_S^n(Q)=\invlim_{k\to\infty}
  \tilde\pi^n_S(Q)/F^n_k(Q).
\end{equation*}

   If $\profinQ$ is a profinite group then we have a fourth quadrant
   Atiyah-Hirzebruch spectral sequence $E^{s,t}_r(\profinQ)$ with $E_2$-term
  \begin{equation*}
E_2^{s,t}(\profinQ)=\tilde H^s(\profinQ,\pi_{-t}^S).
\end{equation*}
If, in addition,
$H^k(\profinQ,\integers/n)$ is finite for every $k\in\naturals_0$,
$n\in\naturals$ (with
trivial action of $\profinQ$ on $\integers/n$), then for   each fixed
$s,t\in\integers$ with $t\ne 0$ and for $k$ sufficiently large $E_k^{s,t}$ 
  does not depend on $k$, and
we have injections of finite abelian groups
  \begin{equation*}
  F^{s+t}_{s}(\profinQ)/F^{s+t}_{s+1}(\profinQ) \hookrightarrow
  E_\infty^{s,t}(\profinQ)
  \quad\text{for }s+t\le 0.
\end{equation*}
\end{proposition}
\begin{proof}
  If $Q$ is finite, then the existence of the spectral sequence with
  the given $E_2$-term is
  standard (compare the remark before \cite[Theorem
  4.2.7]{Kochman(1996)}). The issue is convergence.
  
  Since $Q$ is finite,  $\tilde
  H^s(Q,\pi_{-t}^S)$ is a torsion group \cite[Corollary
  10.2]{Brown(1982)}. Since $\pi_t^S$
  is finitely
  generated (Theorem \ref{stable_homotopy}) and $BQ$ is of finite
  type \cite[Example
  3.4]{Adem-Milgram(1994)}, $\tilde
  H^s(Q,\pi_{-t}^S)$ is a finitely generated torsion group and
  therefore finite for every $s,t\in\integers$.
  Since $E_k^{s,t}(Q)$ is a sub-quotient of $E_{k-1}^{s,t}(Q)$, 
  for each $s,t\in\integers$ these groups
   eventually become constant.

We also have the corresponding Atiyah-Hirzebruch
  spectral sequences for the finite skeleta $BQ^{(n)}$, and these
  converge in the usual sense \cite[Theorem
  4.2.7]{Kochman(1996)}. Now $\tilde H^s(BQ^{(n)},\pi_{-t}^S)=0$ for
  $s>n$ because $BQ^{(n)}$ is $n$-dimensional, and $\tilde H^s(Q,\pi_{-t}^S)=\tilde
  H^s(BQ^{(n)},\pi_{-t}^S)$ for $s<n$. Since $BQ^{(n)}$ is a finite
  CW-complex, Theorem \ref{stable_homotopy} implies immediately that
  $\tilde H^s(BQ^{(n)},\pi^S_{-t})$ is finite whenever $t\ne
  0$. Consequently, $E_2^{s,t}(BQ^{(n)})$ and by Lemma
  \ref{lem:finspecseq}
  $E_\infty^{s,t}(BQ^{(n)})$ are finite except possibly for $(s,t)=(
  n,0)$, and are zero for $s>n$. 

  By Lemma \ref{lem:finspecseq}
  $\tilde \pi_S^{s+t}(BQ^{(n)})$ is finite for $s+t\ne n$, and is zero for
  $s+t>n$.

  Observe that $BQ=\dirlim_{n\to\infty} BQ^{(n)}$. However,
  cohomotopy, as a cohomology theory, is in general not compatible
  with such limits. By \cite[Proposition 7.66 and Remark 1 on
  p.~132]{Switzer(1975)} our finiteness results imply that here the
  Mittag-Leffler condition
  is fulfilled and therefore 
  \begin{equation*}
    \tilde\pi_S^*(Q)=\invlim_{n\to\infty} \tilde\pi_S^*(BQ^{(n)}).
  \end{equation*}
  Observe that we can replace $\tilde\pi^*_S(BQ^{(n)})$ by the image of
  $\tilde\pi^*_S(Q)$ in $\tilde\pi^*_S(BQ^{(n)})$ without changing
  the inverse limit (an obvious general fact about inverse limits), and this
  image of course is isomorphic to $\tilde\pi^*_S(Q)/F^*_n(Q)$
  (since $F^*_n(Q)$ is the kernel of the corresponding projection
  map). Therefore
  \begin{equation*}
    \tilde\pi^*_S(Q) = \invlim_{n\to\infty} \tilde\pi^*_S(Q)/F^*_n(Q).
  \end{equation*}
  Next observe that by definition we have (for $n<N$) the exact sequences
  \begin{equation*}
    \begin{CD}
      0 @>>> F^*_n(Q) @>>> \tilde\pi^*_S(Q) @>>>
      \tilde\pi^*_S(BQ^{(n)})\\
       && @VVV  @VVV   @VV{=}V\\
      0 @>>> F^*_n(BQ^{(N)}) @>>> \tilde\pi^*_S(BQ^{(N)}) @>>>
      \tilde\pi^*_S(BQ^{(n)}).
    \end{CD}
\end{equation*}
    By commutativity, the images of $\tilde\pi^*_S(Q)$ and of
    $\tilde\pi^*_S(BQ^{(N)})$ in $\tilde\pi^*_S(BQ^{(n)})$
    coincide. Now take the inverse limit of the second line for
    $N\to\infty$. The right-most term is constant and the middle term
    converges to the middle term in the first line. Since inverse
    limit is left exact \cite[Proposition 7.63]{Switzer(1975)},
    \begin{equation*}
F^*_n(Q)=\invlim_{N\to\infty} F^*_n(BQ^{(N)}).
\end{equation*}
  
We now want to establish that for $k<l$
\begin{equation}\label{eq:helplimit}
F_k^*(Q)/F^*_l(Q)=\invlim_{n\to\infty}\left( F^*_k(BQ^{(n)})/F^*_l(BQ^{(n)})\right).
\end{equation}
  By the limit result we just obtained and by \cite[Proposition
  7.63]{Switzer(1975)}, we have the exact sequence 
  \begin{equation*}
    0\to F_k^*(Q)\to F^*_l(Q)\to \invlim_{n\to\infty}
    F^*_k(BQ^{(n)})/F^*_l(BQ^{(n)})
    \to {\lim_{n\to\infty}}^1(F_k^*(BQ^{(n)})).
  \end{equation*}
  For $s<n$, $F^s_k(BQ^{(n)})$ is a sub-quotient of the
  cohomotopy group $\pi^s_S(BQ^{(n)})$. We proved  above that the
  latter one is finite. Therefore $F^s_k(BQ^{(n)})$ is finite for $n>s$.
 Now \cite[Remark 1 on p.~132]{Switzer(1975)} implies that
the 
  $\lim^1$-term is zero, hence \eqref{eq:helplimit} follows.

  Last observe that for $n>s$
  \begin{equation}\label{eq:equality}
E_2^{s,t}(Q)=\tilde H^s(Q,\pi^S_{-t}) = \tilde
H^s(BQ^{(n)},\pi^S_{-t}) = 
  E_2^{s,t}(BQ^{(n)}),
\end{equation}
since the inclusion $BQ^{(n)}\to BQ$ is an $n$-connected map.

Fix now $s,t\in\integers$. Since $E^{s,t}_2(Q)$ is finite, there is
$r_0>s$ (depending on
$s,t$) such that $E^{s,t}_r(Q)= E^{s,t}_\infty(Q)$ for every $r\ge
r_0$. By Lemma \ref{lem:specseq} and Equation \eqref{eq:equality}
there exists $n_0\ge 0$ such that
\begin{equation*}
  E^{s,t}_{r_0}(Q)= E^{s,t}_{r_0}(BQ^{(n)})\qquad\text{for all}\quad n\ge n_0.
\end{equation*}
Note that $n_0$ depends on $r_0, s,t$, and a priori, $n_0$ might
become bigger if we choose a larger $r_0$. But, since $r_0>s$, we have
for each $n\ge n_0$
a commutative diagram
\begin{equation*}
  \begin{CD}
    E_\infty^{s,t}(Q) @= E^{s,t}_{r_0}(Q)\\
    @VV\alpha V  @VV=V\\
    E_\infty^{s,t}(BQ^{(n)}) @>\beta>\subset> E_{r_0}^{s,t}(BQ^{(n)}).
  \end{CD}
\end{equation*}
Note that $\beta$ is injective since $r_0>s$. Since $\beta\circ\alpha$
is an isomorphism, $\beta$ is an isomorphism, and the same then is
true for $\alpha$, i.e.~all groups $E^{s,t}_\infty(Q)$,
$E^{s,t}_r(Q)$, $E^{s,t}_{r_0}(Q)$, $E^{s,t}_r(BQ^{(n)})$,
$E^{s,t}_\infty(BQ^{(n)})$ are identified, as long as $r\ge r_0$ and
$n\ge n_0$.

In particular, we can write
  \begin{equation}\label{eq:33}
    E_r^{s,t}(Q)=\invlim_{n\to\infty} E_r^{s,t}(BQ^{(n)})\qquad 2\le r\le
    \infty;\; s,t\in\integers,
  \end{equation}
  since for fixed $s,t,r$ (including $r=\infty$) the sequence of
  groups on the right stabilizes. Since $BQ^{(n)}$ is a finite CW-complex,
  $F_s^{s+t}(BQ^{(n)})/F^{s+t}_{s+1}(BQ^{(n)})\iso
  E^{s,t}_\infty(BQ^{(n)})$. By \eqref{eq:helplimit} and \eqref{eq:33},
  both
  $E_\infty^{s,t}(Q)$ and
  $F^{s+t}_s(Q)/F^{s+t}_{s+1}(Q)$ coincide with the same group
  $\invlim_{n\to\infty}
  E_\infty^{s,t}(BQ^{(n)})$, proving that they are equal as asserted.

  If $\profinQ$ is profinite, we define the $E_r$-term of the Atiyah-Hirzebruch
  spectral sequence for $2\le r<\infty$
  to be the direct limit of the spectral sequences for its finite
  quotients $\{Q_i\}$. This makes sense since the direct limit functor is
  exact, and because of the naturality of the Atiyah-Hirzebruch
  spectral sequence. Since this spectral sequence is concentrated in the fourth
  quadrant, by Lemma \ref{lem:specseq} $E_r^{s,t}(\profinQ)\subset
  E_{s}^{s,t}(\profinQ)$ for $r> s$. In particular,
  $E^{s,t}_\infty(\profinQ)=\lim_{r\to\infty} E_r^{s,t}(\profinQ)$ makes sense, and,
  moreover, $E^{s,t}_\infty(\profinQ)=\bigcap_{r>s} E_r^{s,t}(\profinQ)$.

  Let $Q_i$ be the system of finite quotients of $\profinQ$. Note that
  inverse limits and direct limits do not necessarily commute. This
  means here that the direct limit of
  $E_\infty^{s,t}(Q_i)=\bigcap_{r>s} E_r^{s,t}(Q_i)$ is not necessarily
  $E^{s,t}_\infty(\profinQ)$. Exactness of direct limits implies
  however, using the results about finite groups proved above, that
  \begin{equation*}
\dirlim_{i} E_\infty^{s,t}(Q_i)
=\left(\dirlim_{i}F_s^{s+t}(Q_i)\right)/\left(\dirlim_{i}
  F_{s+1}^{s+t}(Q_i)\right) = F_s^{s+t}(\profinQ)/ F_{s+1}^{s+t}(\profinQ).
\end{equation*}
Assume now that $H^k(\profinQ,\integers/n)$ is a finite group for each
$k\in\naturals_0$, $n\in\naturals$. Since these finite groups are the
direct limits of the 
cohomology groups
of the $Q_i$, this implies that, for fixed $n$, $H^k(Q_i,\integers/n)\to
H^k(\profinQ,\integers/n)$ is an isomorphism for all sufficiently large $i$.

Fix now $s,t\in\integers$ with $t\ne 0$. Observe that $E^{s,t}_r(\profinQ)$
is a sub-quotient of $H^s(\profinQ,\pi_{-t}^S)$ and therefore, since $t\ne 0$,
is finite by Theorem \ref{stable_homotopy} and our finiteness
assumption. Consequently, there is some $R=R(s,t)> s$ such that
$E^{s,t}_r(\profinQ)=E^{s,t}_{R}(\profinQ)$ for $r\ge R$. We require now
additionally that $s+t\le
0$. Following Lemma
\ref{lem:specseq}, choose a finite set $V\subset\integers^2$ such that
$\{E^{x,y}_2(\profinQ) \mid (x,y)\in V\}$ determines
$E^{s,t}_R(\profinQ)$. Since
$s+t\le 0$ and all our
spectral sequences are concentrated in the fourth quadrant, we can
assume that $(x,0)\notin V$ if $x>0$. By the above remark, if $i$
is sufficiently large (depending on $s,t,R$)
\begin{equation*}
E^{x,y}_2(Q_i)=\tilde H^x(Q_i,\pi_{-y}^S)\to \tilde
H^x(\profinQ,\pi_{-y}^S)=E^{x,y}_2(\profinQ)
\end{equation*}
is an isomorphism for all $(x,y)\in V$. Since $\pi_0^S\iso\integers$
and we did only make an assumption about cohomology with finite
coefficients, we carefully have to avoid $y=0$, except that $\tilde
H^0=0$, so $x=0=y$ is permitted.
The choice of $V$ implies that for these sufficiently large $i$ also
$E^{s,t}_R(Q_i)\to
E^{s,t}_R(\profinQ)$ is an isomorphism. We picked $R>s$, hence
$E^{s,t}_\infty(Q_i)\subset E^{s,t}_R(Q_i)=E^{s,t}_\infty(\profinQ)$. Passage
to the limit gives
\begin{equation*}
  \dirlim_i E^{s,t}_\infty(Q_i) = F_s^{s+t}(\profinQ)/ F_{s+1}^{s+t}(\profinQ)
  \subset E_\infty^{s,t}(\profinQ).\hfill\qed
\end{equation*}
\renewcommand{\qed}{}
\end{proof}

Observe that we now essentially have proved Proposition
\ref{prop:AHforfinitegroup}. It only remains to remark that if $Q$ is
a finite $p$-group then $\tilde
H^*(Q,A)$ is a finite (abelian) $p$-group for arbitrary finitely generated
coefficients (compare \cite[Corollary 5.4]{Adem-Milgram(1994)}). By
Lemma \ref{lem:finspecseq},
$\tilde\pi^s_S(Q)/F^s_k(Q)$ is a finite (abelian) $p$-group for each
$k\in\naturals_0$ and $s\in\integers$.

We now study an infinite group $G$.
\begin{proposition}\label{prop:AHforH}
  Assume that $G$ is a discrete group with classifying space $BG$ of finite
  type. Then we have a fourth quadrant
  Atiyah-Hirzebruch
  spectral sequence $E^{s,t}_r(G)$ for stable cohomotopy with
  $E_2^{s,t}(G)=\tilde H^s(G,\pi_{-t}^S)$.

  The following partial convergence results are true: for
  each fixed $s,t\in\integers$ with $t\ne 0$ or for $s=0=t$
  eventually $E_k^{s,t}$ becomes constant
  (this defines $E_\infty^{s,t}$); and if $s+t\le 0$ then 
  \begin{equation*}
E_\infty^{s,t}(G)=
  F^{s+t}_{s}(G)/F^{s+t}_{s+1}(G)  \quad\text{and}\quad
  \tilde\pi_S^{s+t}(G)=\invlim_{k\to\infty}
  \tilde\pi^{s+t}_S(G)/F^{s+t}_k(G).
\end{equation*}
\end{proposition}
\begin{proof}
  We can apply all arguments of the proof of Proposition
  \ref{prop:AHfinite2}. When using finiteness conditions of
  cohomology, we have 
  to be careful to avoid (for $s>0$) $\tilde
  H^s(G,\pi_0^S)=\tilde H^s(G,\integers)$, 
  and $\tilde H^s(G^{(n)},\pi_0^S)=\tilde H^s(G^{(n)},\integers)$
  which were finite (if $s<n$) with $G$
  replaced by the finite group $Q$, but are not necessarily finite
  here. But one checks immediately that the arguments go through for
  the range of
  $s$ and $t$ stated in the proposition. 
\end{proof}

\begin{corollary}\label{corol:localized_AHSS}
Fix a prime number $p$.
  Assume that $G$ is a discrete group with classifying space $BG$ of finite
  type. Then we have a fourth quadrant localized
  Atiyah-Hirzebruch
  spectral sequence $E^{s,t}_r(G)_{p}$ with
  $E_2^{s,t}(G)_{p}=\tilde H^s(G,\pi_{-t}^S)_{(p)}$. For $t\ne 0$ this
  computes as
  \begin{equation}\label{eq:E21}
    E_2^{s,t}(G)_{p} = \tilde H^{s}(G, (\pi^S_{-t})_{(p)}).
  \end{equation}

  The following partial convergence results are true: for
  each fixed $s,t\in\integers$ with $t\ne 0$ or for $s=0=t$,
  if $k$ is sufficiently large then $E_k^{s,t}(G)_p$ becomes
  independent of $k$
  (this defines $E_\infty^{s,t}(G)_p$). If $s+t\le 0$, then 
  \begin{equation*}
E_\infty^{s,t}(G)_p=
  F^{s+t}_{s}(G)_{(p)}/F^{s+t}_{s+1}(G)_{(p)}. 
\end{equation*}

  If $G$ is finite, the last two statements are true for arbitrary
  $s,t\in\integers$.

  If $\profinQ$ is a profinite group, then we have a
  fourth quadrant localized 
  Atiyah-Hirzebruch
  spectral sequence $E^{s,t}_r(\profinQ)_{p}$ with
  $E_2^{s,t}(\profinQ)_{p}=\tilde
  H^s(\profinQ,\pi_{-t}^S)_{(p)}$. 

Assume in addition that
  $H^k(\profinQ,\integers/p)$ is finite for every $p$ and $k$. 
  For $t\ne 0$, we then have
  \begin{equation}\label{eq:E22}
    E_2^{s,t}(\profinQ)_p = \tilde
    H^s(\profinQ,(\pi_{-t}^S)_{(p)}),
  \end{equation}
 and the following
  convergence results are true:
for
  each fixed $s,t\in\integers$ with $t\ne 0$ or for $s=0=t$
  and for $k$ sufficiently large $E_k^{s,t}(\profinQ)_p$ becomes independent
  of $k$
  (this defines $E_\infty^{s,t}(\profinQ)_p$); and if $s+t\le 0$ then we get
  an injection
  \begin{equation*}
  F^{s+t}_{s}(\profinQ)_{(p)}/F^{s+t}_{s+1}(\profinQ)_{(p)} \into
  E_\infty^{s,t}(\profinQ)_p . 
\end{equation*}
\end{corollary}
\begin{proof}
  All the localized spectral sequences are constructed by localizing
  the spectral sequences obtained in Proposition
  \ref{prop:AHforH} or Proposition
  \ref{prop:AHfinite2} at $p$. Observe that this construction produces
  spectral sequences since localization is an exact functor.

  For  the computation of $E_2$-terms \eqref{eq:E21} and
  \eqref{eq:E22}, note that $\pi_{-t}^S$ is the
  (finite) direct sum of the
  $q$-primary components for all primes $q$, if $t<0$. Accordingly, $\tilde
  H^s(G,\pi_{-t}^S)$
  splits as a direct sum, and the $q$-primary part of $\pi_{-t}^S$
  gives rise to the $q$-primary part of the cohomology group. Since
  for abelian torsion groups localization at $p$ is just projection to
  the $p$-primary part,
  $\tilde H^s(G,\pi_{-t}^S)_{(p)}=\tilde
  H^s(G,(\pi_{-t}^S)_{(p)})$. The same argument
  applies in the profinite case.

  If $\profinQ$ is a profinite group with $H^k(\profinQ,\integers/p)$ finite for
  every $k\in\naturals_0$, then because of Lemma \ref{lemma:coeffiso}
  and the statement we have just proved, $\tilde
  H^s(\profinQ,\pi_{-t}^S)_{(p)}=\tilde H^s(\profinQ,(\pi_{-t}^S)_{(p)})$ is
  finite 
  for $t\ne 0$. 
 
  The remaining assertions follow now easily from Proposition
  \ref{prop:AHforH} or Proposition \ref{prop:AHfinite2}, using again the
  fact that localization is an exact functor. 

A priori, one has to be careful with any statement about the
$E_\infty$-terms, because $E^{s,t}_\infty(G)$ is defined as the inverse
limit of the $E^{s,t}_r(G)$, and localization does not in general
commute with inverse limits. However, here every sequence defining the
inverse limits stabilizes (before localization),
i.e.~$E^{s,t}_\infty(G)=E^{s,t}_r(G)$ for $r$ sufficiently large. The
same then is true after localization. This implies
$E^{s,t}_\infty(G)_p=\left(E^{s,t}_\infty(G)\right)_{(p)}$, and,
consequently, the statements for $(E_\infty)_p$ follow also from
Propositions \ref{prop:AHforH} and \ref{prop:AHfinite2} and exactness
of the localization functor.

The arguments given so far apply identically to a profinite group $\profinQ$.
\end{proof}

\begin{remark}
  Note that localization does not in general commute with inverse
  limits. Hence
  we can't conclude that $\pi^s_S(G)_{(p)}=\invlim_{k\to\infty}
  \pi^s_S(G)_{(p)}/F^s_k(G)_{(p)}$ for any
  $s\in\integers$. However, we will see that the spectral sequence
  determines
  enough of the groups
  $\pi^s_S(G)_{(p)}/F^s_k(G)_{(p)}$.
\end{remark}

Observe that all the spectral sequences we have constructed are
natural with respect to group homomorphisms. In particular, the map to
the pro-$p$ completion $c\colon G\to{\hat G}^p$ induces a map of localized
spectral sequences 
\begin{equation*}
E^{s,t}_r(c)\colon E^{s,t}_r({\hat G}^p)_p\to E^{s,t}_r(G)_{p}.
\end{equation*}

If $BG$ is of finite type and $G$ is cohomologically $p$-complete,
Corollary \ref{corol:localized_AHSS} implies that the induced map on
the $E_2$-pages is for $t<0$
\begin{equation*}
\tilde H^s({\hat G}^p,(\pi_{-t}^S)_{(p)})\to
\tilde H^s(G,
(\pi_{-t}^S)_{(p)}).
\end{equation*}
Because of cohomological $p$-completeness
and Lemma \ref{lemma:coeffiso} this map
is an isomorphism for $t<0$ or (trivially) for $s=0=t$. If $t=0$ the
coefficients are not finite so that this can not be extended to other
values of $s$ and $t$. Using Lemma \ref{lem:specseq}, $c$ therefore
induces isomorphisms on
$E_r^{s,t}$ in particular for $s+t\le 0$. By the convergence result of
Corollary
\ref{corol:localized_AHSS}, the map hence induces injections
\begin{multline*}
  F^s_{k-1}({\hat G}^p)_{(p)}/F^s_k({\hat G}^p)_{(p)} \into
  E^{k-1,s-k+1}_\infty({\hat G}^p)_p\\
   \iso E^{k-1,s-k+1}_\infty(G)_p =
  F^s_{k-1}(G)_{(p)}/F^s_k(G)_{(p)}.
\end{multline*}
As in the proof of Lemma \ref{lem:finspecseq}, we get injections
\begin{equation*}
 \pi_S^s(c)_k\colon  \pi_S^s({\hat G}^p)_{(p)}/F^s_k({\hat G}^p)_{(p)} \into
     \pi_S^s(G)_{(p)}/F^s_k(G)_{(p)}\qquad\text{for all }k\ge 0.
   \end{equation*}
This proves Theorem \ref{theorem:cohomotopy_completion_iso} and hence
concludes the proof of Theorem \ref{theorem:splitback}. \qed

\subsection{Splitting of short exact sequences}

We now give a condition on an extension $1\to {H}\to G\to Q\to 1$ as in Theorem
\ref{theorem:splitback} which implies the
splitting of the induced projection ${\hat G}^p\to Q$.

\begin{lemma}\label{lem:completionsplit_criterion}
  Let $Q$ be a finite $p$-group in the exact sequence
  \begin{equation*}
    1\to {H}\to G\xrightarrow{\pi} Q\to 1.
  \end{equation*}
  Assume that, among all
  normal subgroups of $G$ of $p$-power index, there is a cofinal
  system  $U_i\normalsubgroup G$  with
  $U_i\subset {H}$, such 
  that for each $i$ the homomorphism $\pi_i$ in
  \begin{equation*}
    G\xrightarrow{p_i} G/U_i\stackrel{\pi_i}{\onto}Q
  \end{equation*}
  has a split $s_i\colon Q\into G/U_i$. Then the pro-$p$ completion  $ {\hat
    \pi}^p\colon {\hat
  G}^p \onto Q$ has a split $s\colon Q\into {\hat G}^p$, too.
\end{lemma}
\begin{proof}
  For each of the finitely many $q\in Q$ choose $x(q,i)\in G$ with
  \begin{equation*}
p_i(x(q,i))=s_i(q)\in G/U_i.
\end{equation*}
Choose a cofinal subsequence such that
  for every $q\in Q$ the images of $x(q,i)$ in the compact group
  ${\hat G}^p$ converge to an element $x(q)\in{\hat G}^p$ (limit over
  $i$). We claim that the map 
  \begin{equation*}
s:Q\to{\hat G}^p\colon q\mapsto x(q)
\end{equation*}
is a
  split of ${\hat \pi}^p$. Since $\pi_i \circ p_i(x(q,i)) =q\in Q$, by
  continuity 
  also
  ${\hat\pi}^p(x(q))=q$. It remains to check that $s$ is a group
  homomorphism. Observe that for $q_1,q_2\in Q$ 
  \begin{equation*}
    s(q_1q_2^{-1})s(q_2)s(q_1)^{-1} = \lim_i
    x(q_1q_2^{-1},i)x(q_2,i)x(q_1,i)^{-1} .
  \end{equation*}
  The topology of ${\hat G}^p$ and the fact that 
  \begin{equation*}
    p_i(x(q_1q_2^{-1},i)x(q_2,i)x(q_1,i)^{-1}) =
    s_i(q_1q_2^{-1})s_i(q_2)s_i(q_1)^{-1} = 1 \in G/U_i
  \end{equation*}
  implies that this limit is $1$, i.e.~$s$ indeed is a homomorphism.
\end{proof}

\subsection{Enough torsion-free quotients}

In this subsection we establish conditions which guarantee that we can
factorize a projection onto a finite $p$-group through a
torsion-free elementary amenable quotient. This is another key
ingredient in our strategy to prove the Atiyah conjecture for extensions.

\begin{definition}\label{def:tfquotients}
  Let $G$ be a discrete group. We define $\tfsubgroups_G$ to be the
  set of all normal subgroups $U\subgroup G$ such that $G/U$
  is torsion-free and elementary amenable.  $G$ is said to have
  \emph{enough torsion-free
  quotients} if each normal subgroup  $W\normalsubgroup G$, with $G/W$
  a finite
  $p$-group for some prime number $p$, contains a subgroup
  $U_W\in\tfsubgroups_G$, i.e.~we have a factorization 
  \begin{equation*}
    G\onto G/U_W\onto G/W
  \end{equation*}
  with $G/U_W$ torsion-free elementary amenable.

A subset 
   $\tfsubgroups_G'$  of $\tfsubgroups_G$ is called \emph{exhaustive},
   if 
  $U,V\in\tfsubgroups'_G$ implies $U\cap V\in\tfsubgroups'_G$, if for
   every $n\in\naturals$ there is $U_n\in\tfsubgroups'_G$ such that
   $U_n\subgroup \gamma_n(G)$,
  and if, in addition, each $W$ as above contains all but finitely many of
  the $U\in\tfsubgroups_G'$. In particular, all the $U_W$ can be chosen
  to belong to $\tfsubgroups'_G$.

  We say that $G$ has \emph{enough nilpotent torsion-free quotients}
  if we can choose an exhaustive set $\tfsubgroups'_G$ such that, if
  $U\in\tfsubgroups'_G$, then $U$ is a characteristic subgroup of $G$
  and  $G/U$
  is nilpotent. Analogously, we define
  \emph{enough solvable torsion-free quotients} and \emph{enough solvable-by-finite torsion-free quotients}. The corresponding
  $\tfsubgroups'_G$ then is called \emph{nilpotent exhaustive},
  \emph{solvable exhaustive}, or
  \emph{solvable-by-finite exhaustive}, respectively.
\end{definition}

\begin{lemma}\label{lem:tfprime_exists_primitive}
  Suppose $G$ is a finitely generated group. Let $\mathcal{X}$ be a set
  of normal subgroups of $G$ with torsion-free elementary amenable
  quotients, and such that whenever $N\normalsubgroup G$
  and $G/N$ is a finite $p$-group we can find $U\in\mathcal{X}$ with
  $U\subgroup N$. Assume that $U,V\in\mathcal{X}$ implies $U\cap
  V\in\mathcal{X}$. 

  Then there is a nested sequence $U_1\supergroup U_2\supergroup
  U_3\supergroup\cdots$ with $U_k\subgroup \gamma_k(G)$ for every
  $k\in\naturals$ and such that for every finite $p$-group quotient
  $G/N$ there is $k_N$ with $U_{k_N}\subgroup N$. Every $U_k$ is an
  intersection $U_k=V_k\cap N_k$ with $V_k\in\mathcal{X}$ and $G/N_k$
  torsion-free nilpotent.

  In particular, if $G$ has enough torsion-free quotients, an
  exhaustive subset as in Definition \ref{def:tfquotients} exists.
\end{lemma}
\begin{proof}
  Define for a prime number $p$ the subgroup $G_{p,n}$
  to be the intersection of the
  kernels of all  projections to $p$-groups which have order
dividing
  $p^n$. Since $G$ is finitely generated, there are only finitely many
  of these. By Lemma \ref{lem:subgroupintersect},
$G/G_{p,n}$ is a finite
  $p$-group.

  By assumption, for each $G_{p,n}$ we find $U_{p,n}\in\mathcal{X}$
  which is contained in $G_{p,n}$. For $n\in\naturals$ define $V'_n$
  to be the intersection of all $U_{p,k}$ with $p\le n$ and $k\le n$.

  For $n\in\naturals$ we now want to replace
  the subgroup $V'_n$ we just constructed by $U_n$ such that
  $U_n\subgroup\gamma_n(G)$. Since
  $G/\gamma_n(G)$ is nilpotent, the set
  $K_n$ of elements of finite order in $G/\gamma_n(G)$ is a
  characteristic subgroup and the quotient $(G/\gamma_n(G))/K_n$ is
  torsion-free \cite[Theorem 2.22]{Hirsch(1938)}. Let $H_n$ be the
  kernel of the projection $G\to (G/\gamma_n(G))/K_n$. By
  \cite[Theorem 2.1]{Gruenberg(1957)} for each element $1\ne g\in K_n$ we
  find a finite $p$-group $Q_g$ ($p$ depending on $g$) and a
  projection $G/\gamma_n(G)\onto Q_g$ such that $g$ is not mapped to
  $1$. For each such $g$ choose $n_g\in\naturals$ such that
  we get a factorization 
  $G\onto V'_{n_g}\onto Q_g$. Set 
  \begin{equation*}
   V_n:= V'_n\cap \left(\bigcap_{k\le n}\bigcap_{g\in
  K_k} V'_{n_g}\right)\text{ and } U_n:=V_n\cap H_n.
\end{equation*}
This defines a new nested sequence of subgroups. We have seen that
$G/H_n$ is torsion-free elementary amenable. By Lemma
 \ref{lem:subgroupintersect}, the same is true for $G/U_n$. 
 Fix $u\in U_n$. Since $U_n\subset H_n$, the image of $u$ in
  $G/\gamma_n(G)$ is contained in $K_n$. For $1\ne g\in K_n$, $U_n$ is
  mapped to $1$ in each of the quotients $Q_g$.
  But $g$ is not mapped to $1$, hence $u\in\gamma_n(G)$,
 i.e.~$U_n\subset \gamma_n(G)$.

  It is now clear that
$\tfsubgroups_G':=\{ U_n\mid n\in\naturals\}$ is an exhaustive set.
\end{proof}

\begin{lemma}\label{lem:tfprime_exists} Assume that
  $G$ is a finitely generated group.

  If every projection $G\onto Q$ onto a finite $p$-group $Q$ factors
  through a torsion-free nilpotent quotient $G\onto U\onto Q$, then
  $G$ has enough nilpotent torsion-free quotients. In particular, this
  depends only on the (directed) system of nilpotent quotients of $G$,
  ordered by inclusion of the kernels. 

  The corresponding statement holds with ``nilpotent'' replaced by
  ``solvable'' or
  ``solvable-by-finite''.
\end{lemma}
\begin{proof}
  Remember first that  every finite $p$-group is
  nilpotent. This means that the first condition
   in fact only depends on the system of nilpotent quotients.

 Let $\mathcal{X}_{tf}$ be the system of
  normal subgroups of $G$ with nilpotent torsion-free quotients. We
  can apply Lemma \ref{lem:tfprime_exists_primitive} to obtain a nested
  sequence $U_1\supergroup U_2\supergroup\cdots$ with
  $U_k\subgroup\gamma_k(G)$ for every $k\in\naturals$, and such
  that every projection to a finite
  $p$-group factors through $G/U_k$ for suitable $k\in\naturals$.
  By Lemma \ref{lem:subgroupintersect},
  $U_k\in\mathcal{X}_{tf}$ for every $k\in\naturals$.

 Now
  replace $U_k$ by
  $V_k:=\bigcap_{\alpha\in\Aut(G)}\alpha(U_k)$. Because of Lemma
\ref{lem:subgroupintersect},
$V_k\in\mathcal{X}_{tf}$ for every
  $k\in\naturals$, and these subgroups are characteristic.

  The proof for enough solvable or solvable-by-finite torsion-free
  quotients is exactly the
  same, replacing ``nilpotent'' everywhere with ``solvable'' or
  ``solvable-by-finite'', respectively.
\end{proof}

\begin{example}\label{example:series}
  If $G$ is a discrete group such that there are infinitely many lower
  central series quotients $G/\gamma_n(G)$ which are torsion-free, or
  infinitely many derived series quotients $G/G^{(n)}$ which are
  torsion-free, then we can choose the exhaustive subset
  $\tfsubgroups'_G$ to consist of the
  corresponding $\gamma_n(G)$ or $G^{(n)}$, respectively. In
  particular, $G$ has enough nilpotent torsion-free quotients, or
  enough torsion-free solvable quotients, respectively.
\end{example}

\begin{lemma}\label{lem:tf_solvable_series}
  Assume that $G$ is countable, has enough solvable torsion-free quotients and is
  residually
  torsion-free solvable. Then there is a nested sequence $G\supergroup
  U_1\supergroup U_2\supergroup \cdots$ of normal subgroups of $G$ such
  that 
  \begin{enumerate}
  \item $\bigcap_{n=1}^\infty U_n =\{1\}$,
  \item $U_n\subgroup\gamma_n(G)$ for every $n\in\naturals$,
  \item $G/U_n$ is torsion-free solvable for every $n\in\naturals$.
  \end{enumerate}
\end{lemma}
\begin{proof}
  Since $G$ is countable and residually torsion-free solvable, there
  is a nested sequence $G\supset H_1\supset H_2 \supset \cdots$ of normal
  subgroups such that $G/H_n$ is torsion-free solvable for every
  $n\in\naturals$ and such that $\bigcap_{n\in\integers} H_n=\{1\}$.

  On the other hand, since $G$ has enough solvable torsion-free
  quotients, by definition we find $V_n\subset \gamma_n(G)$ such that
  $G/V_n$ is torsion-free solvable for each $n\in \naturals$. Then
  $U_n:=H_n\cap V_n$ will satisfy the assertions.
\end{proof}

\begin{definition}\label{def:Uto_the_H}
Let $G$ be a group and let $U$ be a subgroup of $G$.
Then we set $U^G := \bigcap_{g \in G} U^g$, the largest normal
subgroup of $G$ which is contained in $U$.  Note that
$U^G\cap V^G=(U\cap V)^G$ for two subgroups
  $U,V$ of $G$.

Let $H$ be a subgroup of $G$ containing $U$
as a normal subgroup.  By Lemma \ref{lem:subgroupintersect} if
$H$ is normal in $G$ and $H/U$ is torsion-free, then
so is $H/U^G$.  On the other hand,
if $H$ has finite index in $G$, then not more
than $[G:H]$ of the $U^g$ are different. Consequently, if $H/U$ is
elementary amenable, then so is $G/U^G$.
\end{definition}

\subsection{The main result}

Now essentially everything is in place to prove the main (technical) theorem of
this paper. We need one more lemma:

\begin{lemma}\label{lem:tfquotient_split_criterion}
  Let $Q$ be a finite $p$-group and
let $1\to H\to G\to Q\to 1$ be an
  exact sequence of groups.
  Assume that $H$ is finitely generated and has enough torsion-free
  quotients. Specify an exhaustive
  family $\tfsubgroups'_H$ of subgroups as in Definition
  \ref{def:tfquotients}.  Fix $n\in\naturals$ and for every
  $U\in\tfsubgroups'_H$ assume that
  $G/U^G$ contains a finite subgroup $E_U/U^G$ of order $n$.

  Then there is a subgroup $Q_0\subgroup Q$ of order $n$ splitting
  back to ${\hat
  G_0}^p\subgroup {\hat G}^p$, where $G_0$ is defined by the exact
  sequence $1\to H\to G_0\to Q_0\to 1$.
\end{lemma}
\begin{proof}
  Since $H/U^G$ is torsion-free, the image of $E_U/U^G$ in $Q$ has order
  $n$, too. In particular, it follows that $n$ necessarily is a power
  of the prime number $p$. Let $a_U$ be the (finite) collection of
  subgroups of $Q$ of order 
  $n$ which are images of finite subgroups of $G/U^G$, $U\in \tfsubgroups'_H$.

  If $U,V\in \tfsubgroups'_H$, then
  each finite subgroup of $G/(U^G\cap V^G)$ maps isomorphically to a
  finite subgroup of $G/U^G$, since the kernel $U^G/(U^G\cap V^G)$ of the
  projection $G/(U^G\cap V^G)\onto G/U^G$ is contained in the torsion-free
  group $H/(U\cap V)^G$, and similarly for $G/V^G$. In particular,
  $a_U\cap a_V$ contains $a_{U\cap V}$ and therefore is
  non-empty. Since the set of all subgroups of $Q$ is finite, there
  exists $Q_0\in \bigcap_{U\in\tfsubgroups_H'} a_U$.

  In other words, $Q_0$ splits back to every $G/U^G$ for
  $U\in\tfsubgroups'_H$, and hence also to every quotient of $G/U^G$ which
  projects onto $G/H$. Among these quotients of $G/U^G$
  ($U\in\tfsubgroups'_H$) we can find a cofinal set
  of finite $p$-group
  quotients of $G_0$. This follows as in the proof of Lemma
  \ref{lemma:shortexactcompletion}, since $H$ has enough torsion-free
  quotients. By Lemma \ref{lem:completionsplit_criterion}, $Q_0$ splits
  back to ${\hat G_0}^p$. Last observe that, by Lemma
  \ref{lemma:shortexactcompletion},
  \begin{equation*}
1 \to{\hat H}^p\to{\hat G_0}^p\to
  Q_0\to 1\quad\text{and}\quad 1\to {\hat H}^p\to {\hat G}^p \to Q\to 1
\end{equation*}
are  exact, whence ${\hat G_0}^p$ is a subgroup of ${\hat G}^p$.
\end{proof}



\begin{theorem}\label{theorem:finiteext}
  Let $H$ be a discrete group with finite classifying space which has
  enough torsion-free quotients and which is
  cohomologically complete. Assume that $K H$ fulfills the
strong Atiyah conjecture
 (for some subfield $K=\overline K$
  of $\complexs$).
  \begin{enumerate}
  \item\label{item:main1}
    Let
  \begin{equation*}
1\to H\to G\to Q\to 1
\end{equation*}
be an exact sequence of groups, where $Q$ is a finite group. Then $K
G$ fulfills the strong Atiyah conjecture.
\item\label{item:main2}
Moreover, if $G$ is torsion-free and $\tfsubgroups_H'$ is an
exhaustive family of subgroups of $H$ as in Definition
\ref{def:tfquotients}, then we can find $U\in\tfsubgroups_H'$ such
that $G/U^G$ is torsion-free and elementary amenable.
\end{enumerate}
\end{theorem}


\begin{definition}
  Let $H$ be a group, $K$ a skew field and $DH$ a division ring which
  contains $KH$. We call $DH$ \emph{Hughes free} if for every subgroup
  $U$ of $H$ and set of representatives $X\subset H$ of $H/U$ the set
  $X$ (considered as subset of $DH$) is linearly independent over the
  division closure $DU$ of $KU$ in $DH$.
\end{definition}

\begin{corollary}\label{cor:amext}
  Assume that $H$ fulfills the assumptions of Theorem \ref{theorem:finiteext}
  and 
  \begin{equation*}
1\to H\to G\to A\to 1
\end{equation*}
 is an extension with  $A$ elementary
  amenable. If $\lcm(G)<\infty$ then $K G$ fulfills the strong Atiyah
  conjecture.

    Assume that $k$ is a domain and $k*H$ embeds into a skew field $D_H$
  such that the twisted action of $A$ on $k*H$ extends to
  $D_H$. If $G$ is torsion-free and $H$ has enough solvable-by-finite
  torsion-free quotients, then $k*G$ embeds into a skew field.

  If we only require $H$ to have enough torsion-free quotients (not
  necessarily solvable-by-finite), but
  assume that $D_H$ is Hughes free, then $k*G$
  embeds also into a skew field.
\end{corollary}
\begin{proof}[Proof of Corollary \ref{cor:amext}.]
  It is a direct consequence of
  Proposition \ref{prop:amext} and Theorem \ref{theorem:finiteext}
  that $KG$ fulfills the strong Atiyah
  conjecture. 

For the statement about $k*G$, given a finite
  subgroup $E/H$ of $G/H$ we attempt to construct a suitable skew
  field extension
  $D_E$ of $D_H$ such that $D_H*[E/H]$ embeds into $D_E$. The result
  then follows from Lemma \ref{Tvirabelian}.

  Observe that $H$ is finitely generated (since it has finite
  classifying space) and hence by Lemma \ref{lem:tfprime_exists} an
  exhaustive  family
  $\tfsubgroups_H'$ exists.
  Using Theorem \ref{theorem:finiteext} with $E$ instead of $G$, choose
  $U\in\tfsubgroups_H'$ such that $E/U^E$ 
  is torsion-free. In view of Lemma \ref{lem:subgroupintersect}, $E/U$
  is elementary amenable. 

Let $D_U$ be the division closure of $k*U^E$ in
  $D_H$. This is a skew field. The action of $E$ by conjugation on
  $D_H*[E/H]$ maps $k*U^E$ to 
  itself and therefore maps the division closure $D_U$ to
  itself. Consequently, we can form $D_U*[E/U]$ and $D_U*[H/U]$. Since
  $H/U$ and $E/U$
  are both
  torsion-free elementary amenable groups, and $D_U$ is a skew field,
  by Lemma \ref{Tvirabelian} $D_U*[E/U]$ is an Ore domain with an Ore
  localization (of course a skew field)  which we denote $D_{U,E}$. It has
  a sub-skew field
  $D_{U,H}$ which is the division closure of $D_U*[H/U]$, and of
  course at the same time is the Ore localization of this ring.

  If $D_H$ is Hughes free, then any set of
  representatives for $H/U$ in $D_H$ is by definition linearly
  independent over $D_U$. This implies that $D_{U,H}$ is equal to
  $D_H$, therefore 
  $D_H*[E/H]$ is a domain, since $D_{U,H}*[E/H]$ embeds into
  $D_E$. 

  However, we don't see a reason why $D_U$ should be equal to
  $D_{U,H}$  in
  general.
 Assume now that $H$ has enough solvable-by-finite torsion-free
  quotients and choose $\tfsubgroups_H'$ to be solvable-by-finite
  exhaustive. In view of Lemma \ref{lem:subgroupintersect}, and by
  replacing $U\in\tfsubgroups_H'$ by $U^G$ we may assume that each
  $U\in\tfsubgroups_H'$ is normal in $G$. Then
  we replace $D_H$ by a
  skew field built out of the
  (possibly many
  different) $D_{U,H}$. 

Using a non-trivial ultra-filter $\mathbb {F}$ for the set
  $\tfsubgroups_H'$ (i.e.~an ultrafilter which contains all subsets with finite
complement),
we form a new skew field $D_H''$.  Following
  \cite[\S 2.6]{Jacobson(1989)}
it consists of classes of all sequences $(d_U)_{U\in\tfsubgroups_H'}$ with
$d_U \in D_{U,H}$ (constructed as above), and
\[
(d_U)_{U\in\tfsubgroups_H'} = (d_U')_{U\in\tfsubgroups_H'}
\]
if and only if there is some $F \in \mathbb
{F}$ such that $d_U = d_U'$ for all $U \in F$. We embed $k*H$ in $D_H''$
via the map $x \mapsto (x,x,
\dots )$ and let $D_H'$ be the division closure of $k*H$ in
$D_H''$. In a similar manner we construct $D_E'$ using $D_{U,E}$
instead of $D_{U,H}$. This skew field contains $D_H'$. We should
remark that a finite number of the $D_{U,E}$ not
being a division ring (or not even being defined) does not affect 
$D_E'$ being a division ring;
in fact any finite number of the $D_{U,E}$ does not affect the
isomorphism class of $D_E'$. Our assumptions on $\tfsubgroups_H'$
assure that $E/U$ is torsion-free for all but finitely many of the $U\in\tfsubgroups_H'$.

Since $U$ is normal in $G$, the twisted action of $G/H$
on $k*H$  induces a twisted action on each $D_{U,H}$, hence on
$D_H''$ and consequently on $D_H'$.  It follows that we may form the
crossed product $D_H' *[G/H]$.

Next we show that for each finite subgroup $E/H$ of $G/H$, the
identity map on $k*E$ extends to a monomorphism
from
$D_H'*[E/H]$ to $D_E'$.  We can describe $k*E$ as the free $k$-module
with basis $\{\bar{e} \mid e \in E \}$ and multiplication satisfying
$\overline {e_1} \overline{e_2} = \lambda (e_1, e_2) \overline
{e_1e_2}$ for some map $\lambda \colon E \times E \to k \setminus 0$.
Let $\{x_1, \dots , x_n \}\subset E$ be a set of coset representatives for $H$
and write $x_ix_j = h_{ij}x_{l(i,j)}$ where
$h_{ij} \in H$ and $l(i,j) \in \{1,
\dots, n\}$.  Then $D_H' *[E/H]$ is the ring which is a free
$D_H'$-module with basis $\{\widehat{x_1}, \dots \widehat{x_n} \}$, and
multiplication satisfying $\widehat{x_i} \widehat{x_j}
= \lambda (x_i, x_j) h_{ij} \widehat{x_{l(i,j)}}$.  Now define
$\theta \colon D_H'* [E/H] \to D_E'$ by
\[
\theta \bigl(\sum_{i=1}^n d_i \widehat{x_i} \bigr) =
\sum_{i=1}^n d_i \overline {x_i}.
\]
A routine check shows that $\theta$ is a ring homomorphism.  Also
it is easy to see that $\{ \overline{x_1}, \dots, \overline{x_n} \}$
is linearly independent over $D_H'$ and
it follows that $\theta$ is a monomorphism.

We deduce from the previous paragraph that $D_H'*[E/H]$
is a domain for all finite subgroups $E/H$ of $G/H$.  The result now
follows by Lemma \ref{Tvirabelian}.
\end{proof} 

\begin{remark}
  Natural constructions of a skew field extension $D_H$ usually are
  Hughes free.

  All examples of groups with finite classifying space
known to us with enough torsion-free quotients actually
  have enough solvable-by-finite torsion-free quotients.

  Consequently, the corresponding assumptions in Corollary
  \ref{cor:amext} are no serious restriction of generality.
\end{remark}

\begin{proof}[Proof of Theorem \ref{theorem:finiteext}\\]
\ref{item:main1}  For a prime number $p$ and a $p$-Sylow subgroup $S$ of $Q$,
let $G_S\subgroup G$ be its inverse image
   in $G$. By Lemma
  \ref{lem:psubgroup_check} it suffices to establish the statement 
  for every $KG_S$. To do this, we produce
  a subgroup $U$ of $H$, normal in $G_S$, such that $G_S/U$ is
  elementary amenable, $H/U$ is torsion-free and we have the
  divisibility relation
  $\lcm(G_S/U) \divides \lcm(G_S)$.

  Since $H$ has a finite classifying space, it is torsion-free. 
  Hence by \cite[Proposition 4 and Lemma 3]{Schick(1999)} it and
  its subgroup $U$ fulfill the strong Atiyah
  conjecture if and
  only if $DH$ and $DU$ are skew fields. In particular, we know that
  $DU$ is Artinian. We want to apply Proposition \ref{prop:amext} to
  the extension
  \begin{equation*}
    1\to U \to G_S \to G_S/U \to 1.
  \end{equation*}

  By Lemma \ref{lem:trivial_finite_ext_bound}, if $G_E/U\subgroup
  G_S/U$ is finite, then 
  $G_E\subgroup G_S$ fulfills
  the assumptions of Proposition \ref{prop:amext} with $L=\abs{G_E/U}$,
  and therefore also with $L=\lcm(G_S/U)$. Since $\lcm(G_S/U)$ divides
  $\lcm(G_S)$, by Corollary \ref{corol:atiyah_amext} the strong
  Atiyah conjecture holds for $KG_S$.

  It remains to find the subgroup $U$ with the required
  properties. Since $H$ is torsion-free and $G_S/H$ is a finite
  $p$-group, $\lcm(G_S)=p^n$ for some $n\ge 0$. Assume that
  every $G_S/U^{G_S}$ (defined in Definition \ref{def:Uto_the_H}) for
  $U\in\tfsubgroups_H$ contains a subgroup of order
  $p^k$. Note that each finite subgroup of each $G_S/U^{G_S}$ is a finite
  $p$-group, since $G_S/U^{G_S}$ contains the torsion-free subgroup
  $H/U^{G_S}$ with index a power of $p$. Therefore, it is sufficient to show that in this case $G_S$ itself contains a
  subgroup of order $p^k$. By Lemma
  \ref{lem:tfquotient_split_criterion} we can pass to subgroups $G_0$
  of $G_S$ and $Q_0$ of $G_S/H$ such that $\abs{Q_0}=p^k$ and the
  projection  $\pi$ in
  \begin{equation*}
    1\to {\hat H}^p \to {\hat G_0}^p\xrightarrow{\pi} Q_0\to 1
  \end{equation*}
  splits. Now our important technical Theorem \ref{theorem:splitback}
  applies and yields that $Q_0$ also splits back to $G_0$, hence gives a
  subgroup of order $p^k$ in $G_S$.

\ref{item:main2}
  Assume now that $G$ is torsion-free.
  The above argument implies that there is $U_S$ such that
  $G_S/U_S^{G_S}$ does not
  contain any non-trivial finite subgroup (the trivial group now is the biggest
  finite subgroup
  of $G_S\subgroup G$), i.e.~$G_S/U_S^{G_S}$ is torsion-free. We can
  choose $U_S \in\tfsubgroups_H'$ (since we can restrict ourselves to
  the use of $\tfsubgroups_H'$ throughout).
  
  Let $U$ be the intersection of all $U_S$ for all Sylow subgroups $S$
  of $Q$. Since $G_S/U^{G}$ fits into the exact sequence
  \begin{equation*}
    1\to  U_S^G/U^G \to  G_S/U^G \to G_S/U_S^G \to 1
  \end{equation*}
  with both $G_S/U_S^G$ and $U_S^G/U^G\subgroup H/U^G$ torsion-free
  (by Lemma \ref{lem:subgroupintersect}),
  $G_S/U^G$ is torsion-free.
  
  If $ G/U^G$ would contain a non-trivial torsion element, then by
  raising to an appropriate power we would have a $p$-torsion element
  $1\ne x\in G/U^G$ for some prime number $p$. Since $H/U^G$ is torsion-free, $x$
  would map to a $p$-torsion element in $G/H$, hence to a $p$-Sylow
  subgroup $S$ of $G/H$. Therefore, we
  could assume that $x$ is contained in $G_S/U^G$, which is
  torsion-free. This is a contradiction, hence $G/U^G$ itself is
  torsion-free. Since $H/U$ is elementary amenable,
  by Lemma
  \ref{lem:subgroupintersect} the same is true for $H/U^G$ and
  therefore also for its finite extensions $G/U^G$.
\end{proof}

\section{Extensions of pure braid groups,
  one-relator groups and link groups}
\label{sec:examples}

The next task is to find examples to which Theorem
\ref{theorem:finiteext} 
applies. In particular, we want to show that this is the case for the
pure braid groups and generalizations hereof.

To do this, we will prove that certain types of extensions preserve
the conditions of 
Theorem \ref{theorem:finiteext}, and that we can obtain e.g.~the braid
groups using these constructions.

\subsection{Extensions}

We have introduced a number of properties of groups, in particular, to 
have enough torsion-free quotients, to be residually torsion-free
nilpotent, and to be cohomologically complete. We now study when
these properties are preserved by extensions of groups. This will be
our main tool to produce interesting examples which have these properties. 

\begin{lemma}\label{lem:technical_extensions}
  Let $1\to H\to G\xrightarrow{p} Q\to 1$ be an exact sequence of
  groups. Assume that for each $m\in\naturals$ there is
  $n\in\naturals$ such that $\gamma_n(G)\cap H\subset \gamma_m(H)$.
  \begin{enumerate}
  \item \label{item:enough_quotients} Assume that $H$ has enough torsion-free nilpotent
    quotients. If $Q$ has enough
    torsion-free nilpotent quotients, or enough torsion-free solvable
    quotients, or enough torsion-free solvable-by-finite quotients,
    respectively, then the same is true for $G$.
  \item \label{item:residuality} Assume that $H$ is residually
    torsion-free nilpotent. If $Q$
    is residually torsion-free nilpotent, and for each
    $n\in\naturals$ there is a factorization $Q\onto Q/V\onto
    Q/\gamma_n(Q)$ with $Q/V$ torsion-free nilpotent, then $G$ is
    residually torsion-free nilpotent. If, in addition, for each
    $n\in\naturals$ there is a factorization $H\onto H/W\onto
    H/\gamma_n(H)$ with $H/W$ torsion-free nilpotent, then the
    corresponding statement is true for $G$.

  Corresponding results hold if in the assumptions on $Q$ ``nilpotent''
  is replaced by ``solvable'', or ``solvable-by-finite'',
  respectively, and the lower central series by the derived series or
  the appropriate series for solvable-by-finite quotients. Beware
    that, in the
    assumptions on $H$, we have to require nilpotent throughout.
  \end{enumerate}
\end{lemma}
\begin{proof}
  We first prove \ref{item:enough_quotients}. Let therefore $K\subgroup
  G$ be a normal subgroup with $G/K$ a finite $p$-group for some prime number 
  $p$. Then we get an exact sequence
  \begin{equation*}
    1\to H/(H\cap K)\to G/K\to
  Q/p(K)\to 1
\end{equation*}
  of finite $p$-groups. Let $V\subgroup H\cap K\subgroup
  H$ be a characteristic subgroup of $H$ such that $H/V$ is torsion-free
  nilpotent. Such a $V$ exists since $H$ has enough torsion-free
  nilpotent quotients. Since $H/V$ is nilpotent, there is an $m\in
  \naturals$ with $\gamma_m(H)\subgroup V$. Choose now $n\in\naturals$ 
  with $\gamma_n(G)\cap H\subgroup \gamma_m(H)\subgroup V$, which
  exists by assumption. Choose $n$ in such a way that
  $\gamma_n(G)\subgroup K$ (this is possible since $G/K$ is
  nilpotent). Observe that $p(\gamma_n(G)) = \gamma_n(Q)$.

  Since $Q$ has enough torsion-free nilpotent quotients, we find a
  characteristic subgroup $W$ of $Q$ such that $W\subgroup \gamma_n(Q)$,
  and $Q/W$ is torsion-free nilpotent. Assume that $\gamma_{n'}(Q)\subgroup 
  W$.

  Define now the subgroup $U$ of $G$ as $V\cdot (\gamma_n(G)\cap
  p^{-1}(W))$. Since $V\subgroup K$ and $\gamma_n(G)\subgroup K$, we
  get $U\subgroup K$. 

Since
  $V$ is a characteristic subgroup of $H$, it is normal in $G$. Since
  $p$ is surjective, $p^{-1}(W)$ and $\gamma_n(G)\cap p^{-1}(W)$ are
  also normal in $G$, which implies
  that $U$ is a normal subgroup of $G$.

Moreover,
\begin{equation*}
\gamma_{n'}(G)\subset \gamma_n(G)\cap p^{-1}(W)\subset V\cdot
(\gamma_n(G)\cap p^{-1}(W))=U,
\end{equation*}
i.e.~$G/U$ is nilpotent.

 The map $G\to Q/W$ is
  surjective and has kernel
  $p^{-1}(W)$, which contains $U$. On the other hand, since
  $W\subgroup \gamma_n(Q)$, 
  \begin{equation*}
 p^{-1}(W)\subgroup 
p^{-1}(\gamma_n(Q))=H\cdot\gamma_n(G),
\end{equation*}
i.e.~$p^{-1}(W)=H\cdot(\gamma_n(G)\cap p^{-1}(W))$.
It follows that the induced
  map $G/U\to Q/W$ is defined and has kernel 
  \begin{equation*}
H\cdot (\gamma_n(G)\cap p^{-1}(W))/(V\cdot (\gamma_n(G)\cap
  p^{-1}(W))) \iso H/(H\cap(V\cdot(\gamma_n(G)\cap p^{-1}(W)))).
\end{equation*}

  Now $V\subgroup H\subgroup p^{-1}(W)$, therefore
  \begin{equation*}
V\subgroup H\cap (V\cdot(\gamma_n(G)\cap p^{-1}(W)))= V\cdot(
\gamma_n(G)\cap H)\subgroup V\cdot
  \gamma_m(H)\subgroup V.
\end{equation*}
This means that we get the exact sequence
\begin{equation}\label{eq:sequence_for_G_mod_U}
  1\to H/V\to G/U\to Q/W \to 1.
\end{equation}
Since $H/V$ and $Q/W$ both are torsion-free
the same is true for
$G/U$. If $Q/W$ is only solvable, or solvable-by-finite, respectively, 
instead of nilpotent, \eqref{eq:sequence_for_G_mod_U} shows that the
same is true for $G/U$.

Now we prove \ref{item:residuality}. For each $1\ne g\in G$, we have to
produce a projection $p_g\colon G\to R_g$ with $R_g$ torsion-free
nilpotent and such that $p_g(g)\ne 1$. If $g$ is mapped to a 
non-trivial element of $Q$, this can be done because of the assumption 
that $Q$ is residually torsion-free nilpotent.

Assume therefore that $g\in H=\ker(p)$. Since $H$ is residually
torsion-free nilpotent, there is a (without loss of generality
characteristic) subgroup $V\subgroup H$ with $\gamma_m(H)\subgroup V$
and such that $H/V$ is torsion-free (and of course nilpotent), and such that
$p(g)\ne 1$ for the projection $p\colon H\to H/V$. Choose
now $n\in\naturals$ such that $\gamma_n(G)\cap H\subgroup
\gamma_m(H)\subgroup V$, and a normal subgroup $W\subgroup Q$ with
$W\subgroup \gamma_n(Q)$ and such that $Q/W$ is torsion-free and
elementary amenable (such a $W$ exists by assumption). As in the proof 
of \ref{item:enough_quotients}, define the subgroup
$U:=V\cdot(\gamma_n(G)\cap p^{-1}(W))$. We get the exact sequence \eqref{eq:sequence_for_G_mod_U}
\begin{equation*}
  1\to H/V\to G/U\to Q/W\to 1
\end{equation*}
and as before we conclude that $G/U$ is torsion-free nilpotent. If
$Q/W$ is only solvable or solvable-by-finite, at least these
properties are inherited by $G/U$.

Since $g$ maps to a non-trivial element in $H/V$, and the map $H/V\to
G/U$ induced by the inclusion is injective, $g$ maps to a non-trivial
element in $G/U$, as required.

For the last remaining statement, fix $k\in\naturals$. We assume now
that we find $W$ as above such that in addition $W\subgroup
\gamma_k(H)$. We construct $U$ as before. (Observe that $k\le m\le
n$), and it only remains to show that $U\subgroup\gamma_k(G)$. 

But $U\subgroup V\gamma_n(G)\subgroup\gamma_k(H)\gamma_k(G)=\gamma_k(G)$.
\end{proof}

We need two different types of assumptions in Lemma
\ref{lem:technical_extensions}: assumptions on the building blocks
(which we obviously have to make to get corresponding results for the
big group $G$) and assumptions on the extension itself. We now address
the question of when the latter ones are satisfied.

\begin{lemma}\label{lem:good_nilpotency_conditions}
  Let $1\to H\to G\to Q\to 1$ be an exact sequence of groups such that
  $G$ acts unipotently on $H_1(H;\integers)$ as defined in Definition \ref{def:unipotent_action}. Assume that at least one 
  of the following conditions is satisfied:
  \begin{enumerate}
  \item $Q$ is nilpotent,
  \item The extension $1\to H\to G\to Q\to 1$ is split, i.e.~$G$ is a
    semi-direct product $G=H\semiprod Q$,
  \end{enumerate}
   Then for each $m\in \naturals$ there is $n\in\naturals$ with
   $\gamma_n(G)\cap H\subset \gamma_m(H)$.
\end{lemma}
\begin{proof}
  In the course of the proof we will have to use subgroups of iterated 
  commutators. This can conveniently be dealt with using the
  bracket-less notation as in \cite{Hall(1958)}. Given two subgroups
  $H$ and $K$ of a group $G$ we set $\gamma H K:=[H,K]$. We abbreviate
  $\gamma (\gamma H K) L=: \gamma^2 HKL$. With this notation,
  $\gamma_n(G)=\gamma^{n-1}G^n$.  

  First observe that by definition $G$ acts unipotently on
  $H_1(H;\integers)=H/\gamma H^2$ if and only if there is
  $m\in\integers$ with
  \begin{equation}\label{eq:nilp_action}
\gamma^m
  HG^m\subset \gamma H^2=\gamma_2(H).
\end{equation}
This follows since the action on $H/\gamma H^2$ is induced by conjugation.
Since $\gamma^m HG^m\subset
  \gamma_{m+1}(G)\cap H$, it is therefore necessary that $G$ acts
  unipotently on $H_1(H;\integers)$ for the conclusion of the lemma.
  \cite[Lemma 7]{Hall(1958)} says that Equation \eqref{eq:nilp_action} 
  implies that 
  \begin{equation}
    \label{eq:nilp_conclusion}
    \gamma^{rm} H^r G^{rm-r+1} \subset \gamma^r H^{r+1}=\gamma_{r+1}(H).
  \end{equation}
  Using this inductively, we also see that Equation
  \eqref{eq:nilp_action} implies that for suitable $f(n)\in\naturals$
  \begin{equation}\label{eq:iterated_nilpotent}
    \gamma^{f(n)} H G^{f(n)} \subset \gamma^n H^{n+1}.
  \end{equation}

  Assume first that $Q$ is nilpotent. This means that there is
  $n\in\naturals$ with $\gamma_n(G)\subset H$ (since it is mapped to
  $\gamma_n(Q)=\{1\}$). Combined with \eqref{eq:iterated_nilpotent} this
  implies $\gamma_{n+f(m)}(G)\subset \gamma_{m+1}(H)$, which is what
  we had to prove.

Assume now that the sequence $1\to H\to G\to Q\to 1$ splits. With this 
assumption, we generalize a result of Falk and Randell
\cite[p.~85]{Falk-Randell(1985)} who prove that
$\gamma^n(G)\cap H=\gamma^n(H)$ for a semidirect product $G=H\semiprod 
Q$ if $Q$ acts \emph{trivially} on $H_1(H;\integers)$. 

Following Falk and Randell, let $p\colon G\to Q$ be the projection,
and $j\colon Q\to G$ a split of $p$ (i.e.~$p\circ j=\id_Q$). Using
$j$, we consider $Q$ as subgroup of $G$. Define
\begin{equation*}
  \tau\colon G\to H\colon g\mapsto j(p(g^{-1})) g =p(g^{-1})g.
\end{equation*}
Clearly $p(\tau(g))=1$ for each $g\in G$, i.e.~$\tau$ indeed maps $G$
to $H$. For each $g\in G$, $g=hq$ with $h\in H$ and $q\in Q$ if and
only if $h=\tau(g)$ and $q=p(g)$. In particular, if $g\in H$ then $g=\tau(g)$.

The map $\tau$ is not a homomorphism, but satisfies
\begin{equation}\label{eq:tau_equation}
  \tau(g_1g_2) = p(g_2)^{-1}\tau(g_1)p(g_2) \cdot \tau(g_2)
\end{equation}
(for all these statements, compare \cite[p.~83]{Falk-Randell(1985)}).

We now claim that for each $m\in\naturals$ there is $n\in\naturals$
with
\begin{equation}\label{eq:tau_inclusion}
 \tau(\gamma_{n}(G))\subset \gamma_m(H).
\end{equation}
This implies the
assertion, because $x\in \gamma_n(G)\cap H$ implies $x=\tau(x)\in
\gamma_m(H)$.

We are not proving \eqref{eq:tau_inclusion} directly, but rather prove 
by induction
that for each $m\in\naturals$ 
\begin{equation}
  \label{eq:tau_incl2}
  \tau(\gamma_{2^{m}}(G))\subset \gamma^m H G^m.
\end{equation}
In view of the assumption that $\gamma^m HG^m\subset \gamma H^2$ and
by \eqref{eq:iterated_nilpotent} this implies \eqref{eq:tau_inclusion}.

We need the following lemma.
\begin{lemma}\label{lem:help_interchange}
  If $H$ is a normal subgroup of $G$, then for each $m\in\naturals$
  \begin{equation*}
  \gamma^{m(m-1)/2+1} G^{m(m-1)/2+1} H  \subset    \gamma^m H G^m.
  \end{equation*}
\end{lemma}
\begin{proof}
  This is the statement of \cite[Theorem 2]{Hall(1958)} applied to
  $H/(\gamma^mHG^m)$ and $G/(\gamma^mHG^m)$ inside $G/(\gamma^m HG^m)$.
\end{proof}

The case $m=0$ of \eqref{eq:tau_incl2} is trivial.

For the induction step, assume $m\ge 1$ and choose $g_1\in
\gamma_{2^{m}-1}(G)\subset\gamma_{2^{m-1}}(G)$ and $g_2\in
G$, i.e.~$[g_1,g_2]=
g_1^{-1}g_2^{-1}g_1g_2$ is a 
typical element of $\gamma_{2^{m}}(G)$. Write $h_1=\tau(g_1)$, by
induction $h_1\in \gamma^{m-1}H G^{m-1}$. Set $h_2=\tau(g_2)$. Similarly, set
$q_1:=p(g_1)\in \gamma_{2^{m}-1}(Q)\subset \gamma_{2^{m}-1}(G)$,
$q_2:=p(g_2)$. We write $x^y=y^{-1}xy$. Then, following the
calculation of
\cite[p.~85]{Falk-Randell(1985)} 
\begin{equation*}
  \tau([g_1,g_2]) = [q_1,q_2]^{-1}
  [q_1,h_2]^{h_1}[q_1,q_2]^{h_2h_1}[h_1,h_2][h_1,q_2]^{h_2}. 
\end{equation*}
Using Lemma \ref{lem:help_interchange}, $[q_1,h_2]\in
\gamma^{2^m}G^{2^m}H \subset \gamma^{m}H 
G^{m}$, since $m(m-1)/2+1\le 2^m$ for each $m\in\naturals$, $m\ge 1$. By induction, $[h_1,h_2], [h_1,q_2]\in \gamma^m H
G^{m-1}G=\gamma^m H G^m$. Last, observe that $[[q_1,q_2],h_2h_1] \in
\gamma^{2^m} G^{2^m}H\subset \gamma^{{m}} H G^{{m}}$, so $[q_1,q_2]$ 
commutes with $h_2h_1$ modulo $\gamma^{{m}} H G^{{m}}$. Consequently,
\begin{equation*}
  \tau([g_1,g_2]) \equiv  [q_1,q_2]^{-1} [q_1,q_2] =1 \mod
  \gamma^{m} HG^{m}
\end{equation*}
and so $\tau(\gamma_{2^m }(G))\subset \gamma^m H G^m$.

This finishes the proof of Lemma \ref{lem:good_nilpotency_conditions}.
\end{proof}

\begin{remark}
  Given an extension $1\to H\to G\to Q\to 1$, one can ask whether the
  condition that $G$ acts unipotently on $H_1(H;\integers)$ is
  sufficient to imply that for each fixed $m$, if $n$ is sufficiently
  large then $\gamma_n(G)\cap H\subset
  \gamma_m(H)$. This is not the case in
  general. There are even counterexamples where $H$ is a central
  subgroup of $G$.
\end{remark}

\begin{lemma}\label{lem:iterated_nilpotent_action}
  Assume that $1=G_0\subgroup G_1\subgroup\cdots\subgroup G_n=G$ is a
  nested sequence of groups such that each $G_i$ is normal in
  $G$. Suppose that $G/G_i$ acts unipotently on
  $H_1(G_{i+1}/G_i;\integers)$ for each $i=0,\dots,n-1$. Then $G$
  acts unipotently on $H_1(G_{i};\integers)$ for all $i$.
\end{lemma}
\begin{proof}
  Induction reduces to the situation $1=G_0\subgroup G_1\subgroup
  G_2\subgroup G_3=G$, $G$ acts unipotently on $H_1(G_1;\integers)$
  and $G/G_1$ acts unipotently on
  $H_1(G_2/G_1;\integers)$, and we have to show that $G$ acts
  unipotently on $H_1(G_2;\integers)$. The exact sequence $1\to G_1\to 
  G_2\to G_2/G_1\to 1$ yields an exact sequence
  \begin{equation}\label{eq:H_1_exact_sequence}
    H_1(G_1;\integers)\to H_1(G_2;\integers)\to
    H_1(G_2/G_1;\integers)\to 0
  \end{equation}
  (observe that the first map is not injective in general). Since the
  group homomorphisms $G_1\to G_2\to G_2/G_1$ are compatible with the
  conjugation action by $G$, the same is true for the induced action
  in homology. Set $U:= H_1(G_1;\integers)$, $V:=H_1(G_2;\integers)$
  and $W:=H_1(G_2/G_1;\integers)$. We write $[U,G]\subgroup U$ for the 
  subgroup of $U$ generated by the elements $gu-u$ for $g\in G$, $u\in 
  U$ (using the $G$-action on $U$ induced by conjugation). As in the
  proof of Lemma \ref{lem:good_nilpotency_conditions}, we use the
  notation $\gamma^n UG^n=[\gamma^{n-1}UG^{n-1},G]$, with $\gamma
  UG=[U,G]$.

  Since $G$ acts unipotently on $W$, $\gamma^{n_2} WG^{n_2}=0$ for
  some $n_2\in\naturals$. That means that $\gamma^{n_2} VG^{n_2}$ is
  mapped to zero in the exact sequence \eqref{eq:H_1_exact_sequence},
  i.e.~lies in the image of $U$. Now $\gamma^{n_1} UG^{n_1}=0$ since
  $G$ acts unipotently on $U$ for suitable $n_1\in\naturals$. It
  follows that $\gamma^{n_1} (\gamma^{n_2} VG^{n_2})
  G^{n_1}=\gamma^{n_1+n_2} VG^{n_1+n_2} =0$. In other words, $G$ acts
  unipotently on $H_1(G_2;\integers)$.
\end{proof}

\begin{remark}
  Observe that, in the above lemma, even if $G/G_i$ acts
  trivially on $H_1(G_{i+1}/G_i;\integers)$ for each $i$, we can not 
  conclude that $G_i$ acts trivially on $H_1(G_{i-1};\integers)$, but
  only that it acts unipotently.
\end{remark}

\begin{lemma}\label{lem:nilpotent_Z_to_nilpotent_Zp}
  Assume that $N\subgroup G$ is a normal subgroup. Suppose $G$ acts
  unipotently on $H_1(N;\integers)$. Then $G$ acts unipotently on
  $H_1(N;\integers/p)$ for each prime number $p$.
\end{lemma}
\begin{proof}
  Since $H_0(N;\integers)=\integers$, by the universal coefficient
  theorem we have an isomorphism
  \begin{equation*}
    H_1(N;\integers/p) \iso H_1(N;\integers)\tensor_{\integers}\integers/p.
  \end{equation*}
  This isomorphism is natural and therefore in particular compatible
  with the $G$-action. That $G$ acts unipotently on $H_1(N;\integers)$ 
  means that there is an $n\in\naturals$ such that for $n$ arbitrary elements
  $g_1,\dots,g_n\in G$ the operator $(g_1-1)\cdots(g_n-1)=0$ acting on 
  $H_1(N;\integers)$. The universal coefficient theorem shows that the 
  same is true for the action on $H_1(N;\integers/p)$, i.e.~the action 
  of $G$ on $H_1(N;\integers/p)$ is nilpotent, as well.
\end{proof}

\subsection{The main class of groups}

\begin{definition}\label{def:extendableGroups}
  Let $\extendableGroups$ be the class of groups which fulfill the
  following list of properties:
  \begin{enumerate}
  \item $G$ has a finite classifying space $BG$.
  \item $G$ is cohomologically complete as defined in Definition
    \ref{def:cohom_complete}. 
  \item $G$ has enough nilpotent torsion-free quotients as in
    Definition \ref{def:tfquotients} (for example, all the quotients in the
    lower central series $G/\gamma_n(G)$ are
    torsion-free).
  \end{enumerate}
  Let $\extendableGroups_\infty$ the subclass of groups
  $G\in\extendableGroups$ which are, in addition, residually
  nilpotent.
\end{definition}

Observe that the definition of ``enough nilpotent torsion-free
quotients'' implies immediately the following.
\begin{remark}
  If $G\in\extendableGroups_\infty$ then $G$ is residually
  torsion-free nilpotent.
\end{remark}

\begin{lemma}\label{lem:basic_props_of_class}
    If $H\in\extendableGroups$, then $H$ is finitely generated, of finite
  cohomological dimension, torsion-free, and if $B$ is an
  $H$-module with finitely many elements, then $H^n(H,B)$ is finite
  for each $n$. Moreover, for
  each $k\in\naturals$ there is a normal subgroup $U_k$ of $H$ such
  that $H/U_k$ is torsion-free nilpotent and $U_k\subgroup \gamma_k(H)$.
\end{lemma}
\begin{proof}
    Finiteness of $BH$ immediately implies that $H$ is torsion-free and
  the finiteness assertions.
  The existence of $U_k$ with the required properties is proved in
  Lemma \ref{lem:tfprime_exists_primitive}.
\end{proof}

\begin{example}\label{ex:free_are_good}
  Every finitely generated free group $F_n$ belongs to $\extendableGroups_\infty$.
\end{example}
\begin{proof}
  Clearly $F_n$ has a finite classifying space. It is a classical result that
  $F_n$ is residually torsion-free nilpotent and that even
  $F_n/\gamma_k(F_n)$ is torsion-free for each $k$ (compare
  \cite{magnus37:_bezieh_kommut}). By example
  \ref{ex:free_completeness}, free groups are cohomologically
  complete. The assertion follows.
\end{proof}

\begin{proposition}\label{prop:extension_property}
  Let $1\to H\to G\to Q\to 1$ be a split exact sequence of groups
  (i.e.~$G=H\semiprod Q$)  and
  assume that $Q$ acts unipotently on $H_1(H;\integers)$.

  If $H$ and $Q$ both belong to $\extendableGroups$, then also
  $G\in\extendableGroups$. 
  If even $H,Q\in \extendableGroups_\infty$, then
  $G\in\extendableGroups_\infty$.
\end{proposition}
\begin{proof}
  Since $H$ and $Q$ have finite classifying spaces, the same is true
  for $G$ by Lemma \ref{lemma:fibrationtype}.

  By Lemma \ref{lem:nilpotent_Z_to_nilpotent_Zp}, $Q$ acts unipotently 
  on $H_1(H;\integers/p)$ for each prime number $p$.   Theorem
  \ref{theo:extensions_and_cohomological_completeness} therefore
  implies that $G$ is cohomologically complete, 
  the same being true for $H$ and $Q$ by assumption. 

  We assume that the sequence $1\to H\to G\to Q\to 1$ is split
  exact. By Lemma \ref{lem:good_nilpotency_conditions} we can apply
  Lemma \ref{lem:technical_extensions} \ref{item:enough_quotients} to conclude that $G$ has enough 
  torsion-free nilpotent quotients if $H$ and $Q$ both have the same
  property.
  In particular, if $H$ and $Q$ belong to $\extendableGroups$,
  the same is true for $G$.

  Because of Lemma \ref{lem:basic_props_of_class}, if $H$ and $Q$
  belong to $\extendableGroups_\infty$, i.e.~are residually
  torsion-free nilpotent, we can apply Lemma
  \ref{lem:technical_extensions} \ref{item:residuality} to conclude
  that $G$ is residually torsion-free nilpotent, as well, i.e.~belongs 
  also to $\extendableGroups_\infty$.
\end{proof}

\begin{remark}
  Using Lemma \ref{lem:technical_extensions} and Lemma
  \ref{lem:good_nilpotency_conditions}, other variants of the
  statements of Proposition \ref{prop:extension_property} can be
  proved. We leave this to the interested reader.
\end{remark}

The next theorem shows why the classes $\extendableGroups$ and
$\extendableGroups_\infty$ are important to us: if a group belongs to
them, then the Atiyah conjecture is inherited by extensions with
elementary amenable quotient.

\begin{theorem}\label{theorem:atiyah_for_extensions}
  Assume that $H\in \extendableGroups$ and $H$ fulfills the strong Atiyah
  conjecture over $KH$, where $K=\overline K$ is a subfield of
  $\complexs$. For example, assume that $H\in\extendableGroups_\infty$ 
  and $K=\overline\rationals$, where $\overline\rationals$ is the
  algebraic closure of $\rationals$ in $\complexs$. Let $G$ be an
  extension of $H$ with $\lcm(G)<\infty$ and such that $G/H$ is
  elementary amenable. Then $G$
  fulfills the strong Atiyah conjecture
  over $K G$.

  If the identity map of $KH$ extends to an isomorphism $KH_{\Sigma(H)}\to
  DH$, then the identity map of $K G$
  extends to an isomorphism $K G_{\Sigma(G)}\to DG$.
 
  Let $\tfsubgroups_H'$ be a nilpotent exhaustive set of subgroups of
  $H$. If, in addition, $G$ is torsion-free and $G/H$ is finite, then for
  $n\in\naturals$ sufficiently
  large we find $U_n\in\tfsubgroups_H'$ with $U_n\subgroup
  \gamma_n(H)$ such that $G/U_n$ is torsion-free (and of course virtually
  nilpotent). 

  If, again with $G$ torsion-free, we have a domain $k$ and a
  crossed product $k*G$ such that the subring $k*H$ embeds into a skew
  field $D_H$
  and the twisted action of $Q$ on $k*H$ extends to $D_H$, then $k*G$
  embeds into a skew field.
\end{theorem}
\begin{proof}
  If $G\in\extendableGroups_\infty$, then $\overline \rationals G$
  satisfies the strong Atiyah conjecture by Corollary
  \ref{corol:Atiyah_for_restfnilpotent}. 
  The other statements are immediate consequences of Lemma
  \ref{lem:basic_props_of_class}, Theorem \ref{theorem:finiteext}
  and Corollary \ref{cor:amext}. We use here that all groups in
  $\tfsubgroups_H'$ 
  are characteristic. In the first instance, Theorem \ref{theorem:finiteext}
   produces only one $U\in\tfsubgroups_H'$ with $G/U$
  torsion-free. Obviously $U=: U_{n_0}\subgroup \gamma_{n_0}(H)$
  for some $n_0\in\naturals$. For $n>n_0$, choose
  $V_n\in\tfsubgroups_H'$ with $V_n\subgroup \gamma_n(H)$. Set
  $U_n:=U\cap V_n\in\tfsubgroups_H'$. The exact sequence $1\to
  U/U_n\to G/U_n\to G/U\to 1$ shows that $G/U_n$ is torsion-free since
  $G/U$ and $U/U_n\subgroup H/U_n$ are torsion-free. By Lemma
  \ref{lem:subgroupintersect}, $G/U_n$ is virtually nilpotent.  
\end{proof}

The following example shows that the additional assumptions in Theorem
\ref{theorem:atiyah_for_extensions} often are satisfied.

\begin{example}\label{ex:poly_free_good}
  Suppose $H$ is poly-orderable, i.e.~admits a finite normal series
  with orderable factors. Assume that $G$ is a torsion-free
  extension of $H$ and $G/H$ is elementary 
  amenable.
If $k*G$ is a crossed product
  with a skew field $k$, then the crossed product $k*H$ embeds
  into a canonically defined skew
  field $D_G$, which comes with a (unique) twisted
  action of $G/H$ extending the twisted action on $k*H$.
Furthermore, if $H$ is poly-free,
the identity map of $KH$ extends to an
  isomorphism $KH_{\Sigma(H)}\to DH$. 

  Remember that orderable means two-sided
orderable, and that every residually torsion-free nilpotent group, in
particular every finitely generated free group, is
orderable \cite[Theorems 2.1.1 and 3.1.5]{Mura-Rhemtulla(1977)}.
We want to remark that all the
  groups in $\extendableGroups_\infty$
  are residually torsion-free nilpotent and
  therefore poly-orderable. Many, like the pure braid groups, are also 
  poly-free.
\end{example}
\begin{proof}
  The second statement follows by applying \cite[Lemma
  12.5]{Linnell(1998)} to each
  of the finitely many free extensions in the construction of $H$ as a
  poly-free group, using the fact that poly-free groups fulfill the
  strong Atiyah conjecture (we note that \cite[Lemma
12.5]{Linnell(1998)} remains true when $\mathbb {C}$ is replaced by
an arbitrary subfield $K$ of $\mathbb {C}$ which is closed under
complex conjugation).

For the first statement, we know from \cite[\S 1]{Hughes(1972)} that
$k*H$ has a free division ring of fractions in the sense of Hughes
\cite[\S 2]{Hughes(1970)}.  This means in particular that any
automorphism of $k*H$ which preserves the monomials of $k*H$ (i.e.\
sends an element of  the form $ah$ with $a \in k$ and $h \in H$ to
another such element) extends uniquely to an automorphism of $D_H$.
Thus we get the required (unique) twisted action of $G/H$ on $D_H$
extending the action on $k*H$.
\end{proof}

We continue with some abstract results which show how the methods we
have developed so far give information for certain groups beyond
the class $\extendableGroups$.

\begin{lemma}\label{lem:tf_solvable_root}
  Let $G_1,\dots,G_N$ be torsion-free solvable groups and
  \begin{equation*}
G:=G_1*\cdots* G_N
\end{equation*}
 their free product. Every projection $G\onto Q$
  with $Q$ solvable has a factorization $G\onto Q_1\onto Q$ where
  $Q_1$ is torsion-free solvable.
\end{lemma}
\begin{proof}
  By \cite[Lemma 5.5]{Gruenberg(1957)} we have an exact sequence
  \begin{equation*}
    1\to F\to G_1*\cdots* G_N\to G_1\times\cdots\times G_N\to 1,
  \end{equation*}
  where $F$ is a free group. Let $U$ be the kernel of the projection
  $G\onto Q$. Since $Q$ is solvable  $F/(F\cap U)$ also is
  solvable. Therefore, $F\cap U$ contains the derived series subgroup
  $F^{(k)}$ for a suitable $k\in\naturals$. Since every subgroup of a
  free group is free, the abelianization $F^{(k)}/F^{(k-1)}$ of every
  derived series subgroup of $F$ is free abelian and in particular torsion
  free. By induction, using the exact sequence $1\to
  F^{(k-1)}/F^{(k)}\to F/F^{(k)}\to F/F^{(k-1)} \to 1$, 
  $F/F^{(k)}$ is also torsion-free. Moreover, since $F^{(k)}$ is a
  characteristic subgroup of $F$, it is normal in $G$. Since $G_1,\dots,G_N$
  are solvable, for sufficiently large $k$ we get an extension
  \begin{equation*}
    1\to F/F^{(k)} \to G/F^{(k)} \to G_1\times\cdots\times G_N\to 1
  \end{equation*}
  of $G/F^{(k)}$ by torsion-free solvable groups, therefore
  $G/F^{(k)}$ itself is torsion-free solvable, and by the choice of
  $F^{(k)}$ it maps onto $G/U=Q$.
\end{proof}

\begin{proposition}\label{prop:free_product}
  Assume that $N\in\naturals$ and that
  $G_1,G_2,\dots,G_N\in\extendableGroups$.
  \begin{enumerate}
  \item \label{item:enough_solv_quotients}
    Then the free product $G:=G_1*G_2*\cdots* G_N$ has enough
    solvable torsion-free quotients.
    \item\label{item:free_product_in_F}
    If, more restrictively, $G$ has enough nilpotent torsion-free
    quotients, then $G$ belongs to $\extendableGroups$.
   \item \label{item:res_tf_solvable}
    If $G_i$ is residually torsion-free solvable for each $i$ (e.g.~if
    $G_i\in\extendableGroups_\infty$), then the same is true for $G$.
\end{enumerate}
\end{proposition}
\begin{proof}
  For \ref{item:free_product_in_F}, observe that the classifying space
  of $G$ is the one-point union of the
  classifying spaces of $G_j$ for $j=1,\dots,N$, therefore is compact. By
  Proposition \ref{prop:free_product_complete}, $G$ is cohomologically
  complete. If $G$ has enough nilpotent torsion-free quotients, then it
  follows that $G$ belongs to $\extendableGroups$.

 To prove \ref{item:res_tf_solvable}, remember that the property of
 being torsion-free solvable is a root
  property in the sense of \cite{Gruenberg(1957)}, and free groups are 
  residually torsion-free solvable, even residually torsion-free
  nilpotent (compare \cite{magnus37:_bezieh_kommut}). Therefore, by \cite[Theorem
  4.1]{Gruenberg(1957)}, $G$ is residually torsion-free solvable if
  the same is true for all the $G_i$.

To establish \ref{item:enough_solv_quotients}, given a
normal subgroup $U$ of $G$, if we set 
  $U_j:=G_j\cap U$, then $G/U$ is a quotient of $(G_1/U_1)*\cdots *
  (G_N/U_N)$. If the index of $U$ in $G$ is a power of $p$, the same
  is true for the index of $U_j$ in $G_j$ for
  $j=1,\dots,N$. Be definition of the class $\extendableGroups$, we
  find $V_j\subgroup U_j$ normal in 
  $G_j$ such that $G_j/V_j$ is torsion-free nilpotent and hence
  torsion-free solvable. By Lemma
  \ref{lem:tf_solvable_root} applied to the solvable (since it is a
  finite $p$-group) quotient $G/U$ of $(G_1/V_1)*\cdots * (G_N/V_N)$
  we find a factorization $(G_1/V_1)*\cdots *(G_N/V_N)\to Q_1\to G/U$
  where $Q_1$ is torsion-free solvable. It follows from Lemma
  \ref{lem:tfprime_exists} that $G$ has enough solvable-by-finite torsion-free
  quotients. 
\end{proof}

\begin{proposition}\label{prop:generalext}
  Let $G_1,\dots,G_N$ be groups which belong to
  $\extendableGroups_\infty$. Assume that $Q$ has a finite
  classifying space,
  is cohomologically complete, is residually torsion-free
  solvable-by-finite and
  has enough solvable-by-finite torsion-free quotients. Define
  \begin{equation*}
    G:= G_1\semiprod\bigl(G_2\semiprod(\cdots (G_N\semiprod Q))\bigr),
  \end{equation*}
  and assume that in each semidirect product the quotient acts
  unipotently on the first homology of the kernel. Then $G$ has a
  finite classifying space, is residually torsion-free
  solvable-by-finite, has enough
  torsion-free solvable-by-finite quotients, and is cohomologically
  complete. If
  there is an exact
  sequence 
  \begin{equation*}
    1\to G\to H\to A\to 1
  \end{equation*}
  with $A$ elementary amenable and $\lcm(H)<\infty$, then $\overline\rationals
  H$ fulfills the strong Atiyah conjecture.
\end{proposition}
\begin{proof}
   By an iterative
  application of Lemma \ref{lemma:fibrationtype}, $G$ has a finite
  classifying space, and by Lemma \ref{lem:good_nilpotency_conditions} 
  and Lemma \ref{lem:technical_extensions}, $G$ has enough
  solvable-by-finite 
  torsion-free
  quotients and is residually torsion-free solvable-by-finite (observe 
  that having enough torsion-free solvable-by-finite quotients implies 
  the factorization property for the lower central series quotients as 
  required in Lemma \ref{lem:technical_extensions}
  \ref{item:residuality} because of Lemma
  \ref{lem:tfprime_exists_primitive}). Theorem
  \ref{theo:extensions_and_cohomological_completeness}
 implies that $G$ is
  cohomologically complete. The last assertion follows from Corollary
  \ref{cor:amext} and Corollary \ref{corol:Atiyah_for_restfnilpotent}.
\end{proof}

\begin{example}
  In Proposition \ref{prop:generalext}, $Q=Q_1*\cdots * Q_m$ can
  be a free product of finitely many groups
  $Q_i\in\extendableGroups_\infty$. By Proposition
  \ref{prop:free_product}, $Q$ has
  enough solvable torsion-free quotients and is residually
  torsion-free solvable, and by Proposition
  \ref{prop:free_product_complete} $Q$ is cohomologically
  complete. Its classifying space is the one-point union of the
  $BQ_i$, $i=1,\dots,m$, and therefore is finite.
\end{example}

In the next subsections we describe certain interesting classes of
groups which belong to $\extendableGroups$ or even to
$\extendableGroups_\infty$. 

\subsection{One-relator groups}

\begin{definition}
  A one-relator group $G$ is called \emph{primitive} if it is finitely
  generated and if it has a presentation $G=\generate{x_1,\dots x_d \mid
  r}$ such that the element $r$ in the free group $F$ generated by
  $x_1,\dots,x_d$ is contained in $\gamma_n(F)$ but not in
  $\gamma_{n+1}(F)$ and the image of $r$ in
  $\gamma_n(F)/\gamma_{n+1}(F)$ is not a proper power.

  We adopt the convention that finitely generated free groups are
  primitive one-relator groups.
\end{definition}

\begin{example}\label{ex:primitive_one_relator}
  The following finitely generated one-relator groups are primitive
  \begin{itemize}
  \item one-relator groups where the least common multiple of the
    exponent-sums for the different generators in the relator is one.
  \item fundamental groups of compact orientable two-dimensional
    surfaces (with or without boundary).
  \end{itemize}
\end{example}

\begin{proposition}\label{prop:primitive_one_relator}
  Let $G$ be a primitive one-relator group. Then $G$ belongs to
  $\extendableGroups$. 
\end{proposition}
\begin{proof}
  By \cite[Exemple (2), p.~144]{Labute(1967)} such a group is
  cohomologically complete. The main theorem of \cite{Labute(1970)}
  says in particular
  that $G/\gamma_n(G)$ is torsion-free for each $n\in\naturals$. The
  fact that $G$ is primitive implies that the relator $r$ is not a
  proper power. By \cite{Lyndon(1950)} $G$ has cohomological dimension
  $2$, more precisely, the presentation complex of
  $\generate{x_1,\dots,x_d\mid r}$ with $d$ $1$-cells and one $2$-cell
  is a model for $BG$, i.e.~$BG$ is finite.
\end{proof}

\begin{remark}\label{rem:primitive_one_relator}
  In light of Theorem \ref{theorem:atiyah_for_extensions} and
  Proposition \ref{prop:primitive_one_relator} it is of course
  important to know which one-relator groups satisfy the Atiyah
  conjecture. Examples for this are residually torsion-free solvable
  one-relator groups.
  Unfortunately, we don't know whether all one-relator groups, or at
  least all torsion-free one-relator groups, satisfy the Atiyah conjecture.
\end{remark}

\subsection{Link groups}

\begin{definition} We denote the fundamental group $G$ of the
  complement $M$ of a tubular neighborhood $\nu(L)$ of a tame link $L$ with $d$
  components in
  $S^3$  a \emph{link group with $d$ components}.
  We define the  \emph{linking diagram} to be the edge-labeled graph whose
  vertices
  are the components of the link, and such that any pair of vertices
  is joined by exactly one edge. Each edge is labeled with the linking
  number of the two link
  components involved. 

  We say the link group $G$ is \emph{primitive} if for each prime number $p$
  there is a spanning subtree of the linking diagram such that none of
  the labels of the edges of this subtree is congruent to $0$ modulo
  $p$.

  Observe that in particular every \emph{knot group} (i.e.~a link group 
  with only one component) is primitive.
\end{definition}

\begin{definition}
  A compact orientable $3$-manifold (possibly with boundary) is
  called \emph{irreducible} if every embedded two-sphere bounds an
  embedded three-disc. It is called \emph{Haken} if it is irreducible
  and if there is a properly embedded two-sided surface different from
  a two-sphere such that the induced map of fundamental groups is injective.
\end{definition}

\begin{proposition}\label{prop:primite_links}
  Assume that $G$ is a primitive link group with $d$ components. Then
  $G$ 
has a finite classifying space and $G/\gamma_n(G)$ is
  torsion-free for each $n\in\naturals$. Moreover,
   $G\in
  \extendableGroups$.
\end{proposition}
\begin{proof}
  Let $M\subset S^3$ be the link-complement such that $G=\pi_1(M)$. Then $M$ is
  an orientable manifold with non-empty boundary, and each boundary
  component is a torus.  Moreover, primitivity in particular implies
  that the link is
  non-split, i.e.~if we have a two-sphere embedded into the complement
  of the link in $S^3$, then by Alexander's theorem the link will
  intersect only one of the
  two discs into which this two-sphere splits $S^3$
  (else the linking graph would contain two non-empty parts, connected
  only by edges labeled with zero). This also shows that each two-sphere embedded into the
  complement of the link bounds a three-disc embedded into the
  complement (the one of the above three-discs which does not intersect
  the link). By definition, this means that $M$
  is irreducible.

  Since the boundary of $M$ contains a surface which is not a
  two-disc, $M$ is a Haken manifold by
  \cite[Lemmas
  6.6 and 6.7]{Hempel(1987)}. By \cite[Remark III.19]{Jaco(1980)} the
  finite CW-complex $M$ is the classifying space of $G$ (this is also
  directly a consequence of \cite[Theorem 27.1]{Papakyriakopoulos(1957)}). By
  \cite[Theorems 2 and 1]{Labute(1990)} $G/\gamma_n(G)$ is
  torsion-free for each $n\in\mathbb {N}$.

  In \cite{Linnell-Schick(2000c)} we show that $G$ is
  cohomologically $p$-complete for every prime number
  $p$.  This finishes
  the proof.
\end{proof}

The following proposition is a special case of Proposition
\ref{prop:primite_links}. We single it out here because we can give an
independent and more elementary proof of the key property, namely
cohomological completeness.

\begin{proposition}\label{prop:knots}
  Let $G$ be a knot group, i.e.~the fundamental group of the complement
  of a tame knot in $S^3$, or let $G$ be a primitive link group with
  two components. Then $G\in\extendableGroups$.
\end{proposition}
\begin{proof}
  Above, we have given a proof of the fact that $G$ has a finite
  classifying space, namely the knot or link complement $S^3-\nu(L)$. Let
  $d\in\{1,2\}$ be the number of components of $L$.

  Let $h\colon G\to G^{ab}=H_1(G,\integers)$ be the projection of $G$ onto
  its abelianization. 

  We claim that $G^{ab}$ is isomorphic to $\integers^d$ and that $h$
  induces an isomorphism in homology and
  cohomology (with  integral coefficients).
  The assertion then follows from Proposition \ref{prop:homology_iso} and
  Example \ref{ex:free_ab_cohom}.

  It remains to prove the claim.

  A classical computation in knot theory, using
  Alexander duality, shows that $H_n(G,\integers)$ and  $H^n(G,\integers)$ are
  free abelian for each $n$, and
    \begin{equation*}
\rank
  H^n(G,\integers)=
  \begin{cases}
    1; & \text{if }n=0\\ d; & \text{if } n=1\\ d-1; & \text{if }n=2 \\
    0; & \text{if }n>2.
  \end{cases}
\end{equation*}

Hence for trivial reasons (since $G^{ab}$ is free abelian) $H^n(h)$ is an
isomorphism except possibly for
$n=2$ and $d=2$. In this case $G^{ab}=\integers^2$ and the generator of
$H^2(\integers^2,\integers)$ is the cup product of the two
generators of $H^1(\integers^2,\integers)$. Alexander duality and
the definition of the linking number implies that the cup product of
the two generators of $H^1(G,\integers)$ is the linking number times
the generator of $H^2(G,\integers)$. By assumption, this linking
number is $\pm 1$, so that $G\to G^{ab}$ induces also an isomorphism
on $H^2$. Using the universal coefficient theorem to see that $H_n(h)$
also is an isomorphism for all $n$, this concludes the proof.
\end{proof}

\begin{corollary}
  Let $H$ be either a primitive one-relator group (compare
  Remark \ref{rem:primitive_one_relator})
or a primitive link group, and assume that $H$ is
  residually torsion-free solvable.

  If there is an exact sequence $1\to H\to G\to A\to 1$ with $A$
  elementary amenable and $\lcm(G)<\infty$, then the strong Atiyah
  conjecture is true for $\overline\rationals G$.
\end{corollary}
\begin{proof}
  By Proposition \ref{prop:primitive_one_relator} or by Proposition
  \ref{prop:primite_links}, $H\in\extendableGroups$.  It therefore
   follows in
  both cases by Corollary \ref{corol:Atiyah_for_restfnilpotent} that
  $\overline\rationals H$ satisfies the strong Atiyah conjecture. 
\end{proof}

\subsection{Homology-invariance}

\begin{proposition}\label{prop:homology_iso}
  Assume that $G\in\extendableGroups$, and that $H_1,H_2$ are groups with
  finite classifying spaces and $f_1 \colon H_1\to G$,
$f_2 \colon G\to H_2$ are
  group homomorphisms which induce isomorphisms on homology with
  integral coefficients. Then $H_1$ and $H_2$ both belong to
  $\extendableGroups$, and $f_1$ as well as $f_2$ induce an
  isomorphism between the pro-$p$ completions (for every prime number $p$).
\end{proposition}
\begin{proof}
  By \cite[Theorem 3.4]{Stallings(1965b)}, $f_1$ or $f_2$
  induce isomorphisms between the lower central series quotients
  $G/\gamma_n(G)$ and $H_1/\gamma_n(H_1)$ or $H_2/\gamma_n(H_2)$,
  respectively, for
  every $n\in\naturals$. This implies that the maps induce
  isomorphisms between the inverse systems of nilpotent quotients of
  $G$, $H_1$ and $H_2$. Hence, by Lemma \ref{lem:tfprime_exists},
  since $G$ has enough torsion-free nilpotent quotients, the same is
  true for $H_1$ and for $H_2$.

   Moreover, in particular the
  finite $p$-group quotients of $G$, $H_1$ and $H_2$ coincide (since
  they are nilpotent). This means we
  get the following commuting diagram
  \begin{equation*}
    \begin{CD}
      H_1 @>{f_1}>> G @>{f_2}>> H_2\\
      @VVV @VVV @VVV\\
      {\hat H_1}^p @>{{\hat f_1}^p}>> {\hat G}^p @>{{\hat f_2}^p}>>
      {\hat H_2}^p
    \end{CD}
  \end{equation*}
  and ${\hat f_1}^p$ as well as ${\hat f_2}^p$ are
  isomorphisms. Taking cohomology with $\integers/p$-coefficients, we
  obtain
  \begin{equation*}
    \begin{CD}
      H^k(H_1,\integers/p) @<{H^k(f_1)}<< H^k(G,\integers/p)
      @<{H^k(f_2)}<< H^k(H_2,\integers/p)\\ 
      @A{\alpha_1}AA @A{\beta}AA @A{\alpha_2}AA\\
      H^k({\hat H_1}^p,\integers/p) @<{H^k({\hat f_1}^p)}<< H^k({\hat
      G}^p,\integers/p) @<{H^k({\hat f_2}^p)}<< 
      H^k({\hat H_2}^p,\integers/p)
    \end{CD}
  \end{equation*}
  where now not only $H^k({\hat f_1}^p)$ and $H^k({\hat f_2}^p)$, but
  by assumption, using the universal coefficient theorem for cohomology,
  also $\beta$, $H^k(f_1)$ and $H^k(f_2)$ are isomorphisms. It follows that
  $\alpha_1$ and $\alpha_2$ are isomorphisms, too. This implies that
  $H_1$ and $H_2$ are cohomologically complete and concludes the proof.
\end{proof}

\begin{definition}
  A discrete group $G$ is called \emph{acyclic} if
  $H_k(G,\integers)=0$ for $k>0$.
\end{definition}

\begin{example}\label{ex:free_ab_cohom}
  It is well known that finitely generated free abelian groups belong
  to $\extendableGroups$. Cohomological completeness, which is the
  only non-obvious point, follows from Proposition
  \ref{prop:extension_property}.

  Let $G$ be a group with
  finite classifying space and such that 
  $H_1(G,\integers)$ is free abelian. Assume that the projection $G\to
  G/[G,G]=H_1(G,\integers)$ is a homology isomorphism,
  i.e.~$H_k(G,\integers)\to H_k(G/[G,G],\integers)$ is an isomorphism
  for every $k\in\naturals_0$. It follows from
  Proposition \ref{prop:homology_iso} that
  $G\in\extendableGroups$. Observe that this follows in particular if
  $G$ is acyclic or if $H_*(G,\integers)\iso H_*(\integers,\integers)$
  (as abstract abelian groups).

  We will see in Proposition \ref{prop:knots} that knot groups and
  certain link groups provide examples where the projection to the
  abelianization is a homology isomorphism.

  Acyclic groups with finite classifying space are not hard to come
  by. They play an important role in asphericalization processes:
  given any finite simplicial complex $X$, we can use any acyclic group
  with finite classifying space
  as building block to construct a new finite simplicial complex $BLX$ which
  has the same cohomology as $X$ and is aspherical (i.e.~which is the
  classifying
  space of its fundamental group $LX\colon=\pi_1(BLX)$). For all this, compare
  \cite{Baumslag-Dyer-Heller(1980)}. 

  Fix an acyclic group which has a finite classifying space.
 Following \cite[Section
  9]{Baumslag-Dyer-Heller(1980)}, this gives rise to an
  asphericalization procedure, a functor
  $BL$ from finite simplicial complexes to aspherical finite
  simplicial complexes (more precisely, to $X$ we assign a map $BLX\to
  X$ which induces a surjection on $\pi_1$ and an isomorphism in
  integral homology).

  If we take any acyclic finite simplicial complex $X$ and apply the
  asphericalization construction, then we obtain a new acyclic group
  $LX=\pi_1(BLX)$.

  As a preparation to get non-acyclic examples, let $X$ and $Y$ be two
  spaces whose integral cohomology ring in both cases is isomorphic to
  the cohomology ring of $\integers^n$, i.e.~is a symmetric algebra
  over $\integers$ on $n$ generators which belong to $H^1$. Let
  $f \colon Y\to X$ be a map which induces an isomorphism
  $H^1(f):H^1(X,\integers)\to H^1(Y,\integers)$. Then $f$ is a
  homology isomorphism, because an algebra map like $H^*(f)$ which
  is a bijection between the  generators of two symmetric algebras is
  an isomorphism, and then we can apply the universal coefficient
  theorem.

  Assume next that we apply the asphericalization procedure described above
  to a finite simplicial complex $X$ whose 
  integral cohomology ring is isomorphic to the cohomology ring of the
  free abelian group
  $\integers^n$. We
  get a homology isomorphism $BLX\to X$, in particular
  $H_1(BLX,\integers)=\integers^n$ and $H^*(BLX,\integers)$ is the
  cohomology ring of
  $\integers^n$. The projection $\pi \colon LX\to
  H_1(LX,\integers)=LX/[LX,LX]=\integers^n$ induces by duality an
  isomorphism on the first cohomology with integral
  coefficients. Therefore, it
  is a homology isomorphism. By Proposition
  \ref{prop:homology_iso}, $LX$ belongs to $\extendableGroups$.

  Note that the conditions on the cohomology ring structure are
  trivially fulfilled if $H_*(X,\integers)\iso
  H_*(\integers,\integers)$ (as abelian groups).
\end{example}

\subsection{Braid groups}





\begin{theorem}\label{theo:fiber_type_arrangements}
    Assume that $H$ is the
  fundamental groups of a (complement of a) fiber-type
  arrangement. Then, $H$ belong to $\extendableGroups_\infty$. In particular,
  this is the case for
  Artin's classical pure braid groups and, more generally, for the
  generalized pure braid groups associated to the Coxeter groups
  $A_l$, $C_l$, $G_2$ and $I_2(p)$ \cite[Theorem
  3.1]{Falk-Randell(1988)}.
\end{theorem}
\begin{proof}
    Falk and Randell construct \cite[Proposition
  2.5]{Falk-Randell(1985)} $G$ as
  iterated semidirect product with trivial action on $H_1$ of the
  kernel, starting
  with a finitely generated free group. By Example
  \ref{ex:free_are_good} and Proposition
  \ref{prop:extension_property}, $G\in \extendableGroups_\infty$
  (it is already proved in \cite[Theorem
  2.6]{Falk-Randell(1988)} that $G$ is residually
  torsion-free nilpotent).
\end{proof}

\begin{corollary}\label{corol:Atiyah_for_braid}
  Assume that $H$ is one of the groups of Theorem
  \ref{theo:fiber_type_arrangements}, e.g.~a classical Artin pure
  braid group $P_n$. Then
   every finite extension $G$ of $H$ fulfills the strong Atiyah
   conjecture over $\overline{\rationals}$. Moreover, if such an extension is torsion-free then it is residually
  (torsion-free virtually nilpotent), i.e.~it has plenty  of
  non-trivial torsion-free quotients.

   This applies in particular to Artin's classical
   full braid groups $B_n$
  (as well as the full braid groups associated to $C_n$, $G_2$ and
  $I_2(p)$).

  Specifically, for each 
  classical braid group $B_n$ ($n\in\naturals$ arbitrary) the quotients
  $B_n/\gamma_N(P_n)$ are torsion-free for all $N\ge N(n)$
  sufficiently large, and the
  intersection of the $\gamma_N(P_n)$ is trivial. This implies also that
  the commutator subgroup $B_n'$ has non-trivial torsion-free quotients.

  In view of Example \ref{ex:poly_free_good}, the additional statements
  about localization and crossed products in Theorem
  \ref{theorem:atiyah_for_extensions} apply. In particular, if we have
  an exact sequence  $1\to H\to G\to A\to 1$ where $A$ is elementary
  amenable, and $\lcm(G)<\infty$, then the strong Atiyah conjecture is 
  true for $\overline\rationals G$. If $G$ is torsion-free
  and $k*G$ is a crossed product with a skew field $k$, then $k*G$
  embeds into a skew field.
\end{corollary}
\begin{proof}
  By Theorem 
  \ref{theo:fiber_type_arrangements}, this follows from Theorem
  \ref{theorem:atiyah_for_extensions}. 
\end{proof}

\begin{remark}
  It was a question of Lin \cite[0.9 e,f, Remark 7.3]{Lin(1996)}
  whether the braid groups $B_n$ have non-trivial non-abelian
  torsion-free quotients (actually, he conjectures this not to be the
  case for
  $n>4$). This question is also popularized as question B6 on the
  group theory server \href{http://www.grouptheory.info}{``www.grouptheory.info''}. The question is
  answered with yes for $n\le 6$ in
  \cite{Humphries(1999)} with an explicit torsion-free
  quotient. Humphries' computations are rather involved. In the course
  of them he makes use of computer algebra programs. He has explicit candidates
  for torsion-free quotients of $B_n$ for all $n$. Humphries
  proves 
  that the kernel $K_n$ of his projection contains $\gamma_n(P_n)$
  and, for $n\le
  6$ is contained in $P_n$. Explicit calculations of the ranks of the
  nilpotent abelian quotients of $P_n/K_n$ show that in general
  $K_n\ne \gamma_k(P_n)$
  for any $k\in\naturals$.

  As a by-product of our approach, we get a torsion-free
  quotient as asked for by Lin for each $B_n$ (the only drawback is
  that we only give an existence proof, but without further
work we can't specify the integers $N$ for which $G_n/\gamma_N(P_n)$
is torsion-free).
\end{remark}


\bibliographystyle{plain}
\bibliography{extAtiyah}

\begin{thebibliography}{10}

\bibitem{Adem-Milgram(1994)}
Alejandro Adem and R.~James Milgram.
\newblock {\em Cohomology of finite groups}.
\newblock Springer-Verlag, Berlin, 1994.

\bibitem{Atiyah(1976)}
M.~F. Atiyah.
\newblock Elliptic operators, discrete groups and von {N}eumann algebras.
\newblock In {\em Colloque ``Analyse et Topologie'' en l'Honneur de Henri
  Cartan (Orsay, 1974)}, pages 43--72. Ast\'erisque, No. 32--33. Soc. Math.
  France, Paris, 1976.

\bibitem{Baumslag-Dyer-Heller(1980)}
G.~Baumslag, E.~Dyer, and A.~Heller.
\newblock The topology of discrete groups.
\newblock {\em J. Pure Appl. Algebra}, 16(1):1--47, 1980.

\bibitem{Mura-Rhemtulla(1977)}
Roberta Botto~Mura and Akbar Rhemtulla.
\newblock {\em Orderable groups}.
\newblock Marcel Dekker Inc., New York, 1977.
\newblock Lecture Notes in Pure and Applied Mathematics, Vol. 27.

\bibitem{Bourbaki(1989)}
Nicolas Bourbaki.
\newblock {\em Commutative algebra. {C}hapters 1--7}.
\newblock Springer-Verlag, Berlin, 1989.
\newblock Translated from the French, Reprint of the 1972 edition.

\bibitem{Bratzler(1997)}
Clemens Bratzler.
\newblock {\em $L^2$-Betti Zahlen und Faserungen}.
\newblock Diplomarbeit, Universit{\"a}t Mainz, 1997.
\newblock
  \href{http://wwwmath.uni-muenster.de/u/lueck/group/bratzler.dvi}{http://wwwm%
ath.uni-muenster.de/u/lueck/group/bratzler.dvi}.

\bibitem{Brown(1982)}
Kenneth~S. Brown.
\newblock {\em Cohomology of groups}.
\newblock Springer-Verlag, New York, 1982.

\bibitem{Cohn(1985)}
P.~M. Cohn.
\newblock {\em Free rings and their relations}.
\newblock Academic Press Inc. [Harcourt Brace Jovanovich Publishers], London,
  second edition, 1985.

\bibitem{Dicks-Schick(2002)}
Warren Dicks and Thomas Schick.
\newblock The spectral measure of certain elements of the complex group ring of
  a wreath product.
\newblock {\em Geom. Dedicata}, 93:121--137, 2002.

\bibitem{Dodziuk(1977)}
Jozef Dodziuk.
\newblock De {R}ham-{H}odge theory for ${L}\sp{2}$-cohomology of infinite
  coverings.
\newblock {\em Topology}, 16(2):157--165, 1977.

\bibitem{Dodziuk-Linnell-Mathai-Schick-Yates(2001)}
Jozef D\'odziuk, Peter Linnell, Varghese Mathai, Thomas Schick, and Stuart
  Yates.
\newblock Approximating ${L}^2$-invariants, and the {A}tiyah conjecture.
\newblock {\em Comm. Pure Appl. Math.}, 56(7):839--873, 2003.

\bibitem{Dummit-Foote(1999)}
David~S. Dummit and Richard~M. Foote.
\newblock {\em Abstract algebra}.
\newblock Prentice Hall Inc., Englewood Cliffs, NJ, 2nd. edition, 1994.

\bibitem{Falk-Randell(1985)}
Michael Falk and Richard Randell.
\newblock The lower central series of a fiber-type arrangement.
\newblock {\em Invent. Math.}, 82(1):77--88, 1985.

\bibitem{Falk-Randell(1988)}
Michael Falk and Richard Randell.
\newblock Pure braid groups and products of free groups.
\newblock In {\em Braids (Santa Cruz, CA, 1986)}, pages 217--228. Amer. Math.
  Soc., Providence, RI, 1988.

\bibitem{Fenn-Greene(2000)}
R.~Fenn, M.~T. Greene, D.~Rolfsen, C.~Rourke, and B.~Wiest.
\newblock Ordering the braid groups.
\newblock {\em Pacific J. Math.}, 191(1):49--74, 1999.

\bibitem{Grigorchuk-Linnell-Schick-Zuk(2000)}
Rostislav~I. Grigorchuk, Peter Linnell, Thomas Schick, and Andrzej {\.Z}uk.
\newblock On a question of {A}tiyah.
\newblock {\em C. R. Acad. Sci. Paris S\'er. I Math.}, 331(9):663--668, 2000.

\bibitem{Gromov(1991)}
Michael Gromov.
\newblock K\"ahler hyperbolicity and ${L}\sb 2$-{H}odge theory.
\newblock {\em J. Differential Geom.}, 33(1):263--292, 1991.

\bibitem{Gruenberg(1957)}
K.~W. Gruenberg.
\newblock Residual properties of infinite soluble groups.
\newblock {\em Proc. London Math. Soc. (3)}, 7:29--62, 1957.

\bibitem{Hall(1958)}
P.~Hall.
\newblock Some sufficient conditions for a group to be nilpotent.
\newblock {\em Illinois J. Math.}, 2:787--801, 1958.

\bibitem{Hempel(1987)}
John Hempel.
\newblock Residual finiteness for $3$-manifolds.
\newblock In {\em Combinatorial group theory and topology (Alta, Utah, 1984)},
  pages 379--396. Princeton Univ. Press, Princeton, NJ, 1987.

\bibitem{Hirsch(1938)}
K.~A. Hirsch.
\newblock On infinite soluble groups. {I}{I}.
\newblock {\em Proc. London Math. Soc. (2)}, 44:336--344, 1938.

\bibitem{Hughes(1970)}
Ian Hughes.
\newblock Division rings of fractions for group rings.
\newblock {\em Comm. Pure Appl. Math.}, 23:181--188, 1970.

\bibitem{Hughes(1972)}
Ian Hughes.
\newblock Division rings of fractions for group rings. {I}{I}.
\newblock {\em Comm. Pure Appl. Math.}, 25:127--131, 1972.

\bibitem{Humphries(1999)}
Stephen~P. Humphries.
\newblock Torsion-free quotients of braid groups.
\newblock {\em Internat. J. Algebra Comput.}, 11(3):363--373, 2001.

\bibitem{Jackowski(1988)}
Stefan Jackowski.
\newblock A fixed-point theorem for $p$-group actions.
\newblock {\em Proc. Amer. Math. Soc.}, 102(1):205--208, 1988.

\bibitem{Jaco(1980)}
William Jaco.
\newblock {\em Lectures on three-manifold topology}.
\newblock American Mathematical Society, Providence, R.I., 1980.

\bibitem{Jacobson(1989)}
Nathan Jacobson.
\newblock {\em Basic algebra. {I}{I}}.
\newblock W. H. Freeman and Company, New York, second edition, 1989.

\bibitem{Kochman(1996)}
S.~O. Kochman.
\newblock {\em Bordism, stable homotopy and {A}dams spectral sequences}.
\newblock American Mathematical Society, Providence, RI, 1996.

\bibitem{Kropholler-Linnell-Moody(1988)}
P.~H. Kropholler, P.~A. Linnell, and J.~A. Moody.
\newblock Applications of a new ${K}$-theoretic theorem to soluble group rings.
\newblock {\em Proc. Amer. Math. Soc.}, 104(3):675--684, 1988.

\bibitem{Linnell-Schick(2000c)}
Inga K{\"u}mpel, Peter Linnell, and Thomas Schick.
\newblock Galois cohomology of completed link groups.
\newblock in preparation.

\bibitem{Labute(1967)}
John~P. Labute.
\newblock Alg\`ebres de {L}ie et pro-$p$-groupes d\'efinis par une seule
  relation.
\newblock {\em Invent. Math.}, 4:142--158, 1967.

\bibitem{Labute(1970)}
John~P. Labute.
\newblock On the descending central series of groups with a single defining
  relation.
\newblock {\em J. Algebra}, 14:16--23, 1970.

\bibitem{Labute(1990)}
John~P. Labute.
\newblock The {L}ie algebra associated to the lower central series of a link
  group and {M}urasugi's conjecture.
\newblock {\em Proc. Amer. Math. Soc.}, 109(4):951--956, 1990.

\bibitem{Lin(1996)}
V.~Lin.
\newblock Braids, permutations, polynomials {I}.
\newblock Preprint, Max Planck Institut f{\"u}r Mathematik, Bonn, 112 pages,
  1996.

\bibitem{Linnell(1993)}
Peter~A. Linnell.
\newblock Division rings and group von {N}eumann algebras.
\newblock {\em Forum Math.}, 5(6):561--576, 1993.

\bibitem{Linnell(1998)}
Peter~A. Linnell.
\newblock Analytic versions of the zero divisor conjecture.
\newblock In {\em Geometry and cohomology in group theory (Durham, 1994)},
  pages 209--248. Cambridge Univ. Press, Cambridge, 1998.

\bibitem{Lueck(1989)}
Wolfgang L{\"u}ck.
\newblock {\em Transformation groups and algebraic ${K}$-theory}, volume 1408
  of {\em Lecture Notes in Mathematics}.
\newblock Springer-Verlag, Berlin, 1989.
\newblock Mathematica Gottingensis.

\bibitem{Lueck(1997)}
Wolfgang L{\"u}ck.
\newblock {$L\sp 2$}-invariants of regular coverings of compact manifolds and
  {CW}-complexes.
\newblock In {\em Handbook of geometric topology}, pages 735--817.
  North-Holland, Amsterdam, 2002.

\bibitem{Lueck-Schick-Thielmann(1998)}
Wolfgang L{\"u}ck, Thomas Schick, and Thomas Thielmann.
\newblock Torsion and fibrations.
\newblock {\em J. Reine Angew. Math.}, 498:1--33, 1998.

\bibitem{Lyndon(1950)}
R.~C. Lyndon.
\newblock Two notes on nilpotent groups.
\newblock {\em Proc. Amer. Math. Soc.}, 3:579--583, 1952.

\bibitem{magnus37:_bezieh_kommut}
Wilhelm Magnus.
\newblock {\"U}ber beziehungen zwischen h{\"o}heren {K}ommutatoren.
\newblock {\em J. Reine Angew. Math.}, 177:105--115, 1937.

\bibitem{Papakyriakopoulos(1957)}
C.~D. Papakyriakopoulos.
\newblock On {D}ehn's lemma and the asphericity of knots.
\newblock {\em Ann. of Math. (2)}, 66:1--26, 1957.

\bibitem{Passman(1989)}
Donald~S. Passman.
\newblock {\em Infinite crossed products}.
\newblock Academic Press Inc., Boston, MA, 1989.

\bibitem{Reich(1999)}
Holger Reich.
\newblock {\em Group von Neumann algebras and related algebras}.
\newblock PhD thesis, Universit{\"a}t G\"ottingen, 1999.
\newblock
  \href{http://www.math.uni-muenster.de/u/lueck/publ/diplome/reich.dvi}{http:/%
/www.math.uni-muenster.de/u/lueck/publ/diplome/reich.dvi}.

\bibitem{Rolfsen-Zhu(1998)}
Dale Rolfsen and Jun Zhu.
\newblock Braids, orderings and zero divisors.
\newblock {\em J. Knot Theory Ramifications}, 7(6):837--841, 1998.

\bibitem{Schick(2000a)}
Thomas Schick.
\newblock Finite group extensions and the {B}aum-{C}onnes conjecture.
\newblock preprint, available via
  \href{http://arXiv/math.KT/0209165}{{http://arXiv/math.KT/0209165}}.

\bibitem{Schick(1999)}
Thomas Schick.
\newblock Integrality of ${L}\sp 2$-{B}etti numbers.
\newblock {\em Math. Ann.}, 317(4):727--750, 2000.

\bibitem{MR1894160}
Thomas Schick.
\newblock Erratum: ``{I}ntegrality of {$L\sp 2$}-{B}etti numbers''.
\newblock {\em Math. Ann.}, 322(2):421--422, 2002.

\bibitem{Serre(1994)}
Jean-Pierre Serre.
\newblock {\em Cohomologie {G}aloisienne}.
\newblock Springer-Verlag, Berlin, fifth edition, 1994.

\bibitem{Serre(1997)}
Jean-Pierre Serre.
\newblock {\em Galois cohomology}.
\newblock Springer-Verlag, Berlin, 1997.
\newblock Translated from the French by Patrick Ion and revised by the author.

\bibitem{Stallings(1965b)}
John Stallings.
\newblock Homology and central series of groups.
\newblock {\em J. Algebra}, 2:170--181, 1965.

\bibitem{Switzer(1975)}
Robert~M. Switzer.
\newblock {\em Algebraic topology---homotopy and homology}.
\newblock Springer-Verlag, New York, 1975.
\newblock Die Grundlehren der mathematischen Wissenschaften, Band 212.

\bibitem{Wilson(1998)}
John~S. Wilson.
\newblock {\em Profinite groups}.
\newblock The Clarendon Press Oxford University Press, New York, 1998.

\end{thebibliography}

\end{document}